\magnification=\magstep1 
\overfullrule=0pt
\def\eqde{\,{:=}} \def\leaderfill{\leaders\hbox to 1em{\hss.\hss}\hfill} 
 \def\la{\lambda}  \def\ga{\gamma} 
       \def\i{{\rm i}} 
\def\si{\sigma} \def\eps{\epsilon}    \def\ka{{\kappa}}
         
  \def\L{{\Lambda}} \def\M{{\cal M}}
     
   \def\H{{\cal H}}
\font\huge=cmr10 scaled \magstep2

\font\smcap=cmcsc10     \font\smit=cmmi7  \font\smal=cmr7
\def\z{{\cal Z}}  \def\G{\Gamma}
   
\input amssym.def
\def\Z{{\Bbb Z}} \def\R{{\Bbb R}} \def\Q{{\Bbb Q}}  
\def\C{{\Bbb C}}   \def\M{{\Bbb M}}
\font\smit=cmmi7
\def\sdprod{{\times\!\vrule height5pt depth0pt width0.4pt\,}}

{\nopagenumbers
\rightline{August, 1999}
\bigskip\bigskip
\centerline{{\bf \huge  Monstrous Moonshine}}\medskip
\centerline{{\bf \huge and the Classification of CFT}} \bigskip
\centerline{{\bf (16 lectures given in Istanbul, August 1998)}}
\bigskip \bigskip   \centerline{Terry Gannon}\medskip
\centerline{{\it Department of Mathematical Sciences, University of Alberta,}}
\centerline{{\it Edmonton, Alberta, Canada, T6G 2G1} } \smallskip
\centerline{{e-mail: tgannon@math.ualberta.ca}}
\bigskip\medskip

\centerline{\bf Abstract}\medskip
\noindent
In these notes we give an introduction both to Monstrous Moonshine and to the
classification of rational conformal field theories, using this
as an excuse to explore several related structures and go on a little tour
of modern math. We will discuss Lie algebras, modular functions, the
finite simple group classification, vertex operator algebras, Fermat's Last
Theorem, category theory, (generalised) Kac-Moody
algebras, denominator identities, the A-D-E meta-pattern, representations
of affine algebras, Galois theory, etc. This work is informal and pedagogical,
 and aimed mostly at grad students in math or math phys, but I hope that
many interested nonexperts will find something of value here --- like
any good Walt Disney movie I try not to completely ignore the `grown-ups'.
My emphasis is on ideas and motivations, so these notes are intended to
complement other papers and books where this material is presented with
more technical detail.
The level of difficulty varies significantly from topic to topic. The two
parts --- in fact any of the sections --- can be read independently.
\vfill\eject

\centerline{{\bf Table of Contents}}\bigskip

\centerline{Part 1. The classification of conformal field theories}\medskip

\line{1.1. Informal motivation\leaderfill 1}
\line{1.2. Lie algebras\leaderfill 4}
\line{1.3. Representations of finite-dimensional simple Lie algebras\leaderfill
10}
\line{1.4. Affine algebras and the Kac-Peterson matrices\leaderfill 13}
\line{1.5. The classification of physical invariants\leaderfill 17}
\line{1.6. The A-D-E meta-pattern\leaderfill 19}
\line{1.7. Simple-currents and charge-conjugation\leaderfill 24}
\line{1.8. Galois theory\leaderfill 28}
\line{1.9. The modern approach to classifying physical invariants\leaderfill
31}\bigskip

\centerline{Part 2. Monstrous Moonshine}\medskip

\line{2.1. Introduction\leaderfill 35}
\line{2.2. Ingredient \#1: Finite simple groups and the Monster\leaderfill
38}
\line{2.3. Ingredient \#2: Modular functions and Hauptmoduls\leaderfill 40}
\line{2.4. The Monstrous Moonshine conjectures\leaderfill 46}
\line{2.5. Formal power series\leaderfill 48}
\line{2.6. Ingredient \#3: Vertex operator algebras\leaderfill 51}
\line{2.7. Ingredient \#4: Generalised Kac-Moody algebras\leaderfill 54}
\line{2.8. Ingredient \#5: Denominator identities\leaderfill 59}
\line{2.9. Proof of the Moonshine conjectures\leaderfill 60}

\bigskip
\line{References\leaderfill 64}
\vfill\eject
\centerline{{\bf Glossary}}\medskip

{\settabs\+ empty space&Perron-Frobenius theory\qquad\qquad\qquad&\cr
\+ &24& \S1.6\cr
\+ &$A_\ell$& sl$_{\ell+1}(\C)$\cr
\+ &affine algebra& \S1.4\cr
\+ &algebra& \S1.2\cr
\+ &category theory& \S1.1\cr
\+ &central extension& \S\S1.2, 1.4\cr
\+ &CFT&conformal field theory\cr
\+ &character& \S\S1.3, 1.4\cr
\+ &chiral algebra & \S1.1\cr
\+ &Fermat's Last Theorem& \S2.3\cr
\+ &field & \S1.2\cr
\+ &finite simple group& \S2.2\cr
\+ &genus& \S\S1.1, 2.3\cr
\+ &group& \S\S1.3, 2.2\cr
\+ &Hauptmodul& \S2.3\cr
\+ &highest-weight & \S\S1.3, 1.4\cr
\+ &Kac-Moody algebra& \S2.7\cr
\+ &lattice& \S1.6\cr
\+ &Leech lattice $\L_{24}$& \S2.4\cr
\+ &Lie algebra & \S1.2\cr
\+ &Lie group & \S1.2\cr
\+ &${\Bbb M}$& the Monster finite simple group\cr
\+ &manifold& \S1.2\cr
\+ &modular function& \S2.3\cr
\+ &modular group& PSL$_2(\Z)$ or SL$_2(\Z)$\cr
\+ &$P_+^k$&\S1.4\cr
\+ &partition function&\S\S1.1, 1.2, 1.5\cr
\+ &Perron-Frobenius& \S1.7\cr
\+ &physical invariant& RCFT partition function, \S1.5\cr
\+ &QFT &quantum field theory\cr
\+ &RCFT &rational conformal field theory\cr
\+ &representation& \S1.3\cr
\+ &Riemann surface& \S\S1.1, 2.3\cr
\+ &$S,\,T$& Kac-Peterson matrices, (1.4.2)\cr
\+ &subfactor theory& \S1.9\cr
\+ &topological& \S\S1.1, 2.3\cr
\+ &torus& \S2.3\cr
\+ &$V^\natural$& the Moonshine module, \S\S2.4, 2.6\cr
\+ &Verlinde's formula& (1.4.4)\cr
\+ &Virasoro algebra& (1.2.7)\cr
\+ &VOA&vertex operator algebra\cr
\+ &Weyl group& \S1.3\cr
\+ &WZW model& \S1.1\cr }\vfill\eject}\pageno=1

\centerline{{\huge Part 1. The classification of conformal
field theories}}\bigskip

\noindent{{\smcap 1.1. Informal motivation}}\bigskip

In this section we will sketch a very informal and `hand-wavy' motivation to what we shall call
the classification problem for rational conformal
field theory (RCFT). Much of this material is more carefully treated in
e.g.\ [18].

A CFT is a quantum field theory (QFT), usually with a two-dimensional space-time,
whose symmetries include the conformal transformations. There are  different approaches to
CFT --- for one of these see [26,27]. Another
formulation which has been deeply influential is due to Graeme Segal [52].
It is motivated by string theory and is phrased in an important mathematical language
called {\it category theory}.

A {\it category} consists of two types of things. One are called {\it objects},
and the other are called {\it arrows} (or {\it morphisms}). An arrow, written
$f:A\rightarrow B$, has an initial and a final
object ($A$ and $B$ respectively). Arrows $f,g$ can be composed to yield
a new arrow $f\circ g$, if the final object of $g$ equals the initial object
of $f$. Maps between categories are called {\it functors} if they
take the objects (resp., arrows) of one to the objects (resp., arrows)
of the other, and preserve composition.

The only difficulty people can have in understanding categories is in
realising that there is no real content to them. It's just a language,
highly abstract like the more familiar set theory, but in many
contexts (a great example is the theory of knot invariants [58])
one which  is both natural and
suggestive. It tries to deflect some of our instinctive infatuation
with objects (nouns), to the mathematically more fruitful one with
structure-preserving maps between objects (verbs). A gentle introduction to
the mathematics and philosophy of categories is [43]; we'll give a taste
of this shortly.

The standard example of a category is called {\bf Set}, where the `objects' are sets, and
the arrows from $A$ to $B$ are functions $A\rightarrow B$. A related example
that Segal uses is {\bf Vect}, where the objects are complex vector spaces and
the arrows are linear maps. A rather trivial example of a functor
${\cal F}:${\bf Vect}$\rightarrow${\bf Set} sends a vector space $V$
to its underlying set, also called $V$ --- i.e.\ ${\cal F}$ simply
`forgets' the vector space structure on $V$ and ignores the fact that
the arrows $f$ in {\bf Vect} are linear. The other category Segal uses
he  calls 
{\bf C}; its objects are disjoint unions of (parametrised) circles $S^1$, and
the arrows are (conformal equivalence classes of) {\it cobordisms}, i.e.\ (Riemann)
surfaces whose boundaries are those
circles. Composition of arrows in ${\bf C}$ is accomplished by gluing the surfaces along
the appropriate boundary circles.

Consider the usual definition of a one-to-one function: $f(x)=f(y)$ only
when $x=y$. Category theory replaces this with the following. The
arrow  $f:A\rightarrow B$
is called `monic' if for any arrow $g:C\rightarrow B$, there exists a unique
arrow $h:C\rightarrow A$ such that $f\circ h=g$. So it's a sort of factorisation
property. You can easily verify that in {\bf Set} the notions of `one-to-one'
and `monic' coincide. What does this redefinition gain us? It certainly
doesn't seem any simpler. But it does change the focus from the {\it
argument}  of $f$, to the {\it global} functional behaviour of $f$,
and a change of perspective can never
be bad. And it allows us to transport the idea of one-to-one-ness to
arbitrary categories. For instance, in the Riemann surface category {\bf C},
the `one-to-one functions' are the genus-0 cobordisms.

Or consider the notion of {\it product}. In category theory, we say that the
triple $(P,a,b)$ is a product of objects $A,B$ if $a:P\rightarrow A$ and
$b:P\rightarrow B$ are arrows, and if for any $f:C\rightarrow A$, $g:C\rightarrow
B$, there is a unique arrow $h:C\rightarrow P$ such that $f=a\circ h$ and
$g=b\circ h$. This notion unifies several constructions (each of which is
the `product' in an appropriately chosen category): Cartesian product of
sets; intersection of sets; multiplication of numbers; the logical operator
`and'; direct product; infimum in a partially ordered set; etc. {\it Sum}
can be defined similarly, unifying the constructions of disjoint union,
`or', addition, tensor product, direct sum, supremum, etc.

This generality of course comes with a price: it can wash away all of the
endearing special features of a favourite theory or structure. There
certainly  are
contexts where e.g.\ all human beings should be thought of as equal, but
there are other contexts where the given human is none other than
your mother and must be treated as such. It turns out that category theory
provides a beautiful framework for understanding topological invariants
such as the Jones-Reshetikhin-Turaev-Witten knot invariants (see e.g.\
[58]). And it seems
to be a natural language for formulating CFT axiomatically, as we'll now
see.

According to Segal, a CFT is a functor ${\cal T}$ from {\bf C} to {\bf Vect}, which
obeys various properties. The picture comes from string theory: the fundamental
object is a `vibrating'  loop; a state is given by  a collection
of these loops; each classical Feynman path from the initial
to the final states is a world-sheet, i.e.\ a surface $\Sigma$
whose boundary is all the loops in the initial and final states. QFT assigns a complex number $I$ (the action) to each of these
world-sheets, and the quantum amplitude, written  $\langle {\rm final}|{\rm initial}
\rangle$, will then be the integral over all worldsheets of $e^{\i I/\hbar }$.
(The quantum amplitude is how the theory makes contact with experiment, as it tells us the
probability of the given process $|{\rm initial}\rangle\mapsto|{\rm final}\rangle$
happening.) This is what
Segal is trying to capture formally. The vector spaces in {\bf Vect}
come in in order to handle uniformly and simultaneously 
the various `vibrations' of the strings. In particular there is one
basic vector space $H$ (a Hilbert space of quantum states), and the functor
will take $n$
copies of $S^1$ to $H^n=H\otimes\cdots \otimes H$.

The simplest interesting example here is the `tree-level creation of a string
from the vacuum'. In this case the world-sheet looks like a bowl,
i.e.\ {\it topologically} is a disk $D$ (if we imagine the bowl to be
made of rubber, we could grab its rim and stretch it down flat onto
the table, so we say a bowl and a disk are topologically equivalent
--- see also \S2.3). Segal's
functor gives us a linear map ${\cal T}(D):\C\rightarrow H$ ($H^0$
is just $\C$), which we can think of equivalently as the assignment
of the vector ${\cal T}(D)(1)$ in $H$ to $D$. In the case of the standard
unit disk (i.e.\ where the parametrisation of the boundary $S^1$ is simply
$\theta\mapsto e^{2\pi\i\theta}$), this vector is cal0led the {\it vacuum
state}  $\Omega=|0\rangle$.

For another example, consider the `vacuum-to-vacuum expectation value'. The
initial
and final states (objects) here are both the empty set, so the world-sheets
(arrows) are closed Riemann surfaces. As usual in QFT, we can organise these by
how many internal `loops' are involved (this number is called the {\it genus} of the
surface): topologically, 0-loop (i.e.\ `tree-level') world-sheets are spheres,
1-loop world-sheets are tori, etc. These closed Riemann surfaces are
discussed in more detail in \S 2.3. The 0-loop contribution isn't very interesting
(there is only one conformal equivalence class of spheres),
so let us look at the 1-loop contribution. It will be of the form
$\int {\cal Z}([torus])\, d[torus]$, where ${\cal Z}$ is a complex-valued
function called the partition function, and $[torus]$ is a conformal
equivalence class of tori. In the Segal formalism we recover this in the
following way: the functor takes
$[torus]$ to a linear function from $H^0=\C$ to $H^0=\C$. 
Any such linear function is simply a $1\times 1$
matrix, i.e.\ a complex number, which we call ${\cal Z}([torus])$.

Now there is a nice parametrisation of conformal equivalence classes of
tori, as we will see more explicitly in \S 2.3. Namely, a representative
for each class can be chosen to be of the form $\C/(\Z+\Z\tau)$ where
$\tau$ is in the upper half plane ${\cal H}$. Thus we can write ${\cal Z}$
as a function of a complex variable $\tau$. However, different $\tau$
correspond to the same equivalence class of tori: the redundancy is
exactly captured by the modular group PSL$_2(\Z)$. Namely, $\tau$ and
${a\tau +b\over c\tau+d}$ are equivalent, whenever $a,b,c,d\in\Z$ and
$ad-bc=1$. Thus ${\cal Z}(\tau)={\cal Z}({a\tau+b\over c\tau+d})$.
In other words, the partition function is modular 
invariant!\footnote{$^1$}{{\smal In  higher-dimensional string
theories, a  similar argument shows more
generally that {\smit automorphic}\ {\smit forms}\ will appear naturally.}}

There are two sectors in CFT, a holomorphic one and an antiholomorphic
one, corresponding to the two directions (`left-' and `right-moving')
of motion on a string, or the two components of the group Diff$(S^1)$ of
diffeomorphisms of the circle. This
means that many of the quantities (e.g.\ the partition function) factorise
into parts depending holomorphically and anti-holomorphically on the modular
parameters (e.g.\ ${\tau}$ in genus 1).
In a {\it rational} CFT  there are finitely many `primary fields' $a\in\Phi$ --- the precise meaning
of this is not important here, but it says that the space of states
for the theory decomposes into a finite sum\footnote{$^2$}{{\smal It
seems though that `rational' logarithmic CFT is trying to teach us the
lesson that this familiar requirement can and should be weakened. See
Gaberdiel-Kausch (1999).}}
$H=\oplus_{a,b\in\Phi}M_{ab}H_a\otimes
\overline{H}_b$, where $M_{ab}$ are nonnegative integers which count the
multiplicity of $H_a\otimes\overline{H}_b$ in $H$. The linear maps ${\cal T}
(\Sigma):H^m\rightarrow H^n$ in an RCFT will  factorise similarly;
this `chiral factorisation' is captured by what Segal calls the `modular
functor' [52]. The partition function becomes
$${\cal Z}(\tau)=\sum_{a,b\in\Phi}M_{ab}\,\chi_a(\tau)\,\chi_b(\tau)^*\eqno(1.1.1)$$
for certain holomorphic functions $\chi_a$. One of the primary fields (we'll
denote it
`0') corresponds to the vacuum $\Omega$, and uniqueness of the vacuum
means  that $M_{00}=1$.

$H_0$ is called a {\it chiral algebra}; in the language of \S 2.6,
$H_0$ will be a vertex operator algebra (VOA). $\Phi$
parametrises the irreducible $H_0$-modules and the $\chi$'s are their
characters; in an RCFT we require this
number to be finite. For example, for the Moonshine VOA $V^{\natural}$
discussed in Part 2, $\Phi$ consists of only one element.

The higher-genus behaviour of an RCFT is determined from the
lower-genus behaviour, by composition of `arrows' (i.e.\ the gluing
together of surfaces) in {\bf C}. See Figure 3 of [30] to
find how a genus-2 surface is built up from genus-0 ones. In fact, 
it's generally believed that an RCFT will be uniquely determined by: (i) the
choice of chiral algebra; (ii) the partition function (which tells you the
spectrum of the theory, i.e.\ how the two sectors link up); and (iii) the
structure constants $C_{ab}^c$, which in the Segal formalism correspond
to the surfaces called `pairs-of-pants', equivalently  disks with two interior
disks removed.
Our approach will be to start with a chiral algebra, and find all possible
partition functions. We will thus ignore the important question of
existence and uniqueness of the structure constants, though at least for our
chiral algebras, it
seems to be generally believed that the structure constants will be unique.

Perhaps all chiral algebras come from standard constructions (e.g.\ orbifolds
and the Goddard-Kent-Olive (GKO) coset construction --- see e.g.\
[18])  involving lattices
and affine Kac-Moody algebras. For instance a $\Z_2$-orbifold of the
VOA of the Leech lattice gives us the Moonshine module $V^{\natural}$,
and the so-called minimal models come from GKO cosets involving $A^{(1)}_1$.
This is in line with the spirit of Tannaka-Krein duality (and its
generalisations by Deligne and Doplicher-Roberts), which roughly says
that if a bunch of things {\it act like} they're the set of representations
of a Lie group, then they {\it are} the set of representations of a
Lie group.

In any case, one of the simplest, best understood, and important classes
(called Wess-Zumino-Witten (WZW) models --- see for instance [30,59]
in this volume) of RCFTs 
correspond to affine Kac-Moody algebras at a positive integer level
$k$. We will have much more to say later about these
algebras, but for now let us remark that $\Phi$ here will be the (finite)
set $P_+^k$ of integrable level $k$ highest weights $\la$.
Their chiral algebras were constructed by Frenkel and Zhu. The following
sections concern the attempt to classify the partition functions
corresponding to Kac-Moody algebras --- see especially \S 1.5. 
I will use this theme as an excuse to describe many other things, e.g.\
the A-D-E meta-pattern, Lie theory, Galois, fusion rings, ...
I dedicate these notes to 
the profound friendship developing in recent years between mathematics and
physics. As Victor Kac
said in his 1996 Wigner medal acceptance speech, ``Some of the
best ideas come to my field from the physicists. And on top of this they
award me a medal. One couldn't hope for a better deal.''

\bigskip\bigskip\noindent{{\smcap 1.2. Lie algebras}}\bigskip

{\it Lie algebras} (and their nonlinear partners {\it Lie
groups}) appear in numerous places throughout math and mathematical physics.
A nice introduction is [9]; Lie theory is presented with more of a physics
flavour in [24], as well as [59].

An {\it algebra} is a vector space with a way to multiply vectors which
is compatible with the vector space structure (i.e.\ the vector-valued product
 is {\it bilinear}: $(a\vec{u}+a'\vec{u}')\,(b\vec{v}+b'\vec{v}')=ab\,\vec{u}
 \vec{v}+ab'\,\vec{u}\vec{v}'+a'b\,\vec{u}'\vec{v}+a'b'\,\vec{u}'\vec{v}'$).
 For example, the complex numbers $\C$ can be thought of as a 2-dimensional
 algebra over $\R$ (a basis is 1 and $\i=\sqrt{-1}$; the {\it scalars}
here are real numbers and the {\it vectors} are complex numbers). The
quaternions
 are 4-dimensional over $\R$ and the octonions are 8-dimensional over $\R$.
 Incidentally, these are the only finite-dimensional algebras over $\R$ which
 obey the cancellation law: $\vec{u}\ne 0$ and $\vec{u}\vec{v}=0$ implies
 $\vec{v}=0$ (the reader should try to convince himself why the familiar vector
product on $\R^3$ fails the cancellation law). This important
 little fact makes several unexpected appearances in math. For instance,
 it is trivially possible to `comb the hair' on the circle $S^1$ without `cheating'
 (i.e.\ needing a hair-part or exploiting a bald spot): just comb the hair
 clockwise for example. However it is not possible to comb the hair on
 the sphere $S^2$ (e.g.\ your own head) without cheating. The only other
 $k$-spheres $S^k$ which can be combed (i.e.\ for which there exist $k$
 linearly independent continuous vector fields) are $k=3$ and 7. This
 is intimately connected with the existence of $\C$, the quaternions,
 and octonions (namely, $S^1,S^3,S^7$ can be thought of as the length 1
 complex numbers, quaternions, and octonions, resp.).

In a {\it Lie} algebra ${\frak g}$,  the
product is usually called a `bracket' and is written $[xy]$.
It is required to be `anti-commutative' 
and `anti-associative':
$$\eqalignno{[xy]+[yx]=&\,0&(1.2.1a)\cr
[x[yz]]+[y[zx]]+[z[xy]]=&\,0&(1.2.1b)\cr}$$
(like most other equalities in math, (1.2.1b) is usually
called the {\it Jacobi identity}). Usually
we will consider Lie algebras over $\C$, but sometimes over $\R$. Note
that (1.2.1a) says $[xx]=0$.

One important consequence of bilinearity is that it is enough to know
the values of all the brackets $[x_ix_j]$ for $i<j$, for any basis $\{x_1,x_2,\ldots\}$
of ${\frak g}$. (The reader should convince himself of this before proceeding.)

The simplest example of a Lie algebra is ${\frak g}=\C$ (or ${\frak g}=\R$), with the bracket
$[xy]$ identically 0. In fact, this is the only 1-dimensional Lie algebra.
It is a straightforward exercise for the reader to find all 2- and 3-dimensional
Lie algebras (over $\C$) up to isomorphism (i.e.\ change of basis): there are precisely
2 and 6 of them, respectively (though one of the 6 depends on a complex
parameter). Over $\R$, there are 2 and 9 (with 2 depending on real parameters).
This exercise cannot be continued much further
--- e.g.\ not all 7-dimensional Lie algebras (over $\C$) are known. Nor is it obvious
that this would be an interesting or valuable exercise. We should suspect
that our definition of Lie algebra is probably a little too general for anything
obeying it to be automatically an interesting structure. More often than
not, a  classification turns out to be a stale and useless list.

Two of the 3-dimensional Lie algebras are important in what follows. One
of them is well-known to the reader: consider the vector-product (also
called cross-product) in $\C^3$. Taking the standard basis
$\{e_1,e_2,e_3\}$ of $\C^3$, the bracket can be defined by the relations
$$[e_1e_2]=e_3\ ,\qquad [e_1e_3]=-e_2\ ,\qquad [e_2e_3]=e_1\ .\eqno(1.2.2a)$$
This Lie algebra, denoted
$A_1$ or sl$_2(\C)$, can be called the `mother of all (semi-simple) Lie
algebras'. A more familiar realisation of $A_1$ uses a basis $\{e,f,h\}$
with relations
$$[ef]=h\ ,\qquad [he]=2e\ ,\qquad [hf]=-2f\ .\eqno(1.2.2b)$$
The reader can find the change-of-basis (valid over $\C$ but not $\R$)
showing that (1.2.2) define isomorphic {\it complex} (but
not {\it real}) Lie algebras.

Another important 3-dimensional Lie algebra is called the
Heisenberg algebra\footnote{$^3$}{{\smal Actually, `Heisenberg algebra' refers
to a family of Lie algebras, with (1.2.3) being the one of lowest dimension.}}
 and is the algebra of the canonical commutation relations
in quantum mechanics: choosing a basis $x,p,h$, it is defined by
$$[xp]=h\ ,\qquad [xh]=[ph]=0\ .\eqno(1.2.3)$$

 From our definition, it is far from clear that Lie algebras, as a class, should be
natural and worth studying. After all, there are infinitely many possible
axiomatic systems: why should the one defining a Lie algebra be anything
special {\it a priori}? Perhaps this could have been anticipated by the
following line of reasoning.\smallskip

\noindent{{\bf Axiom.}} Groups are important and interesting.

\noindent{{\bf Axiom.}} Manifolds are important and interesting.\smallskip

Manifolds are structures where calculus is possible; locally, a manifold
looks like a piece of $\R^n$ (or $\C^n$), but these pieces can be bent
and stitched  together to create more interesting shapes. For instance
a circle is a 1-dimensional manifold, while Einstein claimed space-time
is a curved 4-dimensional one.

\smallskip
\noindent{{\bf Definition.}} A Lie group is a manifold with a compatible
group structure.\smallskip

This means that `multiplication' and `inverse' are differentiable maps. $\R$ is
a Lie group, under addition: obviously, $\mu:\R^2\rightarrow\R$ and
$\iota: \R\rightarrow\R$ defined by $\mu(a,b)=a+b$ and $\iota(a)=-a$
are both differentiable. (Why isn't $\R$ a Lie group under
multiplication?) A
circle is also a Lie group: parametrise the points with the angle $\theta$
defined mod $2\pi$ (or mod 360 if you prefer); the `product' of the point at
angle $\theta_1$ with the point at angle $\theta_2$ will be the point
at angle $\theta_1+\theta_2$. Surprisingly, the only other $k$-sphere which is a
Lie group is $S^3$ (the product can be defined using quaternions of unit
length\footnote{$^4$}{{\smal Similarly, the 7-sphere inherits from the octonions
a {\smit nonassociative} (hence nongroup) product, compatible with its
manifold structure.}}, or using the matrix group SU$_2(\C)$). Many but
not all Lie groups can be expressed as matrix groups. Two other examples
are GL$_n$ (invertible $n\times n$ matrices) and SL$_n$ (ones with determinant
1).

A consequence of the above axioms is then surely:

\smallskip
\noindent{{\bf Corollary.}} Lie groups should be important and interesting.
\smallskip

Lie group structure theory can be thought of as a major generalisation of
linear algebra. The basic constructions familiar to undergraduates have
important analogues valid in many Lie groups. For instance, years ago
we were taught to solve linear equations and invert matrices by using
elementary row operations to reduce a matrix to row-echelon form. What
this says is that any matrix $A\in {\rm GL}_n(\C)$ can be factorised
$A=BPN$, where $N$ is uppertriangular with 1's on the diagonal, $P$ is
a permutation matrix, and $B$ is an uppertriangular matrix. This is
essentially what is called the Bruhat decomposition of the Lie group
GL$_n(\C)$. More generally (where it applies to any `reductive' Lie
group $G$), $P$ will be an element of the so-called `Weyl group' of $G$
(of which we'll have much  more to say later), and $B$ will be in a `Borel
subgroup'.

Lie groups appear throughout physics. E.g.\ the orthogonal group SO$_3(\R)$
is the configuration space of a rigid body centred at the origin, while
SU$_2(\C)$ is the set of states of an electron at rest. The gauge group
of the Standard Model of particle physics is SU$_3(\C)\times{\rm SU}_2(\C)
\times{\rm U}_1(\C)$, while the Lorentz group of special relativity is
SO$_{3,1}(\R)$.

There is an important relation between Lie groups and Lie algebras.

\smallskip
\noindent{{\bf Fact.}} The tangent space of a Lie group
is a Lie algebra.
Any (finite-dimensional real or complex) Lie algebra is the tangent space
to some Lie group.\smallskip

More precisely, the tangent space at 1 (i.e.\ the set of all infinitesimal
generators of the Lie group) can be given a natural Lie algebra structure.
A Lie algebra, being a linearised Lie group, is much simpler
and easier to handle.
The Lie algebra preserves the local properties of the Lie group, though it
loses global topological properties (like boundedness). A Lie group has
a single Lie algebra, but a Lie algebra can correspond to many different
Lie groups. The Lie algebra corresponding to both $\R$ and $S^1$
is ${\frak g}=\R$ with trivial bracket. The Lie algebra corresponding to both $S^3
={\rm SU}_2(\C)$ and SO$_3(\R)$
is the cross-product algebra on $\R^3$ (usually called so$_3(\R)$). Given
the above fact, a safe guess would be:\smallskip

\noindent{{\bf Conjecture.}} Lie algebras are important and interesting.\smallskip

 From this line of reasoning, it should be expected that historically
 Lie groups arose first. Indeed that is the case: the Norwegian Sophus Lie
 introduced
them in 1873 to try to develop a Galois theory for ordinary differential equations.
As the reader may be aware, Galois theory is used for instance to show that not all
5th degree (or higher) polynomials can be explicitly `solved' using radicals
--- we will meet Galois theory in \S 1.8. Lie
wanted to study the explicit solvability (integrability) of differential
equations, and this led him to develop what we now call Lie theory.
The importance of Lie groups however have grown well beyond this initial
motivation.

An important class of  Lie algebras are the so-called finite-dimensional
{\it simple} ones. Their definition and motivation will be studied in
\S 2.7 below, but in a certain sense they serve as building blocks for all
other finite-dimensional Lie algebras.

The classification of simple finite-dimensional Lie algebras over $\C$ is
quite important and was accomplished at the turn of the century by Killing
and Cartan. There are 4 infinite families $A_\ell$ ($\ell\ge 1$),
$B_\ell$ ($\ell\ge 3$), $C_\ell$ ($\ell\ge 2$), and $D_\ell$ ($\ell\ge 4)$,
and 5 exceptionals $E_6$, $E_7$, $E_8$, $F_4$ and $G_2$. $A_\ell$ can be
thought of\footnote{$^5$}{{\smal Strictly speaking these are {\smit representations}
(see next section).}} as sl$_{\ell+1}(\C)$, the $(\ell +1)\times(\ell+1)$
matrices with trace 0. The orthogonal algebras $B_\ell$ and $D_{\ell}$
can be identified with so$_{2\ell+1}(\C)$ and so$_{2\ell}(\C)$, resp.,
where so$_n(\C)$ is all $n\times n$ anti-symmetric matrices $A^t=-A$.
The symplectic algebra $C_\ell$ is sp$_{2\ell}(\C)$, i.e.\ all $2\ell\times
2\ell$ matrices $A$ obeying $A\Omega=-\Omega A^t$, where $\Omega=\left(
\matrix{0&I_\ell\cr -I_\ell&0}\right)$ and $I_\ell$ is the identity matrix.
The exceptionals can be constructed e.g.\ using the octonions. In all these cases
the bracket is given by the commutator
$$[AB]=[A,B]:=AB-BA\eqno(1.2.4)$$
(it is a good exercise for the reader to confirm that the commutator satisfies
(1.2.1), and that e.g.\ sl$_{n}(\C)$ is indeed closed under it). To see
that (1.2.2b) truly is sl$_2(\C)$, put
$$e=\left(\matrix{0&1\cr 0&0}\right)\ ,\qquad f=\left(\matrix{0&0\cr 1&0}
\right)\ ,\qquad h=\left(\matrix{1&0\cr 0&-1}\right)\ .\eqno(1.2.5)$$
Incidentally the names $A$, $B$, $C$, $D$ have no significance: since the 4
 series start at $\ell= 1,2,3,4$, it seemed natural to call these $A,B,C,D$,
 resp. Unfortunately a bit of bad luck happened: $B_2$ and $C_2$ are
 isomorphic and so at random that algebra was placed in the orthogonal
series; however affine Dynkin
 diagrams make it clear that it really is a symplectic algebra which accidentally
 looks orthogonal; hence in hindsight
 the names of the $B$- and $C$-series really should have been switched.

This classification changes if the {\it field} --- the choice of scalars=numbers
--- is changed. By a field, we mean
we can add, subtract, multiply and divide, such that all the usual properties
like commutativity and distributivity are obeyed. Fields will make a few
different appearances in these notes. $\C$, $\R$, and $\Q$ are fields,
while $\Z$ is not (you can't always divide an integer by e.g.\ 3, and remain
in $\Z$). The integers mod $n$, which we will write $\Z_n$, are a field
iff $n$ is prime (the reader can verify that in e.g.\ $\Z_4$, it is not possible
to divide by the field element $[2]\in\Z_4$ even though $[2]\ne[0]$ there).
$\C$ and $\R$ are examples of fields of characteristic 0 ---
this means that 0 is the only integer $k$ with the property that $kx=0$
for all $x$ in the field. $\Z_p$
is the simplest example of a field with nonzero characteristic: in $\Z_p$,
multiplying by the integer $p$ has the same effect as multiplying by 0, and so we
say $\Z_p$ has characteristic $p$. Strange fields have important applications
in e.g.\ coding theory and, ironically, in number theory itself --- see e.g.\
\S1.8.

As always, $\C$ is better behaved than e.g.\ $\R$ because every
polynomial can be factorised over $\C$ (we say $\C$ is {\it algebraically
closed}) --- this implies for example that every matrix has an eigenvector over
$\C$ but not necessarily over $\R$. Over $\R$, the difference in the simple
Lie algebra
classification is that each symbol $X_\ell\in\{A_\ell,\ldots,G_2\}$
corresponds to a number of inequivalent algebras (over $\C$, each
algebra has its own symbol). For example, `$A_1$'
corresponds to 3 different real simple Lie algebras, namely the matrix
algebras sl$_2(\R)$, sl$_2(\C)$ (interpreted as a {\it real} vector space),
and so$_3(\R)$. The simple Lie algebra classification has recently been
done in any characteristic $p>7$. It is surprising but very common that
the smaller primes behave very poorly, and the classification for characteristic
2 is probably completely hopeless.

Associated with each simple algebra $X_\ell$ is a Weyl group, and a (Coxeter-)Dynkin
diagram. The Weyl group is a finite reflection group,
e.g.\ for $A_\ell$ it is the symmetric group ${\frak S}_{\ell+1}$. See
Figure 7 of [59] for the Weyl group of $A_2$. The Dynkin diagram
of $X_\ell$ (see e.g.\ [24,36,38] or Figure 6 in [59]) is a graph with
$\ell$  nodes, and with possibly some double
and triple edges. It says how to construct $X_\ell$ abstractly using
generators  and relations --- see \S 2.7. 
We will keep meeting both throughout these notes.

Another source of Lie algebras are the vector fields ${\rm Vect}(M)$ on a manifold
$M$. A vector field $v$ is a choice (in a smooth way) of a tangent
vector $v(p)\in T_pM$ at each point of $M$. It can be thought of as a
(1st order) differential operator, acting on functions $f:M\rightarrow \R$
(or $f:M\rightarrow\C$);
at each point $p\in M$ take the directional derivative of $f$ in the direction
$v(p)$. For example the vector fields on the circle, ${\rm Vect}(S^1)$, can be
thought of as anything of the form $g(\theta){d\over d\theta}$ where $g(\theta)$
can be any function with period 1. We can compose vector fields $u\circ v$,
but this will result in a 2nd order differential operator: e.g.\
$$(f(\theta){d\over d\theta})\circ(g(\theta){d\over d\theta})=f(\theta)\,g(
\theta)\,{d^2\over d\theta^2}+f(\theta)\,g'(\theta)\,{d\over d\theta}\ .$$
Instead, the natural `product' of vector fields is given by their commutator
$[u,v]=u\circ v-v\circ u$, as it always results in a vector field: e.g.\
$$[f(\theta){d\over d\theta},g(\theta){d\over d\theta}]=(f(\theta)\,g'(\theta)
-f'(\theta)\,g(\theta)){d\over d\theta}$$
 in ${\rm Vect}(S^1)$. ${\rm Vect}(M)$ with
this bracket is an infinite-dimensional Lie algebra. In the case where
$M$ is a Lie group $G$, the Lie algebra of $G$ can be interpreted as a
certain finite-dimensional subalgebra of ${\rm Vect}(G)$ given by the `left-invariant
vector fields'.

Simple algebras need not be finite-dimensional. An example of an infinite-dimensional
one is the {\it Witt algebra} ${\cal W}$, which can be defined (over $\C$)
by the basis\footnote{$^6$}{{\smal In infinite dimensions, to avoid
convergence complications, only finite linear combinations of basis vectors are
generally permitted. Infinite linear combinations would involve taking
some `completion'.}} $L_n$, $n\in\Z$, and the relations
$$[L_mL_n]=(m-n)L_{m+n}\ .\eqno(1.2.6)$$
Using the realisation $L_n=-\i e^{-\i n\theta}{d\over d\theta}$, the Witt
algebra can also be interpreted as the polynomial subalgebra of the
complexification $\C\otimes {\rm Vect}(S^1)$ --- i.e.\ change the
scalar field of Vect$(S^1)$ from $\R$ to $\C$. Incidentally, infinite-dimensional
Lie algebras need not have  a corresponding Lie group: e.g.\ the real
algebra ${\rm Vect}(S^1)$ is the Lie algebra of the Lie group ${\rm Diff}^+(S^1)$
of orientation-preserving diffeomorphisms $S^1\rightarrow S^1$, but
$\C\otimes {\rm Vect}(S^1)$ has no Lie group. ${\rm Diff}^+(S^1)$ plays a large role
in CFT, by acting on the objects of Segal's category ${\bf C}$.

The Witt algebra appears naturally in CFT: e.g.\ using the realisation
$L_n=-z^{n+1}{d\over dz}$ it is the polynomial subalgebra of the Lie
algebra ${\rm Vect}(\C/\{0\})$. Very carelessly, ${\rm Vect}(\C/\{0\})$ is often
thought of as the infinitesimal conformal transformations on a suitable
neighbourhood of 0 (yet clearly $L_{-2}$, $L_{-3}$, ... are singular at 0!).
Indeed the CFT literature is very sloppy when discussing the conformal
group in 2-dimensions. The unfortunate fact is that, contrary to claims,
{\it there is no infinite-dimensional conformal group} for $\C\cong\R^2$. The
best we can do is the 3-dimensional group PSL$_2(\C)$ of M\"obius
transformations $z\mapsto {az+b\over cz+d}$, which are orientation-preserving
conformal transformations for the Riemann sphere $\C\cup\{\infty\}$.
There seem to be 2 ways out of this rather embarrassing predicament. One is
to argue that we are really interested in `infinitesimal conformal invariance'
in some meromorphic sense, so the full Witt algebra can appear. The other way
is to argue that it is the conformal group of `Minkowski space'
$\R^{1,1}$ (or  better, its
compactification $S^1\times S^1$) rather than $\R^2\cong \C$ (or its
compactification $S^2$) which is relevant for CFT. That conformal group
{\it is} infinite-dimensional; for $S^1\times S^1$ it consists of 2 copies
of ${\rm Diff}^+(S^1)\times {\rm Diff}^+(S^1)$. For a more careful treatment
of this point, see [51].

For reasons we will discuss in \S 1.4, we are more interested in the {\it
Virasoro algebra} ${\cal V}$ rather than the Witt algebra ${\cal W}$. This is a `1-dimensional
central extension' of ${\cal W}$; as a vector space ${\cal V}={\cal W}\oplus
\C C$ with relations given by
$$\eqalignno{[L_mL_n]=&\,(m-n)L_{m+n}+\delta_{n,-m}{m\,(m^2-1)\over 12}\,C&(1.2.7a)\cr
[L_mC]=&\,0\ .&(1.2.7b)\cr}$$
`1-dimensional central extension' means ${\cal V}$ has one extra  basis
vector $C$, which lies in the {\it centre} of ${\cal V}$ (i.e.\ $[xC]=0$
for all $x\in{\cal V}$), and sending $C\rightarrow 0$ recovers ${\cal
W}$ (i.e.\ takes (1.2.7a) to (1.2.6)).
A common mistake in the physics literature is to regard $C$ as a number:
it is in fact a vector, though in many (but not all) representations it is
 mapped to a scalar multiple of the identity.

The reason for the strange-looking (1.2.7a) is that we have little choice:
${\cal V}$ is the unique nontrivial 1-dimensional central extension of
${\cal W}$. The factor ${1\over 12}$ there is conventional but standard,
and has to do with `zeta-function regularisation' in string theory --- i.e.\
the divergent sum $\sum_{n=1}^\infty n$ is `reinterpreted' as
$\zeta(-1)={-1\over 12}$, where $\zeta(s)=\sum_{n=1}^\infty n^{-s}$ is
the Riemann zeta function. Incidentally $\zeta(s)$ can be written as
the product $\prod_p(1-p^{-s})^{-1}$ over all primes $p=2,3,5,\ldots$
(try to see why);
hence $\zeta(s)$ has a lot to do with primes, in particular their distribution.
In fact the most famous
unsolved problem in math today is the Riemann conjecture, which says that
$\zeta(s)\ne 0$ whenever Re$(s)\ne {1\over 2}$. One researcher recently described
this conjecture as saying that the primes have music in them.

In CFT, $L_0$ is the energy operator. For example the partition function
is given by $\z(\tau)={\rm Tr}_H(q^{L_0-c/24}q^{*\bar{L}_0
-c/24})$ and the (normalised) character $\chi_a$ equals ${\rm Tr}_{H_a}(q^{L_0-c/24})$
for $q=e^{2\pi\i\tau}$. $cI$ is the scalar multiple of the identity to which $C$
gets sent; it has a physical interpretation [18] involving Casimir (vacuum)
 energy, which depends on space-time topology,
and the strange shift by $c/24$ is due to an implicit mapping from
the cylinder to the plane.

\bigskip\bigskip\noindent{{\smcap 1.3. Representations of finite-dimensional
simple Lie algebras}} \bigskip

The representation theory of the simple Lie algebras\footnote{$^7$}{{\smal
See e.g.\ [25] for more details. Historically, representations of Lie algebras were considered even before
representations of finite groups.}} can perhaps be regarded
as an enormous generalisation of trigonometry. For instance the
facts that ${\sin(nx)\over \sin(x)}$ can be written
as a polynomial in $\cos(x)$ for any $n\in\Z$, and that
$${\sin(mx)\,\sin(nx)\over \sin(x)}=\sin((m+n)x)+\sin((m+n-2)x)+\cdots+
\sin((m-n)x)$$
for any $m,n\in\Z_>$, are both easy special cases of the theory.

The classic example of an algebraic structure are the numbers, and they
prejudice us into thinking that commutativity and associativity are
the ideal. We have learned over the past couple of centuries 
that commutativity can often be dropped without losing depth and usefulness,
but most interesting structures seem to obey some form of associativity.
Moreover, true associativity (as opposed to e.g.\ anti-associativity) really
simplifies the arithmetic. Given the happy `accident' that the
commutator $[x,y]:=xy-yx$ in any associative algebra obeys anti-associativity,
it would seem to be both tempting and natural to study the ways (if any)
in which associative algebras ${\frak A}$ can `model' or {\it represent} a given
Lie algebra. Precisely, we are looking for a map $\rho:{\frak g}\rightarrow
{\frak A}$ which preserves the linear structure (i.e.\ $\rho$ is a
linear function), and which sends the bracket $[xy]$ in ${\frak g}$ to the
commutator $[\rho(x),\rho(y)]$ in ${\frak A}$.

In practice groups (resp., algebras) often appear as symmetries (resp.,
infinitesimal generators of symmetries). These symmetries often act
linearly. In other words, in practise the preferred associative algebras
will usually be matrix algebras, and this is the usual form for
a representation and the only kind we will consider. The {\it dimension} of these representations is the
size of the matrices.

Finding all possible representations, even for the simple Lie algebras,
is probably hopeless. However, it is possible to find all {\it finite-dimensional}
representations of the simple Lie algebras, and the answer is easy to
describe. Given a simple Lie algebra $X_\ell$, there is a
representation $L_\la$ for each $\ell$-tuple $\la=(\la_1,\ldots,\la_\ell)$
of nonnegative integers. $\la$ is called a {\it highest-weight}. 
Moreover, we can take direct sums $\oplus_iL_{\la^{(i)}}$
of finitely many of these representations.
 The matrices in such a direct sum will be in
block form. It turns out that, up to change-of-basis, this exhausts all
finite-dimensional representations of $X_\ell$.

It is common to replace `representation $\rho$ of ${\frak g}$' with the equivalent
notion of `${\frak g}$-module $M$' --- i.e.\ we think of the matrices $\rho(x)$
as linear maps $M\rightarrow M$. A ${\frak g}${\it -module} is a vector space
on which ${\frak g}$ acts (on the left). Instead of considering the matrix
$\rho(x)$, we consider `products' $xv$ (think of this as the matrix
$\rho(x)$ times the column vector $v$) for $v\in M$. This product must be
bilinear, and must obey $[xy]v=x(yv)-y(xv)$.

To get an idea of what $L_\la$ looks like, consider $A_1$. Recall its
generators $e,f,h$ and relations (1.2.2b). Choose any 
$\la\in\C$. Define $x_0\ne 0$ to formally obey $hx_0=\la x_0$ and $ex_0=0$.
Define
inductively $x_{i+1}:=fx_i$ for $i=0,1,\ldots$. Define $M_\la$ to be
the span of all $x_i$ --- we will see shortly that they are linearly independent
(so $M_\la$ is infinite-dimensional).
$M_\la$ is a module of $A_1$: the calculations
$hx_{i+1}=hfx_i=([hf]+fh)x_i=(-2f+fh)x_i$ and 
$ex_{i+1}=efx_i=([ef]+fe)x_i=(h+fe)x_i$ show  inductively that
$hx_m=(\la-2m)x_m$ and $ex_m=(\la-m+1)m\,x_{m-1}$. From these the reader can show
that the $x_i$ are linearly independent. $M_\la$ is called a {\it Verma module};
$\la$ is called its highest-weight, and $x_0$ is called a highest-weight
vector.

Now specialise to $\la=n\in\Z_\ge:=\{0,1,2,\ldots\}$.
Note that $ex_{n+1}=0$ and $hx_{n+1}=(-n-2)x_{n+1}$. This means that, for
these $n$, $M_n$ contains a {\it submodule} with highest-weight vector
$x_{n+1}$,  isomorphic to $M_{-n-2}$.
$x_{n+1}$ is called a {\it null vector}.
In other words, we could set $x_{n+1}:=0$ and still have an $A_1$-module.
 We would then get a {\it finite-dimensional} module
which we'll call $L_n:=M_n/M_{-n-2}$ (not to be confused with the
Virasoro generator in (1.2.7)). Its basis is $\{x_0,x_1,\ldots,x_n\}$
and so it has dimension $n+1$.

For example, take $n=1$. Note that what we get in terms of the basis $\{x_0,
x_1\}$ is the familiar representation
sl$_2(\C)$ given in (1.2.5).

The situation for the other simple Lie algebras $X_\ell$ is similar.

It turns out to be hard to compare representations:  $\rho$ and
$\rho'$ could be equivalent (i.e.\ differ merely by a change-of-basis)
but look very different. Or if we are given
a representation, we may want to decompose it into the direct sum of some
$L_{\la^{(i)}}$. When working with representations, it is often very useful to
avoid much of the extraneous basis-dependent detail present in the function $\rho$.
Finite group theory suggests how to do this: we should use {\it characters}.
The character of an $A_1$-module $M$ is given by Weyl:
write $M$ as a direct sum of eigenspaces
$M(m)$ of $h$; then define
$${\rm ch}_M(z):=\sum_m{\rm dim}\ M(m)\ e^{mz}\ ,\eqno(1.3.1)$$
for any $z\in\C$. The $m$ are called {\it weights} and the $M(m)$ {\it
weight-spaces}.
For example, for $L_n$ the weights are $m=n,n-2,\ldots,-n$, the
weight-spaces $L_n(m)$ are $\C \,x_{(n-m)/2}$, and
$${\rm ch}_n(z)=\sum_{i=0}^n  e^{(n-2i)z}={\sin((n+1)\,z)\over\sin(z)}\ .\eqno(1.3.2)$$
Analogous formulas apply to any  algebra $X_\ell$: the character will
then be a function of an $\ell$-dimensional subspace ${\frak h}$ called the {\it Cartan
subalgebra}, spanned by all the $h_i$ (see \S 2.7), so can be thought of
as a complex-valued function of $\ell$ complex variables. The weights $m$
will lie in the dual space to ${\frak h}$ --- i.e.\ are linear maps
${\frak h}\rightarrow \C$ --- so will have $\ell$
components. See for instance Figure 8 in [59]. Incidentally, $\ell$ is
called the {\it rank} of $X_\ell$.

Weyl's definition works: two representations are equivalent iff their
characters are identical, and $M=\oplus_iL_{\la^{(i)}}$ iff
ch$_M(z)=\sum_i{\rm ch}_{\la^{(i)}}(z)$. It also is enormously simpler: e.g.\
the smallest nontrivial representation of $E_8$ is a map from
$\C^{248}$ to the space of $248\times 248$ matrices, while its character
is a function $\C^8\rightarrow \C$. But why is Weyl's definition natural? How
 did he come up with it?

To answer that question, we must remind ourselves of the characters of
finite groups\footnote{$^8$}{{\smal Surprisingly, what we now call the characters of
group representations were invented almost a decade before group representations
were.}}. A representation of a finite group $G$ is a structure-preserving
map $\rho$ (i.e.\ a group homomorphism) from $G$ to matrices. The group's product
becomes matrix product. In these notes we will be exclusively interested
in group representations over $\C$.
Two representations $\rho,\rho'$ are called {\it equivalent} if there exists
a matrix (change-of-basis) $U$ such that $\rho'(g)=U\rho(g)U^{-1}$ for
all $g$. The character ${\rm ch}_\rho$ is the map $G\rightarrow\C$ given by
the trace: ch$_\rho(g)={\rm tr}(\rho(g))$. We see that equivalent
representations will have the same character, because of the fundamental
identity ${\rm tr}(AB)={\rm tr}(BA)$. This identity also tells us that the character
is a `class function', i.e.\ ch$_\rho(hgh^{-1})={\rm tr}(\rho(h)\,\rho(g)\,
\rho(h)^{-1})={\rm ch}_\rho(g)$
so ch$_\rho$ is constant on each `conjugacy class'. Group characters
are also enormously simpler than representations: e.g.\ the smallest
nontrivial representation
of the Monster $\M$ (see Part 2) consists of almost $10^{54}$ matrices,
each of size $196883\times 196883$, while its character consists of 194
complex numbers. Incidentally, finite group representations behave analogously
to the representations of $X_\ell$: the role of the modules $L_\la$ is
played by the irreducible representations $\rho_i$, and any finite-dimensional
representation of $G$ can be decomposed uniquely into a direct sum of
various $\rho_i$. The difference is that there are only finitely many
$\rho_i$ --- their number equals the number of conjugacy classes of $G$.

We can use this group intuition here. In particular, given any Lie algebra
$X_\ell$ and representation $\rho$, we can think of the map  $e^x\mapsto e^{\rho(x)}$
as a representation of a Lie group $G(X_\ell)$ corresponding to $X_\ell$ (the exponential
$e^A$ of a matrix is defined by the usual power series; it will always
converge). The trace of the matrix $e^{\rho(x)}$ will be the {\it group} character
value at $e^x\in G(X_\ell)$, so we'll define it to be the {\it algebra} character
value at $x\in X_\ell$. Again, it suffices to consider only representatives
of each conjugacy class of $G(X_\ell)$, because the character will be a class
function. Now, almost every matrix is diagonalisable (since almost
any $n\times n$ matrix has $n$ distinct eigenvalues), and so it would
seem we aren't losing much by restricting $x\in X_\ell$ to {\it diagonalisable}
matrices. Hence we may take our conjugacy class representatives to
be {\it diagonal} matrices $x\in X_\ell$, i.e.\ (for $X_\ell=A_1$) to
$x=zh$  for $z\in\C$ ($h$ is diagonal in the $x_i$ basis of $L_\la$). 
So the algebra character can be chosen to be a function of
$z$. Finally, the trace of $e^{\rho(x)}=e^{z\rho(h)}$ will be given by
(1.3.1). This completes the motivation for Weyl's character formula.

There is one other important observation we can make. Different diagonal
matrices can belong to the same conjugacy class. For instance,
$$\left(\matrix{0&-1\cr 1&0}\right)\,\left(\matrix{a&0\cr 0&b}\right)
\,\left(\matrix{0&-1\cr 1&0}\right)^{-1}=\left(\matrix{b&0\cr 0&a}\right)\ ,$$
so $e^{zh}$ and $e^{-zh}$ lie in the same $G({A_1})={\rm SL}_2(\C)$ conjugacy class.
Hence ch$_M(z)={\rm ch}_M(-z)$. This symmetry $z\mapsto -z$ belongs to the
{\it Weyl group} for $A_1$. Each $X_\ell$ has similar symmetries, and
the Weyl group plays an
important role in the whole theory, sort of analogous to the modular
group for modular functions we'll discuss in \S 2.3. 

Weyl found a generalisation of the right-side of (1.3.2), valid for
all  $X_\ell$.
The character of $L_\la$ can be written as a fraction (2.8.1): the numerator
will be a alternating sum over the Weyl group, and the denominator
will be a  product over `positive roots'. This formula and its
generalisations have profound consequences, as we'll see in \S 2.8.

Incidentally, the trigonometric identities given at the beginning of this
section are the tensor product formula of representations (interpreted as the
product and sum of characters), and the fact that an arbitrary character
can be written as a polynomial in the fundamental characters, both
specialised to $A_1$ (see (1.3.2) for the $A_1$ characters).

\bigskip\bigskip\noindent{{\smcap 1.4. Affine algebras and the Kac-Peterson
matrices}}\bigskip

The theory of nontwisted affine Kac-Moody algebras (usually called
{\it affine algebras} or {\it current algebras}) is extremely analogous
to that of the finite-dimensional simple Lie algebras. Nothing infinite-dimensional
tries harder to be finite-dimensional than affine algebras. Standard
references for the following material are [38,41,24].

Let $X_\ell$ be any simple finite-dimensional Lie algebra. The affine algebra
$X_\ell^{(1)}$ is essentially the {\it loop algebra} ${\cal L}(X_\ell)$, defined to
be all possible `Laurent polynomials' $\sum_{n\in\Z}a_nt^n$ where
each $a_n\in X_\ell$ and all but finitely many $a_n=0$. $t$ here is an
indeterminant. The bracket in
${\cal L}(X_\ell)$ is the obvious one: e.g.\ $[at^n,bt^m]=[ab]t^{n+m}$.
Geometrically, ${\cal L}(X_\ell)$ is the Lie algebra of polynomial maps $S^1\rightarrow
X_\ell$ --- hence the name (for that realisation, think of
$t=e^{2\pi\i\theta}$).  Hence there are many  generalisations
of the loop algebra (e.g.\ any manifold in place of $S^1$ will do), closely
 related ones called {\it toroidal algebras} being the Lie algebra of maps
 $S^1\times\cdots\times S^1\rightarrow
 X_\ell$. But the loop algebra is simplest
 and best understood, and the only one we'll consider. Note that
 ${\cal L}(X_\ell)$ is infinite-dimensional. Its Lie groups are the {\it loop
 groups}, consisting of all loops $S^1\rightarrow G(X_\ell)$ in a
 Lie group for $X_\ell$.

We saw $S^1$ before, in the discussion of the Witt algebra. Thus the
Virasoro and affine algebras should be related. In fact, the Virasoro
algebra acts on the affine algebras as `derivations', and this
connection plays an important technical role in the theory.

$X_\ell^{(1)}$ is in the same relation to the loop algebra, that the Virasoro
${\cal V}$
is to the Witt ${\cal W}$. Namely, it is its (unique nontrivial 1-dimensional) central
extension --- see e.g.\ (7.7.1) of [38] for the analogue of (1.2.7a)
here.  In addition, for more technical
reasons, a further (noncentral) 1-dimensional extension is usually made:
 the derivation $t{d\over dt}$ is included (see footnote 33). $X_\ell^{(1)}$
is the simplest of the infinite-dimensional Kac-Moody algebras.
The superscript `(1)' denotes
the fact that the loop algebra was twisted by an order-1 automorphism --- i.e.\
that it is untwisted. It is called `affine' because of its Weyl group,
as we shall see.

Central extensions are a common theme in today's infinite-dimensional
Lie theory\footnote{$^9$}{{\smal Incidentally the {\smit finite}-{\smit dimensional}
simple Lie algebras do not have nontrivial central extensions.}}. Their
{\it raison d'\^ etre} is always the same: a richer supply of representations.
 For example, ${\cal W}$ has several representations, but no nontrivial one
is an `irreducible unitary positive-energy representation' --- the kind of greatest
interest in math phys. On the other hand, its central extension ${\cal V}$
has a rich supply of those representations (e.g.\ there's one for each
choice of $c>1,h>0$, namely the Verma module $V_{c,h}$ corresponding
to $L_0x_0=hx_0,Cx_0=cx_0$). At the level of groups, central extensions
allow {\it projective} representations (i.e.\ representations up to a
scalar factor) to become true representations.
Projective representations (hence central extensions) appear naturally
in QFT because a quantum state vector $|v\rangle$ is
physically indistinguishable from any nonzero scalar multiple $\alpha|v\rangle$.

All of the quantities associated to $X_\ell$ have an analogue here: Dynkin
diagram, Weyl group, weights,... For instance, the affine Dynkin
diagram is obtained from the Dynkin diagram for $X_\ell$ by adding one node.
See for example Figure 9 of [59]. The extra node is always labelled by a
`0'. The Cartan subalgebra ${\frak h}$ here will be $(\ell+2)$-dimensional. Many of
these details will be discussed in more detail in \S 2.7 below.

The construction of $X_\ell^{(1)}$ is so trivial that it seems surprising anything interesting
and new can happen here. But a certain `miracle' happens...

No interesting representation of $X_\ell^{(1)}$ is finite-dimensional.
The analogue for $X_\ell^{(1)}$ of the finite-dimensional representations of $X_\ell$ are
called the {\it integrable highest-weight representations}, and will be
denoted $L_\la$. The {\it highest-weight} $\la$ here will be an $(\ell+1)$-tuple
$(\la_0,\la_1,\ldots,\la_\ell)$, $\la_i\in\Z_{\ge}$ (strictly speaking, it
will be an $(\ell+2)$-tuple, but the extra component is not important
and is usually ignored). As for $X_\ell$, the highest-weights can be
thought of as the assignment of a nonnegative integer to each node of the
Dynkin diagram. The construction of $L_\la$ is as in the finite-dimensional
case. They are called {\it integrable} because they are precisely those
highest-weight representations which can be `integrated' to a projective
representation of the corresponding loop group, and hence a representation
of a central extension of the loop group.

We define the character $\chi_\la$ as in (1.3.1), though now the weights
$m$ will be $(\ell+2)$-tuples, and there will be infinitely many of them.
$\chi_\la$ will be a complex-valued
function of $\ell+2$ complex variables $(\vec{z},\tau,u)$ (see (1.4.1a)
below). It can be written
as an alternating sum over the Weyl group $W$, over a `nice' denominator.
The difference here is that  $W$ is now infinite.

Perhaps most of the  interest in affine algebras can be traced to
the `miracle' that their Weyl groups are a semidirect product $Q^\vee\sdprod
\overline{W}$ of
translations in a lattice $\Q^\vee$ (the $\ell$-dimensional `co-root lattice' of $X_\ell$
--- see \S 1.6)
with the (finite) Weyl group $\overline{W}$ of $X_\ell$. See Figure 10
of [59] for the Weyl group of $A_2^{(1)}$. `Semidirect
product'\footnote{$^{10}$}{{\smal This is also discussed briefly in section 2.2.}}
 means that any element of $W$ can be written uniquely as $(t,w)$ for some
translation $t$ and some $w\in\overline{W}$, and $(t,w)\circ(t',w')=
({\rm stuff},w\circ w')$.

One thing this implies is that $\chi_\la$ will be of the form `theta
function'/denominator. Theta functions are classically-studied modular
forms (we will discuss these terms in \S 2.3), and thus the modular
group  SL$_2(\Z)$ will make an appearance!
To make this more precise, consider the highest-weight $\la=(\la_0,\la_1)$
of $A_1^{(1)}$, and write $k=\la_0+\la_1$. Then
$$\chi_\la={\Theta^{(k+2)}_{\la_1+1}-\Theta^{(k+2)}_{-\la_1-1}
\over \Theta^{(2)}_1-\Theta^{(2)}_{-1}}\eqno(1.4.1a)$$
where these functions all depend on 3 complex variables $z,\tau,u$, and
$$\Theta^{(n)}_m( {z},\tau,u):=e^{-2\pi \i nu}\sum_{\ell\in\Z+{m\over 2n}}
\exp[\pi \i n\tau\ell^2-2\sqrt{2}\pi \i n\ell z]\ .\eqno(1.4.1b)$$
In (1.4.1a) we can see the alternating sum over the Weyl group of $A_1$
in the numerator (and denominator,
since we've used the $A_1^{(1)}$ denominator identity in writing (1.4.1a)).
 For general $X_\ell^{(1)}$, the denominator will always be independent
of $\la$, and the theta function (1.4.1b) will become a multidimensional
one involving a sum over $Q^\vee$ shifted by some
weight and appropriately rescaled. The (co-)root lattice of $A_1$ is $\sqrt{2}\Z$.
The key variable in (1.4.1a)  is the modular one $\tau$, which will lie in the
upper half complex plane ${\cal H}$ (in order to have convergence).
In the applications to CFT, the other variables are often set to 0.

The number $k$ introduced in (1.4.1a) plays an important role in the general
theory. In the representation $L_\la$, the central term $C$ will get
sent to some multiple of the identity ---  the multiplier is labelled $k$
and is called the {\it level} of the representation. For any $X_\ell^{(1)}$
there is a simple formula expressing the level $k$ in terms of the
highest-weight $\la$; e.g.\ for $A_\ell^{(1)}$ and $C_\ell^{(1)}$ it is
given by $k=\la_0+\la_1+\cdots+\la_\ell$. Write $P_+^k$ for the (finite) set
of level $k$ highest-weights (so the size of $P_+^k$ for $A_\ell^{(1)}$
is $\left(k+\ell\atop\ell\right)$). An important weight in $P_+^k$
is $(k,0,\ldots,0)$. We will denote this `0'. In RCFT it corresponds to the vacuum.

The modular group SL$_2(\Z)$ acts on the Cartan subalgebra ${\frak h}$
of $X_\ell^{(1)}$ in the following way:
$$\left(\matrix{a&b\cr c&d\cr}\right)\,(\vec{z},\tau,u)=({\vec{z}\over
c\tau+d},{a\tau+b\over c\tau+d},u-{c\,\vec{z}\cdot\vec{z}\over 2(c\tau+d)})$$
Under this action, the characters $\chi_\la$ also transform nicely:
in particular we find for any level $k$ weight $\la$
$$\eqalignno{\chi_\la({\vec{z}\over \tau},{-1\over \tau},u-{\vec{z}\cdot
\vec{z}\over
2\tau})=&\,\sum_{\mu\in P_+^k}S_{\la\mu}\,\chi_\mu(\vec{z},\tau,u)&(1.4.2a)\cr
\chi_\la({\vec{z}},{\tau+1},u)
=&\,\sum_{\mu\in P_+^k}T_{\la\mu}\,\chi_\mu(\vec{z},\tau,u)&(1.4.2b)\cr}$$
where $S$ and $T$ are complex matrices called the {\it Kac-Peterson
matrices}. $S$ will always be symmetric and unitary, and has many
remarkable properties as we shall see. Its entries are related to
Lie group characters at elements of finite order (see (1.4.5) below).
$T$ is diagonal and unitary; its entries are related to the eigenvalues
of the quadratic Casimir.

For example, consider $A_1^{(1)}$ at level $k$. Then $S$ and $T$ will
be $(k+1)\times(k+1)$ matrices given by
$$ S_{\la\mu}=\sqrt{{2\over k+2}}\,\sin(\pi\,{(\la_1+1)(\mu_1
+1)\over k+2})\ ,
\qquad T_{\la\mu}=\exp[\pi\i\,{(\la_1+1)^2\over 2(k+2)}
-{\pi\i\over 4} ]\,\delta_{\la,\mu}\ .\eqno(1.4.3)$$

One important place $S$ appears is the famous {\it Verlinde formula}
$$N_{\la\mu}^\nu=\sum_{\kappa\in P_+^k}{S_{\la\ka}\,S_{\mu\ka}\,S_{\nu\ka}^*
\over S_{0,\ka}}\eqno(1.4.4)$$
for the fusion coefficients $N_{\la\mu}^\nu$ of the corresponding RCFT.
We will investigate some consequences of this formula in a later
section. The fusion coefficients for the affine algebras are
well-understood; see e.g.\ Section 4 of [59] for their interpretation (usually
called the Kac-Walton formula) as
`folded tensor product coefficients'.

We will see in \S 1.7 that symmetries of the extended Dynkin diagram have
consequences for $S$ and $T$ (simple-currents, charge-conjugation).
There is a `Galois action' on $S$ which we will discuss in \S 1.8.
There is a strange property of $S$ and $T$ called {\it rank-level
duality} (see e.g.\ [45]): the matrices for $A_\ell^{(1)}$ at level $k$ are closely
related to those of $A_{k-1}^{(1)}$ at level $\ell+1$, and similar
statements hold for $B_\ell^{(1)}$, $C_\ell^{(1)}$ and $D_\ell^{(1)}$.
Another reason $S$ is mathematically interesting is the formula
$${S_{\la\mu}\over S_{0\mu}}={\rm ch}_{\bar{\la}}(-2\pi \i\,{\overline{\mu+\rho}
\over k+h^\vee})\ .\eqno(1.4.5)$$
The right-side is a character of $X_\ell$, and $\bar{\la}=(\la_1,\ldots,\la_\ell)$
means ignore the extended node. $\rho$ is the `Weyl vector' $(1,1,\ldots,1)$
and $h^\vee$ is called the {\it dual Coxeter number} and is the level of
$\rho$. For $A_\ell^{(1)}$, $h^\vee=\ell+1$. Of course the right-side can
also be regarded as a character for a Lie group associated to $X_\ell$,
in which case the argument would have to be exponentiated and would
correspond to an element of finite order in the group. These numbers (1.4.5) have been studied
by many people (most extensively by Pianzola) and have some nice properties. For
instance Moody-Patera (1984)
have argued that exploiting them leads to some  quick algorithms
for computing e.g.\ tensor product coefficients. Kac [37] found a curious
application for them: a Lie theoretic proof of `quadratic reciprocity'!

Quadratic reciprocity is one of the gems of classical number theory.
It tells us that the equations
$$\eqalignno{x^2\equiv &\,a\ ({\rm mod}\ b)&\cr
y^2\equiv &\,b\ ({\rm mod}\ a)&\cr}$$
are related; more precisely, for fixed $a$ and $b$ (for simplicity take
them to both be primes$\,\ne 2$) the questions of whether there is a solution $x$
to the first equation and a solution $y$ to the second, are related.
They will both have the same yes or no answer, unless $a\equiv b\equiv 3$ (mod 4),
in which case they will have opposite answers. E.g.\ take $a=23$ and $b=3$,
then we know the first equation does not have a solution (since $a\equiv
2$ (mod 3) and $x^2\equiv 2$ doesn't have a solution mod 3), and hence
the second equation must have a solution (indeed, $y=7$ works).
There are now many proofs for quadratic reciprocity, and Kac used Lie
characters at elements of finite order to find another one.

What is interesting here is that Kac's proof uses only certain special
weights for $A_\ell$. The natural question is: is it possible to find any
generalisations of quadratic reciprocity using other weights and algebras?
Many generalisations of quadratic reciprocity are known; will generalising
Kac's argument recover them, or will they perhaps yield new reciprocity laws?
It seems no one knows.

The relation (1.4.5) is important because it connects finite-dimensional
Lie data with infinite-dimensional Lie data. The `conceptual arrow' can
be exploited both ways: in the generalisations of the arguments of \S 1.9 to
other algebras, (1.4.5) allows us to
use our extensive knowledge of finite-dimensional algebras to squeeze out
some information in the affine setting;
but also it is possible to use the richer symmetries of the affine data
to see `hidden' symmetries in finite-dimensional data. For example it
can be used (Gannon-Walton 1995) to find a sort
of Galois symmetry of dominant weight multiplicities in $X_\ell$, which would be
difficult or impossible to anticipate without (1.4.5).

\bigskip\bigskip\noindent{{\smcap 1.5. The classification of physical invariants}}
\bigskip

We are interested in the following classification problem. Choose any
affine algebra $X_\ell^{(1)}$ and level $k\in\Z_\ge$. Find all matrices $M=(M_{\la
\mu})_{\la,\mu\in P_+^k}$ such that\smallskip

\item{(P1)} $MS=SM$ and $MT=TM$, where $S,T$ are the 
Kac-Peterson matrices (1.4.2);

\item{(P2)} each entry $M_{\la\mu}\in\Z_{\ge}$;

\item{(P3)} $M_{00}=1$.\smallskip

Any such $M$, or equivalently the corresponding partition function ${\cal Z}=\sum_{\la,\mu}
M_{\la\mu}\chi_\la\chi_\mu^*$, is called a {\it physical invariant}.

The first and most important classification of physical invariants
was the Cappelli-Itzykson-Zuber A-D-E classification for $A_1^{(1)}$ at
all levels $k$ [8]. We will give their result shortly. This implies for instance the
minimal model RCFT classification, as well as the $N=1$ super(symmetric)conformal
minimal models. The other classifications of comparable magnitude are
$A_2^{(1)}$ for all $k$; $A_\ell^{(1)}$, $B_\ell^{(1)}$ and $D_\ell^{(1)}$
 for all $k\le 3$; $(A_1\oplus A_1)^{(1)}$ for all levels $(k_1,k_2)$;
 and $(u(1)\oplus\cdots\oplus u(1))^{(1)}$ for all (matrix-valued) levels
 $k$. See e.g.\ [29] for references. The most difficult of these classifications is
 for $A_2^{(1)}$, done by Gannon (1994).

In other words, very little in this direction has been accomplished in the 15 or so
years this problem has existed. But this is not really a good measure of progress.
The effort instead has been directed primarily towards the full classification;
most of these partial results are merely easy spin-offs from that more serious
and ambitious assault.

The proof in [8] was very complicated and followed the following lines.
First, an explicit basis was found for the vector space (called the
`commutant') of all matrices obeying (P1). Then
(P2) and (P3) were imposed.
Unfortunately their proof was long and formidable. 
Others tried to apply their approach to $A_2^{(1)}$, but without success.
The eventual proof for $A_2^{(1)}$ was completely independent of the
[8] argument, and exploited more of the structure implicit in the problem.
As the $A_2^{(1)}$ argument became more refined, it became the model
for the general assault. In \S 1.9 we sketch this new approach.

 From this more general perspective, of these completed
classifications  only the level 2 $B_\ell^{(1)}$ and $D_\ell^{(1)}$
ones will have any lasting value (the orthogonal algebras at level
2 behave very peculiarly, possess large numbers of exceptional physical
invariants, and must be treated separately). The others behave more generically
and will fall out as special cases once the more general classifications are
concluded. Other
classifications which should be straightforward with our present
understanding are $C_2^{(1)}$ at all $k$; $G_2^{(1)}$ at all $k$;
and $B_\ell^{(1)}$ and $D_\ell^{(1)}$ at $k=4$. The $C_2^{(1)}$  should be
easiest and would imply the $C_\ell^{(1)}$ level 2, as well as the
$B_\ell^{(1)}$ and $D_\ell^{(1)}$ level 5, classifications. A very safe conjecture
is that the only {\it exceptional} physical invariants (we define this
term in \S 1.7) for $C_2^{(1)}$ occur at $k=3,7,8,12$ --- this is known to be true
for all $k\le 500$. $G_2^{(1)}$ would be more difficult but also much
more valuable; its only known exceptionals occur at $k=3,4$, and these
are the only exceptionals for $k\le 500$, and a very safe conjecture is
that there are no other $G_2^{(1)}$ exceptionals. $B_\ell^{(1)}$ and
$D_\ell^{(1)}$ at level 4 will also be more difficult, but also would
be valuable; less is understood about its physical invariants and there
is a good chance new exceptionals exist there.

The most surprising thing about the known physical invariant classifications
is that there so few surprises: almost every physical invariant is
`generic'. We will see that the symmetries of the extended Dynkin
diagram give rise to general families of physical invariants. We will
call any physical invariants which do not arise in these generic ways
(i.e.\ using what are called simple-currents or conjugations), {\it exceptional}.
Many exceptionals have been found, and now we are almost at the point where
we can safely conjecture the complete list of physical invariants
for $X_\ell^{(1)}$ at any $k$, for $X_\ell$ a simple algebra.

Unfortunately
the classification for semi-simple algebras $X_{\ell_1}\oplus\cdots\oplus
X_{\ell_s}$ does not reduce to the one for simple ones.
In fact, {\it any explicit classification of the physical invariants
for $X^{(1)}$, for all semi-simple $X$, would easily be one of the
greatest  accomplishments
in the history of math,}  for it would include as a small part such
monumental things as an explicit classification of all 
positive-definite integral lattices. Thus we unfortunately
cannot expect an {\it explicit} classification for the {\it semi-}simple
algebras.

To make this discussion more concrete and explicit, consider $A_1^{(1)}$.
For convenience drop $\la_0$, so $P_+^k=\{0,1,\ldots,k\}$. Write $J$
for the permutation (called a simple-current) $Ja:=k-a$. Then
the complete list of
physical invariants for $A_1^{(1)}$ is 
$$\eqalignno{{\cal A}_{k+1}=&\,\sum_{a=0}^{k}\,|\chi_a|^2\ ,\qquad\qquad
\qquad {\rm for\ all}\ k\ge 1&
\cr{\cal D}_{{k\over 2}+2}=&\,\sum_{a=0}^{k}\,\chi_a\,\chi_{J^aa}^*\ ,\qquad 
\qquad\qquad{\rm whenever}\ {k\over 2}\ {\rm is\ odd}&\cr
{\cal D}_{{k\over 2}+2}=&\,|\chi_0+\chi_{J0}|^2+|\chi_2+\chi_{J2}|^2+\cdots
+2|\chi_{{k\over 2}}|^2\ ,\qquad {\rm whenever}\ {k\over 2}\ {\rm is\ even}&\cr
{\cal E}_6=&\,|\chi_0+\chi_6|^2+|\chi_3+\chi_7|^2+|\chi_4+\chi_{10}|^2\ ,
\qquad\qquad {\rm for}\ k=10&\cr
{\cal E}_7=&\,|\chi_0+\chi_{16}|^2+|\chi_4+\chi_{12}|^2+|\chi_6+\chi_{10}|^2
&\cr&\,+\chi_8\,(\chi_2+\chi_{14})^*+(\chi_2+\chi_{14})\,\chi^*_8+|\chi_8|^2
\ ,\qquad\qquad {\rm for}\ k=16&\cr
{\cal E}_8=&\,|\chi_0+\chi_{10}+\chi_{18}+\chi_{28}|^2+|\chi_6+\chi_{12}+
\chi_{16}+\chi_{22}|^2\ ,\qquad {\rm for}\ k=28\ .&\cr}$$\smallskip

The physical invariants ${\cal A}_n$ and ${\cal D}_n$ are generic, corresponding
respectively to the order 1 (i.e.\ identity) and order 2 (i.e.\ the simple-current $J$)
Dynkin diagram symmetries, as we shall see in \S 1.7. Physically,
they are the partition functions of  WZW models on SU$_2(\C)$ and SO$_3(\R)$
group manifolds, resp. The exceptionals ${\cal E}_6$ and ${\cal E}_8$
are best interpreted as due to the $C_{2,1}\supset A_{1,10}$ and $G_{2,1}
\supset A_{1,28}$ conformal embeddings (see \S 1.7; standard notation is to
write `$X_{\ell,k}$' for `$X_\ell^{(1)}$ and level $k$'). The ${\cal E}_7$ exceptional
is harder to interpret, but can be thought of as the first in an infinite
series of exceptionals involving rank-level duality and $D_4$ triality.

Around Christmas 1985, Zuber wrote Kac about the $A_1^{(1)}$
physical invariant problem, and mentioned the physical invariants he
and Itzykson  knew at that point (what we now call ${\cal A}_\star$ and ${\cal D}_{even}$).
A few weeks later, Kac wrote back saying he found one more invariant,
and jokingly pointed out that it must be indeed quite exceptional as the exponents of $E_6$
appeared in it. ``I must confess that I didn't pay much attention to that
last remark (I hardly knew what Coxeter exponents were, at the time!)''
[63]. By spring 1986, Cappelli arrived in
Paris and got things moving again; together Cappelli-Itzykson-Zuber
found  ${\cal E}_7$, ${\cal D}_{odd}$,
and then ${\cal E}_8$, and struggled to find more. ``And it is only in
August [1986], during a conversation with Pasquier, in which he was showing
me his construction of lattice models based on Dynkin diagrams, that I
suddenly remembered this cryptic but crucial!\ observation of Victor,
rushed to the library to find a list of the exponents of the other algebras...
and found with the delight that you can imagine that they were matching
our list'' [63]. Thus the A-D-E pattern to these physical invariants was
discovered.

\bigskip\bigskip\noindent{{\smcap 1.6. The A-D-E meta-pattern}}\bigskip

Before we discuss {\it meta-patterns} in math, let's introduce the
notion of {\it lattice}\footnote{$^{11}$}{{\smal There are
many words in math which have several incompatible meanings. For example,
there are {\smit vector} fields and {\smit number} fields, and modular {\smit
forms} and modular {\smit representations}. `Lattice' is another of these
words. Aside from the geometric meaning we will use, it also refers to
a `partially ordered set'.}}, a simple geometric structure we'll keep 
returning to in these
notes. The standard reference for lattice theory is [13].

 Consider the real vector
space $\R^{m,n}$: its vectors look like $\vec{x}=(\vec{x}_+;\vec{x}_-)$
where $\vec{x}_+$ and $\vec{x}_-$ are $m$- and $n$-component vectors respectively,
and dot products are given by $\vec{x}\cdot \vec{y}=
\vec{x}_+\cdot\vec{y}_+-\vec{x}_-\cdot\vec{y}_-$. The dot products
$\vec{x}_{\pm}\cdot\vec{y}_{\pm}$ are given by the usual product and
sum of components. For example, the familiar Euclidean (positive-definite)
space is $\R^n=\R^{n,0}$, while Minkowski space is  $\R^{3,1}$.

Now choose any basis $B=\{\vec{x}_1,\ldots,\vec{x}_{m+n}\}$ in $\R^{m,n}$.
So $\R^{m,n}=\R\vec{x}_1+\cdots+\R\vec{x}_{m+n}$. Define the set
$\L(B):=\Z\vec{x}_1+\cdots+\Z\vec{x}_{m+n}$. This is a {\it lattice},
and all lattices can be formed in this way\footnote{$^{12}$}{{\smal In most
presentations a lattice is permitted to have smaller dimension
than its ambient space, however that freedom gains no real generality.}}. So 
a lattice is discrete and is closed under
sums and integer multiples. For example, $\Z^{m,n}$ is a lattice
(take the standard basis in $\R^{m,n}$). A more interesting lattice
is the hexagonal lattice (also called $A_2$), given by the basis
$B=\{({\sqrt{2}\over 2},{\sqrt{6}\over 2}),({\sqrt{2}},0)\}$
 of $\R^2$ --- try to plot several points. If you wanted to slide a
bunch of coins on a table together as tightly as possible, their
centres would form this hexagonal lattice. Another important lattice
is $II_{1,1}\subset\R^{1,1}$, given by $B=\{({1\over \sqrt{2}};\,{1\over
\sqrt{2}}),\,({1\over \sqrt{2}};\,{-1\over \sqrt{2}})\}$; equivalently
it can be thought of as the set of all pairs $(a,b)\in\Z^2$ with dot product
$$(a,b)\cdot(c,d)=ad+bc\ .\eqno(1.6.1)$$

It is important to note that different choices of basis may or may not
result in a different lattice. For a trivial example, consider $B=\{1\}$
and $B'=\{-1\}$ in $\R=\R^{1,0}$: they both give the lattice $\Z=\Z^{1,0}$.
Two lattices are called {\it equivalent} if they only differ by a
change-of-basis. E.g.\ $B=\{({1\over \sqrt{2}},\,{1\over
\sqrt{2}}),\,({1\over \sqrt{2}},\,{-1\over \sqrt{2}})\}$ in $\R^2$ yields
a lattice equivalent to $\Z^2$.

The {\it dimension} of the lattice is $m+n$. The lattice is
called {\it positive-definite} if it lies in some $\R^m$ (i.e.\ $n=0$).
The lattice is called {\it integral} if all dot products $\vec{x}\cdot
\vec{y}$ are integers, for $\vec{x},\vec{y}\in\L$. A lattice $\L$ is
called {\it even} if it is integral and in addition all norms
$\vec{x}\cdot\vec{x}$ are {\it even} integers. For example, $\Z^{m,n}$
is integral but not even, while $A_2$ and $II_{1,1}$ are even. The {\it
dual} $\Lambda^*$ of a lattice $\L$ consists of all vectors $\vec{x}\in
\R^{m,n}$ such that $\vec{x}\cdot\L\subset\Z$. So a lattice is integral
iff $\L\subseteq\L^*$. A lattice is called {\it self-dual} if $\L=\L^*$.
$\Z^{m,n}$ and $II_{1,1}$ are self-dual but $A_2$ is not.

There are lots of  `meta-patterns' in math,
i.e.\ collections of seemingly different problems which have similar answers. Once
one of  these meta-patterns is identified it is always helpful to
understand what is responsible for it. For example, while I
was writing up my PhD thesis I noticed in several places the numbers
1, 2, 3, 4, and 6. For instance $\cos(2\pi r)\in\Q$ for $r\in\Q$ iff
the denominator of $r$ is 1, 2, 3, 4, or 6. This pattern was easy to explain:
they are precisely those positive integers $n$ with Euler totient $\phi(n)
\le 2$, i.e.\ there are at most 2 positive numbers less than $n$
coprime\footnote{$^{13}$}{{\smal We say {\smit m},{\smit n} are {\smit
coprime} if any prime {\smit p} which divides {\smit m} does not divide
{\smit n}, and vice versa.}}  to $n$.
The other incidences of these numbers can usually be reduced to this
$\phi(n)\le 2$ property (e.g.\ the dimension of the number field $\Q[\cos(2\pi
{a\over b})]$ (see \S1.8)
considered as a vector space over $\Q$ will be $\phi(b)/2$).

A more interesting meta-pattern involves the number 24 and its divisors. One sees 24
wherever  modular forms naturally appear. For instance, we see it
in the critical dimensions in string theory: $24+2$ and $8+2$. Another example:
the dimensions of even
self-dual positive-definite lattices must be a multiple of 8 (e.g.\
the $E_8$  root lattice defined shortly
has dimension 8, while the Leech lattice discussed in \S2.4 has
dimension 24).  The meta-pattern
24 is also understood: the fundamental problem for which it is the answer
is the following one. Fix $n$, and consider the congruence $x^2\equiv 1$
(mod $n$). Certainly in order to have a chance of satisfying this,
$x$ and $n$ must be coprime. The extreme situation is when {\it every}
number $x$ coprime to $n$ satisfies this congruence:
$${\rm gcd}(x,n)=1\qquad\Longleftrightarrow\qquad x^2\equiv 1\ ({\rm mod}\ n)\ .
\eqno(1.6.2)$$
 The reader can try to
verify the following simple fact: $n$ obeys this extreme situation (1.6.2) iff 
$n$ divides 24.

What does this congruence property have to do with these other occurrences
of 24? Let $\L$ be an even self-dual positive-definite lattice of dimension
$n$. Then an elementary argument shows that there will exist an $n$-tuple
$\vec{a}=(a_1,\ldots,a_n)$ of odd integers with the property that 8
must divide $\vec{a}\cdot\vec{a}=\sum_ia_i^2$. But $a_i^2\equiv 1$ (mod 8),
and so we get $8|n$.

A much deeper and still not-completely-understood meta-pattern is
called A-D-E (see [1] for a discussion and examples).
The name comes from the so-called {\it simply-laced algebras}, i.e.\
the simple finite-dimensional Lie algebras whose Dynkin diagrams ---
see Figure 6 in [59] ---  contain only single edges (i.e.\ no
arrows). These  are the $A_\star$- and $D_\star$-series,
along with the $E_6$, $E_7$ and $E_8$ exceptionals. The claim is that
many other problems, which don't seem to have anything directly in common
with simple Lie algebras, have a solution which falls into this A-D-E
pattern (for an object to be meaningfully labelled $X_\ell$, some of
the data  associated to the algebra $X_\ell$ should
reappear in some form in that object). Let's look at some examples.

Consider even positive-definite lattices $\L$. The smallest possible
nonzero norm in $\L$ will be 2, and the vectors of
norm 2 are special and are called {\it roots}. The reason they are special
is that reflecting through them will always be an automorphism of $\L$.
That is, the reflection $\vec{u}\mapsto
\vec{u}-2{\vec{u}\cdot\vec{\alpha}\over \vec{\alpha}\cdot \vec{\alpha}}
\,\vec{\alpha}$ through $\vec{\alpha}\ne \vec{0}$  won't in general
map $\L$ to itself, unless $\vec{\alpha}$ is a root of $\L$.
It is important in lattice theory to know the lattices which are spanned
by their roots; it turns out these are precisely the orthogonal direct sums of
lattices called $A_n$, $D_n$, and $E_6$, $E_7$ and $E_8$. They carry those
names for a number of reasons. For example, the lattice called $X_n$ will have a basis
$\{\vec{\alpha}_1,\ldots,\vec{\alpha}_n\}$ with the property that the matrix
$A_{ij}:=\vec{\alpha}_i\cdot\vec{\alpha}_j$ is the Cartan matrix (see \S 2.7)
for the Lie
algebra $X_n$! Also, the reflection group generated by reflections in
the roots of the lattice $X_n$ will be isomorphic to the Weyl group of
the Lie algebra $X_n$. Finally, to any simple Lie algebra there is
canonically associated a  lattice called the root  lattice; for the
simply-laced algebras,
these will equal the corresponding
lattice of the same name. Incidentally, the root lattices for the non-simply-laced
simple algebras will (up to rescalings) be direct sums of the simply-laced
root lattices.

We have already
met the $A_2$ lattice: it is the densest packing of circles in the plane. It has
long been believed that the obvious pyramidal way to pack oranges is also
the densest possible way --- the centres of the oranges form the $A_3$ root lattice.
A controversial proof for this famous conjecture has been offered by
W.-Y.\ Hsiang in 1991; in 1998 a new proof by Hale {\it et al} has been
proposed. The densest known packings in dimensions 4,5,6,7,8 are $D_4,D_5,
E_6,E_7,E_8$, resp.\ $E_8$ is the smallest even self-dual positive-definite
lattice.

A famous A-D-E example is called the McKay\footnote{$^{14}$}{{\smal He is the
same John McKay we will celebrate in section 2.1.}} correspondence. Consider any
finite subgroup $G$ of the Lie group SU$_2(\C)$ (i.e.\ the $2\times 2$
unitary matrices with determinant 1). For example, there is the cyclic
group $\Z_n$ of $n$ elements generated by the matrix
$$M_n=\left(\matrix{\exp[2\pi \i/n]&0\cr 0&\exp[-2\pi \i/n]\cr}\right)$$
Let $R_i$ be the irreducible representations of $G$. For instance,
for $\Z_n$, there are precisely $n$ of these,
all 1-dimensional, given by sending the generator $M_n$ to $\exp[2k\pi \i/n]$
for each $k=1,2,\ldots,n$. Now consider the tensor product $G\otimes R_i$,
where we interpret $G\subset {\rm SU}_2(\C)$ here as a 2-dimensional
 representation. We can decompose that product into
a direct sum $\oplus_j m_{ij}R_j$ of irreducibles (the $m_{ij}$ here
are multiplicities). Now create a graph with one node for each $R_i$,
and with the $i$th and $j$th nodes ($i\ne j$) connected with precisely $m_{ij}$
directed edges $i\rightarrow j$. If $m_{ij}=m_{ji}$, we agree to erase
the  double arrows
from the $m_{ij}$ edges. Then McKay observed that the graph of any $G$
will be a distinct extended Dynkin diagram of A-D-E type! For instance, the
cyclic group with $n$ elements corresponds to the extended graph of
$A_{n-1}$.

How was McKay led to his remarkable correspondence? He knew that the
sum of the `marks' $a_i=1,2,3,4,5,6,4,2,3$ associated to each node of
the extended $E_8$ Dynkin diagram equaled 30, the Coxeter number of
$E_8$. So what did their {\it squares} add to? 120, which he
recognised as the cardinality of one of the exceptional finite
subgroups of SU$_2(\C)$, and that got him thinking...

Another famous example of A-D-E, due to Arnol'd, are the `simple critical points'
of smooth complex-valued functions $f$, on e.g.\ $\C^3$. For example,
both $x^2+y^2+z^{n+1}$ and $x^2+y^3+z^5$ have singularities at $(0,0,0)$
(i.e.\ their first partial derivatives all vanish there), and they are
assigned to $A_n$ and $E_8$, respectively. The SU$_2(\C)$ subgroups
can be related to singularities as follows. The group SU$_2(\C)$ acts
on $\C^2$ in the obvious way (matrix multiplication). If $G$ is a discrete
subgroup, then consider the (ring of) polynomials in 2 variables $w_1,w_2$
invariant under $G$. It turns out it will have 3 generators $x(w_1,w_2)$,
$y(w_1,w_2)$, $z(w_1,w_2)$, which are connected by 1 polynomial relation (syzygy).
For instance, take $G$ to be the cyclic group $\Z_n$, then we're interested in
polynomials $p(w_1,w_2)$ invariant under $w_1\mapsto\exp[2\pi\i/n]w_1$,
$w_2\mapsto\exp[-2\pi\i/n]w_2$. Any such invariant $p(w_1,w_2)$ is clearly
generated by (i.e.\ can be written as a polynomial in) 
$w_1w_2$, $w_1^n$ and $w_2^n$. Choosing instead the generators $x={w_1^n-
w_2^n\over 2}$, $y=\i\,{w_1^n+w_2^n\over 2}$, $z=w_1w_2$, we get the
syzygy $z^n=-(x^2+y^2)$. For any $G$, generators $x,y,z$ can always be found
so that the syzygy will be one of the polynomials associated to a
simple singularity, and in fact will give the equation of the algebraic
surface $\C^2/G$ as a 2-dimensional complex surface in $\C^3$ (e.g.\
the complex surfaces $\C^2/\Z_n$ and $\{(x,y,z)\in\C^3\,|\,x^2+y^2+z^n=0\}$
are equivalent).

Arguably the first A-D-E classification goes back to Theaetetus, around 400 B.C.
He classified the regular solids. For instance the tetrahedron can be
associated to $E_6$ while the cube is matched with $E_7$. This A-D-E
is only partial, as there are no regular solids assigned to the A-series,
and to get the D-series one must look at `degenerate regular solids'.

The closest thing to an explanation of the A-D-E meta-pattern
would seem to be the notion of `additive
assignments' on graphs (which is a picturesque way of describing the
corresponding eigenvalue problem). Consider any graph ${\cal G}$ with
undirected edges,
and none of the edges run from a node to itself. We can also assume
without loss of generality that ${\cal G}$ is connected. Assign a positive number
$a_i$ to each node. If this assignment has the property that for each $i$,
 $2a_i=\sum a_j$ where the sum is over all nodes $j$ adjacent to $i$
 (counting multiplicities of edges), then we call it `additive'.
For instance, for the graph $\circ\!\!\!=\!\!\!\circ$, the assignment $a_1=1=a_2$
is additive, but the assignment $a_1=1,a_2=2$ is not. The question is,
which graphs have an additive assignment? The answer is: precisely the extended
Dynkin diagrams of A-D-E type! And their additive assignments are unique (up to
constant proportionality) and are given by
the {\it marks} $a_i$ of the algebra (see e.g.\ the Table on
p.54 of [38]). For example the extended $A_n$ graph 
consists of $n+1$ nodes arranged in a circle, and its marks $a_i$ all equal 1. 

What do additive assignments have to do with the other A-D-E classifications?
Consider a finite subgroup $G$ of SU$_2(\C)$. Take the dimension of
the equation $G\otimes R_i=\oplus_jm_{ij}R_j$: we get $2 d_i=\sum_j
m_{ij}d_j$ where $d_j={\rm dim}(R_j)$. Hence the dimensions of the irreducible
representations define an additive assignment for each of McKay's graphs,
and hence those graphs must be of A-D-E type (provided we know $m_{ij}=m_{ji}$).

As Cappelli-Itzykson-Zuber observed, the physical invariants for $A_1^{(1)}$ also
 realise the A-D-E pattern, in the following sense.
The {\it Coxeter number} $h$ of the name ${\cal X}_\ell$ (i.e.\ the sum
$\sum_ia_i$ of the marks) equals 
$k+2$, and the {\it exponents} $m_i$ of ${ X}_{\ell}$ equal those
$a\in P_+^k$ for which $M_{aa}\ne 0$ (for the simply-laced algebras, the $m_i$
are defined by writing the eigenvalues of the corresponding Cartan
matrix (see \S 2.7) as $4\sin^2({\pi m_i\over 2h})$ --- the $m_i$
are integers and the smallest is always 1). Probably what first led Kac to
his observation about the $E_6$ exponents  was that $k+2$
(this is how $k$ enters most formulas)
for his exceptional  equalled the Coxeter number 12 for
$E_6$. More recently, the operator algebraists Ocneanu [48] and independently
 B\"ockenhauer-Evans 
[4] found an A-D-E interpretation for the off-diagonal entries $M_{ab}$
of the $A_1^{(1)}$ physical invariants, using subfactor theory.

We are not claiming 
 that this $A_1^{(1)}$ classification is `equivalent' to any other A-D-E
 one --- that would miss the point of meta-patterns. What we really
 want to do is to identify some critical combinatorial part of an
 $A_1^{(1)}$  proof with critical parts in other A-D-E
 classifications --- this is what we did with the other
 meta-patterns. A considerably simplified proof of the $A_1^{(1)}$  
classification is
now available [29], so hopefully this task will now be easier.

There has been some progress  at understanding this $A_1^{(1)}$
A-D-E. Nahm [46] constructed the invariant ${\cal X}_\ell$ in terms of the
compact simply-connected Lie group of type $X_\ell$, and in this way could
interpret the $k+2=h$ and $M_{m_im_i}\ne 0$ coincidences.
A very general explanation for A-D-E has been suggested by
Ocneanu [48] using his theory of path algebras on graphs; although his work
 has never been published, others are now rediscovering (and publishing!)
 similar work (see e.g.\ [4]).
Nevertheless, the A-D-E in CFT remains almost as mysterious now as it did a 
dozen years ago --- for example it still isn't clear how it directly relates
to additive assignments.

There are 4 other claims for A-D-E classifications of families of 
RCFT physical invariants, and all of them inherit their (approximate)
A-D-E pattern from the more fundamental $A_1^{(1)}$ one. One 
is the $c<1$ minimal models, also proven in [8], and another is the
$N=1$ superconformal minimal models, proved by Cappelli (1987). In
both cases the
physical invariants are parametrised by pairs of A-D-E diagrams.
The list of known
$c=1$ RCFTs  also looks like A-D-E (two series parametrised by $\Q_+$,
and three exceptionals), but the completeness of that list
has never been rigourously established.

The fourth classification often quoted as A-D-E, is the $N=2$
superconformal minimal models. Their classification was done by Gannon
(1997). The  connection here with A-D-E turns out to
be rather weak: e.g.\ 20, 30, and 24 distinct invariants would 
have an equal right to be  called
${\cal E}_6$, ${\cal E}_7$, and ${\cal E}_8$ respectively. It would appear
that the frequent claims that the $N=2$ minimal models fall
into an A-D-E pattern are rather dubious.

Hanany-He [35] suggest that the $A_1^{(1)}$ A-D-E pattern can be related
to subgroups $G\subset {\rm SU}_2(\C)$ by orbifolding 4-dimensional
$N=4$ supersymmetric gauge theory by $G$, resulting in an $N=2$ superCFT
whose `matter matrix' can be read off from the Dynkin diagram corresponding
to $G$. The same game can be played with finite subgroups of SU$_3(\C)$,
resulting in $N=1$ superCFTs whose matter matrices correspond to graphs
very reminiscent of the `fusion graphs' of Di Francesco-Petkova-Zuber
(see e.g.\ [62]) corresponding to $A_2^{(1)}$ physical invariants.
[35] use this to conjecture a McKay-type correspondence between
singularities of type $\C^n/G$, for $G\subset {\rm SU}_n(\C)$, and the
physical invariants of $A_{n-1}^{(1)}$. This in their view would be the
form A-D-E takes for higher rank physical invariants.
Their actual conjecture though is still somewhat too vague.

For a final example of meta-pattern, consider `modular function' (see \S2.3).
After all, they appear in a surprising variety of places and
disguises. Maybe we shouldn't regard their ubiquity as fortuitous,
instead perhaps there's a deeper common `situation' which is the
source for that ubiquity. Just as `symmetry' yields `group', or
`rain-followed-by-heat' breeds mosquitos. Math is not above
metaphysics; like any area it grows by asking questions, and
changing your perspective --- even to a metaphysical one --- should
suggest new questions.

\bigskip\bigskip\noindent{{\smcap 1.7. Simple-currents and
charge-conjugation}}\bigskip

The key properties\footnote{$^{15}$}{{\smal A good exercise for the
reader is to prove that if {\smit S} is unitary and symmetric, and
obeys (1.7.1), then there will be at most finitely many physical invariants
{\smit M} for that {\smit S},{\smit T}.}} of the matrix $S$ are that 
it's unitary and symmetric (so $M$ in \S1.5 equals $SMS^*$),
$$S_{0\mu}>0\ {\rm for\ all}\ \mu\in P_+^k\ ,\eqno(1.7.1)$$
 and that the numbers $N_{\la\mu}^\nu$
defined by Verlinde's formula (1.4.4) are nonnegative integers. These are
obeyed by the matrix $S$ in any (unitary) RCFT. From these
basic properties, we will obtain here some elementary consequences
which have important applications.

But first, let's make an observation which isn't difficult to prove,
but doesn't appear to be generally known. 
\smallskip\centerline{\it Verlinde's formula looks
strange, but it is quite generic,}\smallskip

\noindent and we can see it throughout math and mathematical physics. Consider the following.

Let ${\frak A}$ be a commutative associative algebra, over $\R$ say.
Suppose ${\frak A}$ has a finite basis $\Phi$ (over $\R)$ containing the unit 1.
Define the `structure constants' $N_{ab}^c\in\R$, for $a,b,c\in\Phi$,
by $ab=\sum_{c\in\Phi}N_{ab}^c c$. Suppose there is an algebra homomorphism
$*$ (so $*$ is linear, and $(xy)^*=x^*y^*$) which permutes the basis vectors
(so $\Phi^*=\Phi$),
and we have the relation $N_{ab}^1=\delta_{b,a^*}$. We call any such
algebra ${\frak A}$ a {\it fusion algebra}. Then any fusion algebra will
necessarily have a unitary matrix $S$ with $S_{1a}>0$ and with the
structure constants given by Verlinde's formula. Algebraically, the relation
$S=S^t$ holds if ${\frak A}$ is `self-dual' in a certain natural sense.

Define the `fusion matrices' $N_\la$ by $(N_\la)_{\mu\nu}=N_{\la\mu}^\nu$.
Then Verlinde's formula  says that
the $\mu$th column $S_{\uparrow\mu}$ of $S$ is an
eigenvector of each fusion matrix $N_\la$, with eigenvalue ${S_{\la\mu}
\over S_{0\mu}}$.\smallskip

\noindent{{\bf Useful Fact.}} If ${v}$ is a simultaneous eigenvector 
of each fusion matrix $N_\la$, then there exists a constant $c\in\C$ and
a $\la\in P_+^k$ such that ${v}=c\,S_{\uparrow\la}$.\smallskip

For one consequence, take the complex conjugate of the eigenvector equation
$N_\la\,S_{\uparrow\mu}={S_{\la\mu}\over S_{0\mu}}\,S_{\uparrow\mu}$: we get
that the vector $S_{\uparrow\mu}^*$ is a simultaneous eigenvector of all
$N_\la$, and hence must equal $c\,S_{\uparrow\ga}$ for some number $c$ and
weight $\ga\in P_+^k$, both depending on $\mu$. Write $\ga=C\mu$; then
$C$ defines a permutation of $P_+^k$. The reader can verify
that unitarity of $S$ forces $|c|=1$, while (1.7.1) forces $c>0$. Thus
$c=1$ and we obtain the formula
$$S_{\la\mu}^*=S_{\la,C\mu}=S_{C\la,\mu}\ .\eqno(1.7.2)$$
Also, unitarity and symmetry of $S$ forces $C=S^2$, while conjugating
twice shows $C^2=id.$ $C$ is an important matrix in RCFT, and is called
{\it charge-conjugation}. When $C=id.$, then the matrix $S$ is real.

Note that (1.7.1) now implies $C0=0$. Also $CT=TC$. Hence $M=C$ will
always define a physical invariant, and if $M$ is any other physical
invariant, the matrix product $MC=CM$ will define another physical
invariant. Also, $N_{C\la,C\mu}^{C\nu}=N_{\la\mu}^\nu$ and
$(N_\la)^t=CN_{\la}=N_{\la}C=N_{C\la}$.

For the WZW (=affine) case, $C$ has a special meaning:
$C\la$ is the highest-weight `contragredient' to $\la$.
$C$ corresponds to an order 2 (or 1) symmetry of the (unextended) Dynkin
diagram. For example, for $A_\ell^{(1)}$, we have $C(\la_0,\la_1,\ldots,\la_{\ell-1},
\la_\ell)=(\la_0,\la_\ell,\la_{\ell-1},\ldots,\la_1)$. For 
$A_1^{(1)}$ then, $C=id.$, which can also be read off from (1.4.3).

The algebras $D^{(1)}_{even}$ all have at least
one nontrivial symmetry of the (unextended) Dynkin diagram which isn't
the charge-conjugation. The most interesting example is $D_4^{(1)}$, which
has 5 of these. By a {\it conjugation}, we will mean any symmetry of the
unextended Dynkin diagram.

To go much further, we need a fascinating tool called {\it
Perron-Frobenius theory} --- a collection of results concerning
the eigenvalues and eigenvectors of nonnegative matrices
(i.e.\ matrices in which every entry is a nonnegative real number).
Whenever you have such matrices in your problem, and it is natural
to multiply them, then there is a good chance
Perron-Frobenius theory will tell you something interesting. The basic
result here is that if $A$ is a nonnegative matrix, then there will
be a nonnegative eigenvector ${x}\ge 0$ with eigenvalue $\rho\ge 0$,
such that if $\la$ is any other eigenvalue of $A$, then $|\la|\le\rho$.
There are lots of other results (see e.g.\ [44]), e.g.\ $\rho$ must be at least as large as
 any diagonal entry of $A$, and there must be a row-sum of $A$ no bigger than
 $\rho$, and another row-sum no smaller than $\rho$.

For instance, consider
$$A=\left(\matrix{1&1&1\cr 1&1&1\cr 1&1&1}\right)\ ,\qquad
B=\left(\matrix{1&1&1\cr 1&0&0\cr 1&0&0}\right)\ .$$
Perron-Frobenius eigenvectors for $A$ and $B$ are $\left(\matrix{
1\cr 1\cr 1}\right)$ and $\left(\matrix{2\cr 1\cr 1}\right)$,
with eigenvalues 3 and 2 resp. The other eigenvalue of $A$
is 0 (multiplicity 2), while those of $B$ are 0 and $-1$.

Fusion matrices $N_\la$ are nonnegative, and it is indeed
natural to multiply them:
$$N_\la N_\mu=\sum_{\nu\in P_+^k}N_{\la\mu}^\nu N_\nu\ .$$
So we can expect Perron-Frobenius to tell us something interesting.
This is the case, and we obtain the curious-looking inequalities
$$S_{\la 0}\,S_{0\mu}\ge|S_{\la\mu}|\,S_{00}\ .\eqno(1.7.3)$$
Squaring both sides, summing over $\mu$ and using unitarity, we get
that $S_{\la 0}\ge S_{00}$. In other words, the ratio
${S_{\la 0}\over S_{00}}$, called the {\it quantum-dimension} 
 of $\la$, will necessarily be $\ge 1$.

The term `quantum-dimension' comes from quantum groups, where ${S_{\la
0}\over S_{00}}$ is the quantum-dimension of the module labelled by $\la$
of the quantum group $U_q(X_\ell)$.

The borderline case then \footnote{$^{16}$}{{\smal This seems to be a standard
trick in math: when some sort of bound is established, look at the
extremal cases which realise that bound. If your bound is a good one, it
should be possible to say something about those extremal cases, and having
something to say  is always of paramount importance. This trick is 
used for instance in the definition of 24 last section, and the definition of
normal subgroup in section 2.2.}}
is when a quantum-dimension {\it equals} 1.
Any such weight is called a {\it simple-current}. The theory of
simple-currents was developed most extensively by Schellekens and
collaborators (see e.g.\ [50]). The simple-currents for the affine algebras were classified
by J.\ Fuchs (1991), and the result is that (with one unimportant exception:
$E_8^{(1)}$ at level 2) they all correspond to symmetries of
the extended Dynkin diagrams. In particular, applying any such symmetry
to the vacuum $0=(k,0,\ldots,0)$ gives the list of simple-currents.
For instance, the $\ell+1$ weights of the form $(0,\ldots,k,\ldots,0)$
($k$ in the $i$th spot) are the simple-currents for $A_\ell^{(1)}$.
There are 2 simple-currents for $B_\ell^{(1)}$, $C_\ell^{(1)}$ and $E_7^{(1)}$,
3 for $E_6^{(1)}$, and 4 for $D_\ell^{(1)}$. Simple-currents play a large
role in RCFT, as we shall see.

Let $j$ be any simple-current. Then (1.7.3) becomes $S_{0\mu}\ge|S_{j\mu}|$
for all $\mu$, so  unitarity forces $S_{0\mu}=|S_{j\mu}|$, that is
 $$S_{j\mu}=\exp[2\pi \i \,Q_j(\mu)]\,S_{0\mu}\qquad\forall \mu\in P_+^k
\eqno(1.7.4)$$
for some rational numbers $0\le Q_j(\mu)<1$. Hence by diagonalising,
we get $N_jN_{Cj}=I$. But the inverse of a nonnegative  matrix $A$ is
itself nonnegative, only if $A$ is a `generalised permutation matrix',
i.e.\ a permutation matrix except the 1's can be replaced by any positive
numbers. But $N_j$ and $N_{Cj}$ are also integral, and so they must in fact
be permutation matrices. Write $(N_j)_{\la\mu}=\delta_{\mu,J\la}$ for
some permutation $J$ of $P_+^k$. So $j=J0$. Then
$$\delta_{\la,\mu}=N_{j,\la}^{J\mu}=\sum_{\nu}\exp[2\pi \i \,Q_j(\nu)]\,
S_{\la\nu}\,S_{J\mu,\nu}^*\ ,$$
so taking absolute values and using the triangle inequality and unitarity
of $S$, we find that (1.7.4) generalises:
$$S_{J\la,\mu}=\exp[2\pi \i \,Q_j(\mu)]\,S_{\la\mu}\ .\eqno(1.7.5)$$

The simple-currents form a finite abelian group, corresponding to the
composition of the permutations $J$. For any simple-currents $J,J'$,
we get the symmetry $N_{J\la,J'\mu}^\nu=N_{\la\mu}^{JJ'\nu}$.
The $\Q/\Z$-valued functions
$Q_j$ define gradings on the fusion rings, and conversely any grading
corresponds to a simple-current in this way.

For example, the simple-current $j=(0,k)$ of $A_1^{(1)}$
at level $k$ corresponds to $Q_j(\la)=\la_1/2$ and the permutation
$J\la=(\la_1,\la_0)$. We can see this directly from (1.4.3).
For $A_2^{(1)}$ level $k$, there are 2 nontrivial
simple-currents, $(0,k,0)$ and $(0,0,k)$. The first of these corresponds
to triality $Q(\la)=(\la_1+2\la_2)/3$ and $\la\mapsto(\la_2,\la_0,\la_1)$,
while the second to $(2\la_1+\la_2)/3$ and $\la\mapsto(\la_1,\la_2,\la_0)$.
Similar statements hold for all affine algebras: e.g.\
 for $B_\ell^{(1)}$ level $k$,
the nontrivial simple-current has $Q_j(\la)=\la_\ell/2$ and $J\la=
(\la_1,\la_0,\la_2,\ldots,\la_\ell)$.

One of the applications of simple-currents is that physical invariants
can be built from them in generic ways.
These physical invariants all obey the selection rule
$$M_{\la\mu}\ne 0\quad\Longrightarrow\quad \mu=J\la\quad{\rm for\ some\
simple\!-\!current}\ J=J(\la,\mu)\ .\eqno(1.7.6)$$
We will call any such physical invariant $M$ a {\it simple-current invariant}.
A special case 
is the ${\cal D}_{{k\over 2}+2}$ physical invariant for $A_1^{(1)}$ at even
level $k$. 
Up to a fairly mild assumption, all simple-current invariants have been
classified for any RCFT by Schellekens and collaborators; given that
assumption, they can all be constructed by generic methods. The basic
construction is due to Bernard [3], though it has been generalised
by others. In the WZW case, all simple-current invariants 
(except some for $D_\ell^{(1)}$)  correspond
to strings on nonsimply-connected Lie groups. 

By a {\it generic physical invariant} of $X_\ell^{(1)}$ we mean
one of the form $M=C'M'C''$ where $C',C''$ are (charge-)conjugations, and $M'$
is a simple-current invariant. In other words, $M$ is
constructed in generic ways from symmetries of the extended Dynkin diagram
of $X_\ell$. Any other $M$ are called {\it exceptional}.

All known results point to the validity of the following guess:

\medskip\noindent{{\bf Conjecture.}} Choose a simple algebra $X_\ell$. Then
for all sufficiently large $k$, all physical invariants of $X_\ell^{(1)}$
at level $k$ will be generic.\medskip

In other words, any given $X_\ell^{(1)}$ will have only finitely many
exceptionals.
For instance, for $A_1^{(1)}$ and $A_2^{(1)}$ at any $k>28$ and $k>21$
resp., all physical invariants are generic. For $C_2^{(1)}$ and $G_2^{(1)}$, $k>12$ and
$k>4$ resp.\ should work.

The richest source of exceptionals are {\it conformal embeddings}.
In some cases the affine representations $L_\la$ for some algebra
$X_\ell^{(1)}$ (necessarily at level 1) can be decomposed  into  {\it finite}
direct sums of representations of some affine subalgebra $Y^{(1)}_m$ (at
some level $k$).
In this case, a physical invariant for $X_\ell^{(1)}$ level 1 will yield a
physical invariant for $Y_m^{(1)}$ level $k$, obtained by replacing every
$X_\ell^{(1)}$ level 1 character $\chi_\la$ by the appropriate finite sum
of $Y_m^{(1)}$ level $k$ characters. An example will demonstrate this
simple idea: $A_1^{(1)}$ level 28 is a conformal subalgebra of $G_2^{(1)}$ level
1, and we have the character decompositions
$$\eqalignno{\chi_{(1,0,0)}=&\,\chi'_{(28,0)}+\chi'_{(18,10)}+\chi'_{(10,18)}+\chi'_{(0,28)}&\cr
\chi_{(0,0,1)}=&\,\chi'_{(22,6)}+\chi'_{(16,12)}+\chi'_{(12,16)}+\chi'_{(6,22)}\ .&\cr}$$
Thus the unique level 1 $G_2^{(1)}$ physical invariant $|\chi_0|^2+|\chi_{(0,0,1)}|^2$
yields what we call the ${\cal E}_8$ physical invariant of $A_1^{(1)}$.
All level 1 physical invariants are known, as are all conformal embeddings
and the corresponding character decompositions (branching rules).

\bigskip\bigskip\noindent{{\smcap 1.8. Galois Theory}}\bigskip

Evariste Galois was a brilliantly original French mathematician. Born shortly
before Napoleon's ill-fated invasion of Russia, he died shortly before the
ill-fated 1832 uprising in Paris. His last words: ``Don't cry, I need all my courage
to die at 20''.

Galois grew up in a time and place confused and excited by revolution.
He was known to say ``if I were only sure that a body would be enough
to incite the people to revolt, I would offer mine''.
On May 2 1832, after frustration over failure
in love and failure to convince the Paris math establishment of the
depth of his ideas, he made his decision. A duel was arranged with a friend,
but only his friend's gun would be loaded. Galois died the day after a bullet
perforated his intestine. At his funeral it was discovered that a famous general
 had also just died, and the revolutionaries decided to use the
general's death  rather than Galois'
as a pretext for an armed uprising. A few days later the streets of Paris
were blocked by barricades, but not because of Galois' sacrifice: his death
had been pointless [56].

Galois theory in its most general  form is the study of relations between
objects defined implicitly by some conditions. For example, the objects
could be the solutions to a given differential equation. In the incarnation
of Galois we are interested in here, the objects are numbers, namely the
zeros of certain polynomials. We will sketch this theory below, but
see e.g.\ the article by Stark in [60] for more details.

Gauss seems to have been the first to show that `weird' (complex)
numbers could tell us about the integers.
For instance, suppose we are interested in the equation $n=a^2+b^2$.
Consider $5=2^2+1^2$. We can write this as $5=(2+\i)(2-\i)$, so we are
led to consider complex numbers of the form $a+b\i$, for $a,b\in\Z$.
These are now called `Gaussian integers'. Suppose we know the following
theorem:

\medskip\noindent{{\bf Fact.}} Let $p\in\Z$ be any prime number. Then
$p$ factorises over the Gaussian integers iff $p=2$ or $p\equiv 1$ (mod 4).
\medskip

By `factorise' there, we mean $p=zw$ where neither $z$ nor $w$ is
a `unit': $\pm 1$, $\pm\i$.

Now suppose $p$ is a prime, $=2$ or $\equiv 1$ (mod 4), and we write 
$p=(a+b\i)(c+d\i)$.
Then $p^2=(a^2+b^2)(c^2+d^2)$, so $a^2+b^2=c^2+d^2=p$. Conversely,
suppose $p=a^2+b^2$, then $p=(a+b\i)(a-b\i)$. Thus:

\medskip\noindent{{\bf Consequence.}}
\footnote{$^{17}$}{{\smal This result was first stated by Fermat in one of his infamous margin notes
(another is discussed in Section 2.3), and was finally proved a century later
by Euler. A remarkable 1-line proof was found  by Zagier [61].}}
 Let $p\in\Z$ be any prime number. Then
$p=a^2+b^2$ for $a,b\in\Z$ iff $p=2$ or $p\equiv 1$ (mod 4).
\medskip

Now we can answer the question: can a given $n$ be written as a sum of 2
squares $n=a^2+b^2$? Write out the prime decomposition $n=\prod p^{a_p}$.
Then $n=a^2+b^2$ has a solution iff $a_p$ is {\it even} for every
$p\equiv 3$ (mod 4). For instance $60=2^2\cdot 3^1\cdot 5^1$ cannot be
written as the sum of 2 squares, but $90=2^1\cdot 3^2\cdot 5^1$  can.
We can also find (and count) all solutions: e.g.\ $90=2\cdot 3^2\cdot 5
=\{(1+\i)3(1+2\i)\}\{(1-\i)3(1-2\i)\}$, giving $90=(-3)^2+9^2$.

This problem should give the reader a small appreciation for the power
of using nonintegers to study integers. Nonintegers often lurk in the shadows,
secretly watching their more arrogant brethren the integers strut.
 One of the consequences of their presence can be the existence of certain
`Galois' symmetries. Such  happens in RCFT, as we will show below.

Look at complex conjugation: $(wz)^*=w^*z^*$ and $(w+z)^*=w^*+z^*$.
Also, $r^*=r$ for any $r\in\R$. So we can say that $*$ is a structure-preserving
map $\C\rightarrow\C$ (called an {\it automorphism} of $\C$)
fixing the reals. We will write this $*\in{\rm Gal}(\C/\R)$. `Gal$(\C/\R)$'
is the Galois group of $\C$ over $\R$; it turns out
to contain only $*$ and the identity.

A way of thinking about the automorphism $*$ is that it says that, as far as the real
numbers are concerned, $\i$ and $-\i$ are identical twins.

Let ${\Bbb F}$ be any {\it field} containing $\Q$
(we defined `field' in \S1.2). The {\it Galois group}
Gal$({\Bbb F}/\Q)$ then will be the set of all automorphisms=symmetries
of ${\Bbb F}$ which fix all rationals.

For example, take ${\Bbb F}$ to be the set of all numbers of the form
$a+b\sqrt{5}$, where $a,b\in\Q$. Then ${\Bbb F}$ will be a field, which
is commonly denoted $\Q[\sqrt{5}]$ because it is generated by $\Q$ and
$\sqrt{5}$. Let's try to find its Galois group. Let $\si\in{\rm Gal}({\Bbb
F}/\Q)$. Then $\si(a+b\sqrt{5})=\si(a)+\si(b)\si(\sqrt{5})=a+b\si(\sqrt{5})$,
so once we know what $\si$ does to $\sqrt{5}$, we know everything about
$\si$. But $5=\si(5)=\si(\sqrt{5}^2)=(\si(\sqrt{5}))^2$, so
$\si(\sqrt{5})=\pm\sqrt{5}$ and there are precisely 2 possible Galois
automorphisms here (one is the identity). As far as $\Q$ is concerned, $\pm\sqrt{5}$
are interchangeable: it cannot see the difference.

For a more important example, consider the {\it cyclotomic field}
${\Bbb F}=\Q[\xi_n]$, where $\xi_n:=\exp[2\pi \i/n]$ is an $n$th root of 1.
So $\Q[\xi_n]$ consists of all complex numbers which can be expressed
as polynomials $a_m\xi_n^m+a_{m-1}\xi_n^{m-1}+\cdots+a_0$ in $\xi_n$ with
rational coefficients $a_i$. Once again, to find the Galois group
Gal$(\Q[\xi_n]/\Q)$, it is enough to see what an automorphism $\si$ does
to the generator $\xi_n$. Since $\xi_n^n=1$, we see that it must send it to another
$n$th root of 1, $\xi_n^\ell$ say; in fact it is easy to see that $\si(\xi_n)$
must be another `primitive' $n$th root of 1, i.e.\ $\ell$ must be coprime
to $n$. So Gal($\Q[\xi_n]/\Q)$ will be isomorphic to the multiplicative
group $\Z_n^\times$ of numbers between 1 and $n$ coprime to $n$.
The rationals can't see any difference between the primitive $n$th roots
of 1 --- for instance $\Q$ can't tell that $\xi_n^{\pm 1}$ are `closer
to 1'  than the other primitive roots.
So any $\si\in\,$Gal($\Q[\xi_n]/\Q)$ will correspond to some $\ell\in\Z_n^\times$,
and to see what $\si$ does to some $z\in\Q[\xi_n]$ what we do is write
$z$ as a polynomial in $\xi_n$ and then replace each occurrence of $\xi_n$
with $\xi_n^\ell$.
For example,
$$\si\bigl(\cos(2\pi a/n)\bigr)=\si\left({\xi_n^a+\xi_n^{-a}\over 2}\right)
={\xi_n^{a\ell}+\xi_n^{-a\ell}\over 2}=\cos(2\pi a\ell/n)\ .$$

So to summarise, Galois automorphisms are a massive generalisation
of the idea of complex conjugation. If in your problem complex conjugation
seems interesting, then there
is a good chance more general Galois automorphisms will play an interesting
role. This is what happens in RCFT, as we now show.

\medskip\noindent{{\bf Fact.}} [14]
Suppose $S$ is unitary and symmetric and each $S_{0a}>0$.

\item{(a)} If in addition the numbers
$N_{ab}^c$ given by Verlinde's formula (1.4.4) are rational, then the
entries $S_{ab}$ of $S$ must lie in a cyclotomic field.

\item{(b)} The numbers $N_{ab}^c$ will be rational iff for any $\si\in{\rm
Gal}(\Q[S]/\Q)$,
 there is a permutation $a\mapsto
a^\si$, and a choice
of signs $\eps_\si(a)\in\{\pm 1\}$, such that
$$\si(S_{ab})=\eps_\si(a)\,S_{a^\si,b}=\eps_\si(b)\,S_{a,b^\si}\ .\eqno(1.8.1)$$
\smallskip

`$\Q[S]$' in part (b) denotes the field generated by $\Q$ and all matrix
entries $S_{ab}$. The argument follows the one given for the charge-conjugation
$C$ at the beginning of the last section. The kinds of complex numbers
which lie in cyclotomic fields are $\sin(\pi r)$, $\cos(\pi r)$, $\sqrt{r}$
and $r\i$ for any $r\in\Q$. Almost all complex numbers fail to lie in any
cyclotomic field: e.g.\ generic cube roots, 4th roots, ..., of rationals,
as well as transcendental numbers like $e$, $\pi$ and $e^\pi$.

Of course the affine algebras satisfy the conditions of the Fact, as
does more generally the modular matrix $S$ for any unitary RCFT, and so
these will possess the Galois action. For the affine algebras this action has a
geometric interpretation in terms of multiplying weights by an integer
$\ell$ and applying Weyl group elements --- see [14] for a description.

This Fact is useful in both directions: as a way of testing whether a
conjectured matrix $S$ has a chance of producing the integral fusions we
want it to yield; and more importantly as a source of a symmetry of the RCFT which
generalises charge-conjugation. Any statement about charge-conjugation
seems to have an analogue for any of these Galois symmetries, although
it is usually more complicated.

As an example, consider $A_1^{(1)}$: (1.4.3) shows explicitly that $S_{\la\mu}$
lies in the cyclotomic field $\Q[\xi_{4(k+2)}]$.
Write $\{x\}$ for the number congruent to $x$ mod $2(k+2)$
satisfying $0\le \{x\}<2(k+2)$. Choose any Galois automorphism $\si$, and let
$\ell\in\Z_{4(k+2)}^\times$ be the corresponding integer. Then if 
$\{\ell (a+1)\}<k+2$, we will have $a^\si=\{\ell (a+1)
\}-1$, while if $\{\ell (a+1)\}>k+2$, we'll have $a^\si=2(k+2)-\{\ell (a+1)\}-1$.
The sign $\eps_\si(a)$ will depend on a contribution from $\sqrt{{2\over k+2}}$
(which for most purposes can be ignored), as well as the sign $+1$ or $-1$,
resp., depending on whether or not $\{\ell(a+1)\}<k+2$.

Consider specifically $k=10$, and
the Galois automorphism $\si_5$ corresponding to $\ell=5$. Then
the permutation is
$0\leftrightarrow(6,4)$, $(9,1)\leftrightarrow(1,9)$, $(8,2)\leftrightarrow
(2,8)$, $(4,6)\leftrightarrow(0,10)$, while $(7,3)$ and (5,5) are fixed.

This Galois symmetry has been used to find certain exceptional physical
invariants, but its greatest use so far is as a powerful selection rule
we will describe next section.

\vfill\eject\noindent{{\smcap 1.9. The modern approach to classifying
physical invariants}} \bigskip

In this final section we include some of the basic tools belonging
to the `modern' classifications of physical invariants, and we give a flavour
of their proofs.
We will state them for the $A_1^{(1)}$ level $k$ problem given above, but
everything generalises without effort. See [29] and references
therein for more details. Recall the matrices $S,T$ in (1.4.3).

First note that commutation of $M$ with $T$ implies the selection rule
$$M_{\la\mu}\ne 0\quad\Longrightarrow\quad(\la_1+1)^2\equiv(\mu_1+1)^2\
\quad({\rm mod}\ 4(k+2))\ .\eqno(1.9.1)$$
It is much harder to squeeze information out of the commutation with $S$,
but the resulting information turns out surprisingly to be much more
useful. In fact, commutation with $S$ is almost incompatible
with the constraint $M_{\la\mu}\in\Z_{\ge}$.

Note that the vacuum $0\in P_+^k$ is both physically and mathematically
special; our strategy will be to find all possible 0th rows and columns
of $M$, and then for each of these possibilities to find the remaining
entries of $M$.

The easiest result follows by evaluating $MS=SM$ at $(0,\la)$ for any $\la
\in P_+^k$:
$$\sum_{\mu\in P_+^k}M_{0\mu}\,S_{\mu\la}\ge 0\ ,\eqno(1.9.2)$$
with equality iff the $\la$th column of $M$ is identically 0. (1.9.2)
has two uses: it severely constrains the values of $M_{0\mu}$ (similarly 
$M_{\mu 0}$), and it says precisely which columns (and rows) are nonzero.

Next, let's apply the triangle inequality to sums involving (1.7.5).
Choose any $i,j\in\{0,1\}$. Then
$$M_{J^i0,J^j0}=\sum_{\la,\mu}(-1)^{\la_1 i}\,S_{0\la}\,M_{\la\mu}\,(-1)^{\mu_1
j}\,S_{0\mu}\ .$$
Taking absolute values, we obtain
$$M_{J^i0,J^j0}\le\sum_{\la,\mu}S_{0\la}\,M_{\la\mu}\,S_{0\mu}=M_{00}=1\ .$$
Thus $M_{J^i0,J^j0}$ can equal only 0 or 1. If it equals 1, then we obtain
the selection rule: 
$$\la_1i\equiv \mu_1j\ ({\rm mod}\ 2)\ {\rm whenever}\ M_{\la\mu}\ne
0\ ;$$
this implies the symmetry $M_{J^i\la,J^j\mu}=M_{\la\mu}$ for all $\la,\mu\in P_+^k$.
We can see both of these in the list of physical invariants for $A_1^{(1)}$
level $k$. This explains a lot of the properties of those invariants.
For instance, try to use this selection rule to explain why no
$\chi_{odd}$ appears in the exceptional called ${\cal E}_8$.

Our
$M$ is nonnegative, and although multiplying $M$'s may not give us back a
physical invariant,  it will give us a matrix commuting with
$S$ and $T$. In other words, the commutant is much more than merely a
vector space, it is in fact an algebra. Thus we should expect Perron-Frobenius
to tell us something here. A first application is the following.

 {Suppose $M_{\la 0}=\delta_{\la,0}$ --- i.e.\ the 0th column of $M$ is all
zeros except for $M_{00}=1$. Then Perron-Frobenius implies (with a little
work) that $M$ will be a permutation matrix --- i.e.\
there is some permutation $\pi$ of $P_+^k$ such that $M_{\la\mu}=\delta_{\mu,
\pi \la}$, and $S_{\pi \la,\pi \mu}=S_{\la\mu}$.}
This nice fact applies directly to the ${\cal A}_{\star}$ and
${\cal D}_{odd}$ physical invariants of $A_1^{(1)}$.

This is proved by studying the powers $(M^t\,
M)^L$ as $L$ goes to infinity: its diagonal entries will grow exponentially
with $L$, unless there is at most one nonzero entry on each row of $M$,
and that entry equals 1.

More careful reasoning along those lines tells us about the other
generic situation here. Namely, suppose $M_{\la 0}\ne 0$ only 
for $\la=0$ and $\la=J0$, and similarly for
$M_{0\la}$ --- i.e.\ the 0th row and column of $M$ are all
zeros except for $M_{J^i0,J^j0}=1$. Then the $\la$th row (or column) of $M$
will be identically 0 iff $\la_1$ is odd. Moreover, let $\la,\mu$ be
any non-fixed-points of $J$, and suppose $M_{\la\mu}\ne 0$. Then
$$M_{\la\nu}=\left\{\matrix{1&{\rm if}\ \nu=\mu\ {\rm or}\ \nu=J\mu\cr 0&{\rm otherwise}
\cr}\right.$$
with a similar formula for $M_{\nu\mu}$. This applies to the ${\cal D}_{
even}$ and ${\cal E}_{7}$ invariants of $A_1^{(1)}$.

Our final ingredient is the Galois symmetry (1.8.1) obeyed by $S$. Choose
any Galois automorphism $\si$. It will correspond to some
 integer $\ell$ coprime to $2(k+2)$. From (1.8.1) and $M=SMS^{*}$
 we get, for all $\la,\mu$, the important relation
$$M_{\la\mu}=\epsilon_\si(\la)\,\epsilon_\si(\mu)\,M_{\la^\si,\mu^\si}\ .
\eqno(1.9.3)$$
 From (1.9.3) and the positivity of $M$, we obtain the powerful 
{\it Galois selection rule}
$$M_{\la\mu}\ne 0\quad\Longrightarrow\quad\epsilon_\si(\la)=
\epsilon_\si(\mu)\ .\eqno(1.9.4)$$

Next let us quickly sketch how these tools are used to obtain the $A_1^{(1)}$
classification. For details the reader should consult [29].

The first step will be to find all possible values of $\la$ such that
$M_{0\la}\ne 0$ or $M_{\la 0}\ne 0$. These $\la$ are severely constrained.
We know two generic possibilities: $\la_1=0$ (good for all $k$), and $\la=J0$
(good when ${k\over 2}$ is even). We now ask the question, what other 
possibilities for $\la$ are there? Our goal is to prove (1.9.7). Assume $\la\ne
0, J0$, and write $a=\la_1+1$ and $n=k+2$.

There are only two constraints on $\la$ which we will need. One is (1.9.1):
$$(a-1)\,(a+1)\equiv 0\ ({\rm mod}\ 4n)\ .\eqno(1.9.5)$$
More useful is the Galois selection rule (1.9.4), which we 
can write as $\sin(\pi\ell {a\over n})\sin(\pi\ell{1\over n})>0$, for
all those $\ell$. But a product
of sines can be rewritten as a difference of cosines, so 
$$\cos(\pi\,\ell\,{a-1\over n})>\cos(\pi\,\ell\,{a+1\over n})\ .\eqno(1.9.6)$$
(1.9.6)  is strong and easy to solve; the reader should try to find her
own argument.

What we get is that, provided $n\ne 12,30$, $M$ obeys the
strong condition 
$$M_{\la 0}\ne 0\ {\rm or}\ M_{0\la}\ne 0\qquad\Longrightarrow\qquad \la\in\{
0,J0\}\ .\eqno(1.9.7)$$
Consider first {\bf case 1}: $M_{\la 0}=\delta_{\la,0}$. 
From above, we know $M_{\la\mu}= \delta_{
\mu,\pi \la}$ for some permutation $\pi$ of $P_+^k$ obeying $S_{\la\mu}=
S_{\pi \la,\pi \mu}$.
We know $\pi 0=0$;  put $\mu:=\pi (k-1,1)$. Then $\sin(\pi {2\over n})=\sin(\pi\,
{\mu_1+1\over n})$, and so we get either $\mu=(k-1,1)$ or $\mu=J(k-1,1)$.
By $T$-invariance (1.9.1),
the second possibility can only occur if $4\equiv (n-2)^2$ (mod $4n$), i.e.\
4 divides $n$. But for those $n$,  ${\cal D}_{{n\over 2}+1}$ is also a permutation
matrix, so replacing $M$ if necessary with the matrix product $M\,{\cal 
D}_{{n\over 2}+1}$, we can always require $\mu=(k-1,1)$, i.e.\ $\pi$ also
fixes $(k-1,1)$. It is now easy to show $\pi$ must fix any $\la$,
 i.e.\ that $M$ is the identity matrix ${\cal A}_{n-1}$.

The other possibility, {\bf case 2}, is that both $M_{0,J0}\ne 0$ and $M_{J0,
0}\ne 0$. (1.9.1) says $1\equiv (n-1)^2$ (mod $4n$),
i.e.\ ${n\over 2}$ is odd. The argument here is similar to that of {\bf
case 1}, but with $(k-2,2)$ playing the role of $(k-1,1)$. We can show that
$M_{(k-2,2),(k-2,2)}\ne 0$, except possibly for $k=16$, where we find the
exceptional ${\cal E}_7$. Otherwise we get $M ={\cal D}_{{n\over 2}+1}$.

For more general $X_\ell^{(1)}$ level $k$, the approach is

\smallskip\item{(i)} to look at all the constraints on the
$\la\in P_+^k$ for which $M_{0\la}\ne 0$ or $M_{\la 0}\ne 0$. Most
important here are $TM=MT$ (which will always be some sort of norm selection
rule) and the Galois selection rule (1.9.4). Generically, what we will find
is that such a $\la$ must equal $J0$ for some simple-current $J$, as
in (1.9.7) for $A_1^{(1)}$.

\smallskip\item{(ii)} Solve this generic case (in the $A_1^{(1)}$
classification, these were the physical invariants ${\cal A}_\star$, 
${\cal D}_\star$ and ${\cal E}_7$).

\smallskip\item{(iii)} Solve the nongeneric case. The worst of these are
the orthogonal algebras at $k=2$, as well as the places where conformal
embeddings (see \S1.7) occur.\smallskip

(ii) has recently been completed for all simple $X_\ell$, as has the
$k=2$ part of (iii). (i) is the main
remaining task in the physical invariant classification for simple
$X_\ell$.

A natural question to ask is whether A-D-E has been observed in e.g.\
the $A_2^{(1)}$ classification. The answer is no, although the fusion
graph theory of Di Francesco-Petkova-Zuber [62] is an attempt to
assign to these physical invariants graphs reminiscent of the A-D-E
Dynkin diagrams. Also, there is related
work trying to understand the $A_2^{(1)}$ classification in terms
of subgroups of SU$_3(\C)$ (as opposed to SU$_2(\C)$ for $A_1^{(1)}$) ---
see e.g.\ [35]. Finding the $A_3^{(1)}$, $A_4^{(1)}$,... classifications
would permit the clarification and testing of this vaguely conjectured relation
between the $A_n^{(1)}$ physical invariants, and singularities $\C^{n+1}
/G$ for $G$ a finite subgroup of SU$_{n+1}(\C)$.

However, a few years ago Philippe Ruelle
 was walking in a library in Dublin. He spotted a
yellow book in the math section, called {\it Complex Multiplication}
by  Lang. A strange title for a book by Lang! After all, there can't
be all that much even Lang could really say about complex multiplication! Ruelle flipped
it to a random page, which turned out to be p.26. On there he found
what we would call the Galois selection rule for $A_2^{(1)}$, analysed
and solved for the cases where $k+3$ is coprime to 6. Lang however didn't
know about physical invariants; he was reporting on work by Koblitz and
Rohrlich on decomposing the Jacobians of the Fermat curve $x^n+y^n=z^n$
into their prime pieces, called `simple factors' in algebraic
geometry. $n$ here corresponds to $k+3$. Similarly,  Itzykson discovered
traces of the $A_2^{(1)}$ exceptionals --- these occur when $k+3=8,12,24$ ---
in the Jacobian of $x^{24}+y^{24}
=z^{24}$. See [2] for further observations along these lines.
These `coincidences' are still far from
understood. Nor is it known if, more generally, the $A_\ell^{(1)}$ 
level $k$ classification
will somehow be related to the hypersurface $x_1^n+\cdots+x_\ell^n=z^n$,
for $n=k+\ell+1$.

The $(u(1)\oplus\cdots\oplus u(1))^{(1)}$ classification
has connections to rational points on Grassmannians. The Grassmannian
is (essentially) the moduli space for the Narain compactifications of the
(classical) lattice string. It would be very
interesting
to interpret other large families of physical invariants as special
points on other moduli spaces.

These new connections
relating various physical invariant classifications to other areas of
math seem to indicate that although the physical invariant classifications
are difficult, they could be well worth the effort and be of interest
outside RCFT.  Once the physical invariant
lists are obtained, we will still have the fascinating task of explaining
and developing all these mysterious connections. These thoughts keep
me going!

Another motivation for completing these lists comes from their relation
to subfactor theory
in von Neumann algebras\footnote{$^{18}$}{{\smal For reasons of necessity,
in the following discussion I'll take more liberties than usual
in the presentation.}}.
These algebras (see e.g.\ [22]) can be thought of
as symmetries of a (generally infinite) group. Their building blocks
are called {\it factors}. Jones initiated the combinatorial study of
{\it subfactors} $N$ of $M$ (i.e.\ inclusions $N\subseteq M$ where $M,N$
are factors), relating it to e.g.\ knots, and for this won a
Fields medal in 1990. Jones assigned to each subfactor $N\subseteq M$
a numerical invariant called an `index', a sort of (generally
irrational) ratio of dimensions.
Graphs (called principal and dual principal) are also associated
to subfactors. A much more refined subfactor invariant, called a
`paragroup',  has been introduced
by Ocneanu. It is essentially equivalent to a
(2+1)-dimensional topological field theory. Moreover, any RCFT can be
assigned a paragroup, and any paragroup (via a process called asymptotic
inclusion which is akin to Drinfeld's quantum doubling of Hopf algebras)
yields an RCFT. See [22] for details.

B\"ockenhauer-Evans [4] have recently developed this much further, and have
clarified the fusion graph $\leftrightarrow$ physical invariant relation.
The fusion graphs 
will correspond to subfactor principal graphs. In the work of Di
Francesco-Petkova-Zuber, that relation seems to be only empirical
(i.e.\ nonconceptual).
 
Subfactor theory
together with singularity theory is our best hope at present for
understanding and generalising the A-D-E meta-pattern.

\vfill\eject
\centerline{{\bf \huge Part 2. Monstrous Moonshine}}\bigskip

\noindent{{\smcap 2.1. Introduction}}\bigskip 

In 1978, John McKay made a very curious observation. One of the
well-known\footnote{$^{19}$}{{\smal `Well-known' is math euphemism for
`a basic result of which until recently we were
utterly ignorant.' As Conway later said, ``the {\smit j}-function was
`well-known' to other people, but not `well-known' to me.''}} functions
of classical number theory is the $j$-function\footnote{$^{20}$}{{\smal
This and other technical terms used in this introduction will be carefully
explained in the following subsections. This section is merely offered
as a quick overview.}}, given by
$$\eqalignno{j(\tau):=&\,{(1+240\sum_{n=1}^\infty\sigma_3(n)\,q^n)^3\over q\prod_{n=1}^\infty
(1-q^n)^{24}}={\Theta_{E_8}(\tau)^3\over\eta(\tau)^{24}}&\cr
=&\, q^{-1}+744+196\,884\,q+21\,493\,760\,q^2+864\,299\,970\,q^3+\cdots&(2.1.1)
\cr}$$
Here as elsewhere in this paper, $q=\exp[2\pi\i\,\tau]$. Also, $\sigma_3(n)=\sum_{d|n}
d^3$, $\Theta_{E_8}$ is the theta function of the $E_8$ root lattice,
and $\eta$ is the Dedekind eta. What is important here are the values of
the first few coefficients. What McKay noticed was that $196\,884\approx
196\,883$. Closer inspection shows $21\,493\,760\approx 21\,296\,876$, and $864\,299\,970\approx
842\,609\,326$. In fact,
$$\eqalignno{196\,884=&\,196\,883+1&(2.1.2a)\cr
21\,493\,760=&\, 21\,296\,876+196\,883+1&(2.1.2b)\cr
864\,299\,970=&\,842\,609\,326+ 21\,296\,876+2\cdot 196\,883+2\cdot 1&(2.1.2c)
\cr}$$
The numbers on the right-side are the dimensions
of the smallest irreducible representations of the Monster finite simple
group ${\Bbb M}$ (in 1978 it still wasn't certain that ${\Bbb M}$ even existed
so back then these numbers were merely conjectural).
The same game could be played with other coefficients of the $j$-function.
With numbers so large, it seemed to him doubtful that this numerology was
merely a coincidence. On the other hand,  it was hard to
imagine any deep conceptual connection between
the Monster and the $j$-function: they seem completely unrelated.

In November 1978 he mailed the `McKay equation' (2.1.2a) to
John  Thompson. At first Thompson
dismissed this as nonsense, but after checking
the next few coefficients he became convinced. He then added a
vital piece to the puzzle. It should be well-known that when one sees a
nonnegative integer, it often helps to try to interpret it as the dimension
of some vector space. Essentially, that is what McKay was proposing here.
(2.1.2) are really hinting that there is a `graded' representation $V$ of $\M$:
$$V=V_{-1}\oplus V_1\oplus V_2\oplus V_3\oplus\cdots$$
where $V_{-1}=\rho_0$, $V_{1}=\rho_1\oplus\rho_0$, $V_2=\rho_2\oplus\rho_1
\oplus\rho_0$, $V_3=\rho_3\oplus\rho_2\oplus\rho_1\oplus\rho_1\oplus\rho_0
\oplus\rho_0$, etc, where $\rho_i$ are the irreducible representations
of $\M$ (ordered by dimension), and that
$$j(\tau)-744={\rm dim}_q(V):={\rm dim}(V_{-1})\,q^{-1}+\sum_{i=1}^\infty {\rm dim}(V_i)
\,q^i\ ,\eqno(2.1.3)$$
the graded dimension of $V$.

Thompson suggested that we twist dim$_q(V)$, i.e.\ that more generally we
consider the series (now called the {\it McKay-Thompson series})
$$T_g(\tau):= 
{\rm ch}_{V,q}(g)={\rm ch}_{V_{-1}}(g)\,q^{-1}+\sum_{i=1}^\infty {\rm ch}_{V_i
}(g)\,q^i\ ,\eqno(2.1.4)$$
for each element $g\in \M$. The point is that, for any group representation
$\rho$, the character value ${\rm ch}_\rho
(id.)$ equals the dimension of $\rho$, and so
$T_{id.}(\tau)=j(\tau)-744$ and we recover (2.1.2) as special cases. But
there are many other possible choices of $g\in\M$. Thompson couldn't guess
what these functions $T_g$ would be, but he suggested that they too might
be interesting. This is a nice thought: when we see a positive integer,
we should try to interpret it as a dimension of a vector space; if there
is a symmetry present, then it may act on the vector space --- i.e.\ our
vector space may carry a representation of that symmetry group --- in which
case we can apply the Thompson trick and see what if any significance the
other character values have in our context.

Conway and Norton [12] did precisely what Thompson asked. Conway called it
``one of the most exciting moments in my life'' [11] when he opened Jacobi's
foundational (but 150 year old!) book on elliptic and
modular functions and found that the first few terms of the McKay-Thompson
series agreed perfectly with the first few terms of certain special
functions, namely the Hauptmoduls of various
genus 0 modular groups.   Monstrous Moonshine was officially born.

The word `moonshine' here is English slang for `unsubstantial or unreal'.
 It was chosen by Conway to convey as well the feeling that 
things here are dimly lit, and that Conway-Norton were `distilling
 information illegally' from the Monster character table.

In fact the first incarnation of Moonshine goes back to
Andrew Ogg in 1975. He was in France describing his result that the
primes $p$ for which the group $\Gamma_0(p)+$ has genus 0, are $\{2,3,5,11,
13, 17,19,23,29,31,41,47,59,71\}$. $\Gamma_0(p)+$ is the group generated by
$\left(\matrix{0&1\cr -p&0}\right)$ and $\Gamma_0(p)$, and is the
normaliser of $\Gamma_0(p)$ in SL$_2(\R)$ (this sentence will make a little
more sense after \S 2.3, but it isn't important here to understand it).
He also attended a lecture
by Jacques Tits, who was describing a newly conjectured simple group.
When Tits wrote down the prime decomposition of the order of that group
(see (2.2.1) below), Ogg noticed its prime factors precisely equalled his
list of primes. Presumably as a joke, he offered a bottle of Jack Daniels'
whisky to the first person to explain the coincidence.

The next step was accomplished by Griess in 1980, with the construction
of the Monster\footnote{$^{21}$}{{\smal Griess also came up with the symbol
for the Monster; Conway came up with the name.}} $\M$, and with it the proof that the
conjectured character table for
$\M$ was correct. Griess did this by explicitly constructing the 196883-dimensional
representation $\rho_1$; it turns out to have a (commutative nonassociative)
 algebra structure, now
called the {\it Griess algebra}. Though this paper was clearly
important, the construction was artificial and 100 pages long: since
the Monster is presumably a natural mathematical object (see \S 2.2),
 an elegant construction for it should exist. This was ultimately
accomplished in the mid 1980s with the construction by Frenkel-Lepowsky-Meurman
[23] of the Moonshine module $V^{\natural}$ and its interpretation by
Borcherds as a {\it vertex operator algebra}. The Griess algebra appears
naturally in $V^{\natural}$, as we shall see. $V^{\natural}$ does indeed seem
 to be a `natural' mathematical structure, and $\M$ is its automorphism
group: in fact $V^{\natural}$ is  the graded representation $V$
of $\M$ conjectured by McKay and Thompson.

Connections with physics (CFT) go back to
Dixon-Ginsparg-Harvey [19] in 1988, in a paper titled ``Beauty and the
beast: Superconformal symmetry in a Monster module''. The Moonshine module $V^{\natural}$
can be interpreted as the string theory for a $\Z_2$-orbifold of free bosons compactified on
the torus $\R^{24}/\L_{24}$ ($\L_{24}$ is the Leech lattice). Many aspects
of Moonshine make complete sense within CFT, but some
(e.g.\ the genus zero property) remain more obscure.
(Though in 1987 Moore speculated
that the 0-genus of $\Gamma_0(a)+$ could be related to the vanishing of the
cosmological constant in certain string theories related to $\M$, and
Tuite [57] related genus-zero with the conjectured uniqueness of $V^\natural$.)
 Nevertheless this helps
make the words of Dyson ring prophetic: ``I have a sneaking hope, a hope
unsupported by any facts or any evidence, that sometime in the twenty-first
century physicists will stumble  upon the Monster group, built in some
unsuspected way into the structure of the universe'' [21].

Finally, in 1992 Borcherds [5] completed the proof of the Conway-Norton
conjectures by showing $V^\natural$ is the desired representation $V$. 
The full conceptual relationship between the Monster and the
Hauptmoduls (like $j$) seems to remain `dimly lit', although much progress has been
realised. This is a subject where it is much easier to conjecture than
to prove, and we are still awash in unresolved conjectures. 

McKay also noticed in 1978 that similar coincidences hold if $\M$ and
$j(\tau)$ are replaced with the Lie group $E_8(\C)$ and
$(qj(q))^{{1\over 3}}=1+248q+\cdots$. This turns out to be much easier
to explain, and in 1980 both Kac and Lepowsky remarked that the unique
level 1 highest-weight representation of the affine algebra
$E_8^{(1)}$ has graded dimension $(qj(q))^{{1\over 3}}$.

Moonshiners have a little chip on their shoulders. Modern math,
they say, tends to be a little too infatuated with the pursuit of generalisations
 for
generalisations' sake. Surely a noble goal for math is to find interesting and
fundamentally new theorems. It can be argued that both history and common-sense
suggest that to this end it is most profitable to look simultaneously at both exceptional structures
and generic structures, to understand the special features of the former
in the context of the latter, and to be led in this way to a new generation of
exceptional and generic structures. Moonshiners would sympathise with those
biologists who study
the duck-billed platypus and lungfish rather than hide them in the closet
as monsters: BECAUSE they appear to be unique, those animals presumably have much to
teach us about our general understanding of evolution, etc.

It often seems to people that Moonshine can't be very deep: the Conway-Norton
conjectures seem to be so finite\footnote{$^{22}$}{{\smal Indeed the Moonshine
conjectures {\smit are} finite (it is enough to check the first 1200
coefficients), and a slightly weaker form was quickly
proved on a computer by Atkin, Fong and Smith [54]. However this sort
of argument adds no light to Moonshine, and tells us nothing of {\smit
V} except that it exists.}}
 and specialised. There only are 171 distinct
McKay-Thompson series $T_g$ in Monstrous Moonshine, after all. The whole point though
is to try to understand {\it why} the Monster and the Hauptmoduls are so related,
and then to try to extend and apply this understanding to other contexts.
Moonshine is still young, and our understanding remains incomplete. But
already math has benefitted: e.g.\ we now have a natural definition of $\M$
(as the automorphism group of $V^\natural$),
and Moonshine helped lead us to the rich structures of generalised Kac-Moody
algebras and vertex operator algebras.

We will see that Moonshine involves the interplay between {\it exceptional
structures} such as the number 24, the Leech lattice $\L_{24}$, the Monster
group $\M$, and the Moonshine module $V^{\natural}$, and {\it generic structures} such as
modular functions, vertex operator algebras,  generalised Kac-Moody algebras,
and conformal field theories. The following sections will introduce the
reader to many of these structures, as we use Moonshine as another
happy excuse to take a second little tour through modern mathematics.

\bigskip\bigskip\noindent{{\smcap 2.2. Ingredient \#1: Finite simple groups and
the Monster}}\bigskip

A readable introduction to the basics of finite group representation theory is
[25]. The finite simple groups are described in [33]; see also [11].
Group representations were introduced in \S 1.3.

A {\it normal} subgroup $H$ of a group is one obeying $gHg^{-1}=H$ for all
$g\in G$. These are important because the set $G/H$ of `cosets' $gH$ has a
natural group structure
precisely when $H$ is normal. Every group has two trivial normal subgroups:
itself and $\{1\}$. If these are the only normal subgroups, the group is called
{\it simple}. It is conventional to regard the trivial group $\{1\}$ as not
simple (just as 1 is conventionally regarded as not prime). An alternate
definition of a (finite) simple group $G$ is that if $\varphi:G\rightarrow H$
is any group homomorphism (i.e.\ structure-preserving map: $\varphi(gg')
=\varphi(g)\varphi(g')$), then $\varphi$ is either constant (i.e.\
$\varphi(G)=\{1\}$), or $\varphi$ is one-to-one.

The importance of simple groups is provided by the {\it Jordan-H\"older Theorem}.
By a {`composition series'} for a group $G$, we mean a nested sequence
$$G=H_0\supset H_1\supset H_2\supset\cdots\supset H_k\supset H_{k+1}=\{1\}$$
of groups such that $H_i$ is normal in $H_{i-1}$, and $H_{i-1}/H_i$ (called
a `composition factor') is
simple. Any finite group $G$ has at least one composition series. If
$H'_0\supset\cdots\supset H'_{\ell+1}=
\{1\}$ is a second composition
series for $G$, then Jordan-H\"older says that $k=\ell$ and, up to a reordering
$\pi$, the simple groups $H_{i-1}/H_i$ and $H'_{\pi j-1}/H_{\pi j}'$ are isomorphic.

For example, the cyclic group $\Z_n$ of order(=size) $n$ --- you can think of
it as the integers modulo $n$ under addition --- is simple iff $n$ is prime.
Consider the group $\Z_{12}=\langle 1\rangle$. Two composition series are
$$\eqalignno{\Z_{12}\supset\langle 2\rangle\supset\langle
4\rangle\supset\langle 0\rangle&&\cr\Z_{12}\supset\langle 
3\rangle\supset\langle 6\rangle\supset\langle 0\rangle&&\cr}$$
corresponding to composition factors $\Z_2$, $\Z_2$, $\Z_3$, and $\Z_3$, $\Z_2$,
$\Z_2$. Of course this is consistent with Jordan-H\"older. This is
reminiscent of the fact that
 $2\cdot 2\cdot 3=3\cdot 2\cdot 2$ are both prime factorisations of 12.

There is some value to regarding finite groups as a massive generalisation
of the notion of number. The number $n$ can be identified with the cyclic
group $\Z_n$. The divisor of a number corresponds to a normal subgroup,
so a prime number corresponds to a simple group. The Jordan-H\"older Theorem
generalises the uniqueness of prime factorisations. That you can build up
any number by multiplying primes, is generalised to building up a group
by semi-direct products (more generally, by group extensions): if $H$
is a normal subgroup
of $G$, then $G$ will be an extension of $H$ by the quotient group $G/H$.

Note however that $\Z_6\times\Z_2$ and ${\frak S}_3\times \Z_2$ ---
both  different
from $\Z_{12}$ --- will also have $\Z_2,\Z_2,\Z_3$ as composition factors:
unlike for numbers, `multiplication' here does not give a unique answer.
The semidirect product $\Z_2\sdprod \Z_2$ can equal either $\Z_4$
or $\Z_2\times\Z_2$, depending
on how the product is taken. More precisely, the notation $G\sdprod G'$
means a group where every element can be written uniquely as a pair
$(g,g')$, for $g\in G$ and $g'\in G'$, and where the group operation is
$(g,g')(h,h')=({\rm stuff},gh)$.

Thus simple groups have an importance for group theory approximating what
primes have for number theory. One of the greatest accomplishments of
twentieth century math is surely the classification of the finite simple
groups. (On the other hand, group extensions turn out to be technically
 quite difficult and leads one into group cohomology.) This work, completed
in the early 1980s (although gaps are continually being discovered and filled
in the
arguments), runs to approximately 15 000 journal pages, spread over 500
individual papers, and is the work of a whole generation
of group theorists. A modern revision is currently underway (see e.g.\
[34]) to simplify the proof and find and fill all gaps, but the final
proof is still expected to be around 4000 pages long. The resulting list is:
\smallskip

\item{$\bullet$} the cyclic groups $\Z_p$ ($p$ a prime);

\item{$\bullet$} the alternating groups ${\frak A}_n$ for $n\ge 5$;

\item{$\bullet$} 16 families of Lie type;

\item{$\bullet$} 26 sporadics.\smallskip

We've already met the cyclic groups. The alternating group ${\frak A}_n$ consists of the even
permutations in the symmetric group ${\frak S}_n$, and so has order(=size) ${1\over 2}\,n!$.
The groups of Lie type are
essentially Lie groups defined over finite fields\footnote{$^{23}$}{{\smal
There is a finite field with {\smit q} elements, iff {\smit q} is a
power of a prime. For each such {\smit q}, there is only 1 field of
that size. The field with prime {\smit p} elements is the integers
taken mod {\smit p}.}} ${\Bbb F}_q$ (such as
$\Z_p$), sometimes `twisted' in certain senses. The simplest
example is PSL$_n({\Bbb F}_q)$, which consists of the $n\times n$ matrices
with entries in ${\Bbb F}_q$, with determinant 1, quotiented out by the
centre of SL$_n({\Bbb F}_q)$ (namely the scalar matrices diag$(a,a,\ldots,a)$
for $a^n=1$) (except for PSL$_2(\Z_2)$ and PSL$_2(\Z_3)$, which aren't
simple).

Note that the determinant $|\rho(g)|$ for any representation $\rho$ of
any (noncyclic) simple group must be 1, otherwise we would violate the
homomorphism definition of simple group (try to see why). Also, the
centre of  any (noncyclic)
simple group must be trivial (why?).
The smallest noncyclic simple group is ${\frak A}_5$, with order 60.\footnote{$
^{24}$}{{\smal This implies, incidentally, that if {\smit G} and {\smit H}
are any two groups with the same order below 60, then they will have
the same composition factors.}} It is the same as (isomorphic to) PSL$_2(\Z_5)$ and
PSL$_2({\Bbb F}_4)$, and can also be expressed as the group of all rotations
(reflections have determinant $-1$ and so cannot belong to any simple group)
of $\R^3$ that bring a regular icosahedron back to itself.

The smallest sporadic group is the Mathieu group $M_{11}$, order 7920,
discovered in 1861\footnote{$^{25}$}{{\smal Although his arguments
apparently weren't very convincing. In fact some people, including the
Jordan of Jordan-H\"older fame, argued in later papers that the
largest of Mathieu's sporadic groups couldn't exist.}}. The largest
is the Monster $\M$, conjectured by
Fischer and Griess in 1973 and finally proved to exist by Griess in 1980.
Its order is\footnote{$^{26}$}{{\smal The inquisitive reader, hungry for
more `coincidences', may have noticed that 196883 and 21296876 --- see
(2.1.2) --- exactly divide
the order of the Monster. Indeed this will hold for any finite group:
the dimensions of the irreducible representations of a finite group will
always divide its order.}}
$$\|\M\|=2^{46}\cdot 3^{20}\cdot 5^9\cdot 7^6\cdot 11^2\cdot 13^3\cdot 17
\cdot 19\cdot 23\cdot 29\cdot 31\cdot 41\cdot 47\cdot 59\cdot 71\approx
8\times 10^{53}\ .\eqno(2.2.1)$$
20 of the 26 sporadics are involved in (i.e.\ are quotients of subgroups of)
 the Monster. Some
relations among $\M$,  the Leech lattice $\L_{24}$ and the largest
Mathieu group $M_{24}$ are given in Chapters 10 and 29 of [13].

Moonshine hints at a tantalising connection between the classification
of finite simple groups, and the classification of RCFTs discussed in 
Part 1. Speculates [23]
(page xli): ``One can certainly hope for a uniform description of the
finite simple groups as automorphism groups of certain vertex operator
algebras --- or conformal quantum field theories. If such a quantum field
theory could somehow be attached a priori to a finite simple group, the
classification of such theories, a problem of great current interest
among string theorists, might some day be part of a new approach to the
classification of the finite simple groups. On the other hand, can the
known classification of the finite simple groups help in the classification
of conformal field theories?''

\bigskip\bigskip\noindent{{\smcap 2.3. Ingredient \#2: Modular functions
and Hauptmoduls}}\bigskip

A readable introduction to some of the topics discussed in this section is
[16,42,60].

We know from complex analysis that the group SL$_2(\R)$ of $2\times 2$
matrices with real entries and determinant 1, acts on the upper-half plane
$\H=\{\tau\in\C\,|\,{\rm Im}(\tau)>0\}$ by fractional linear (or {\it M\"obius})
transformations:
$$\left(\matrix{a&b\cr c&d}\right)\cdot \tau\eqde {a\tau+b\over c\tau+d}\ .
\eqno(2.3.1)$$
For example the matrix
$S\eqde\left(\matrix{0&-1\cr 1&0}\right)$
corresponds to the function $\tau\mapsto -1/\tau$, 
 while the matrix $T\eqde\left(\matrix{1&1\cr
 0&1}\right)$ corresponds to the translation $\tau\mapsto \tau+1$. 
Since $\pm\left(\matrix{1&0\cr 0&1}\right)$
 correspond to the same M\"obius transformation, strictly
speaking our group here is PSL$_2(\R)={\rm SL}_2(\R)/\{\pm I\}$.

The only reason this action (2.3.1) of the $2\times 2$ matrices on
complex numbers
(or more precisely the Riemann sphere $\C\cup\{\infty\}$) might not
look  strange to us, is because
familiarity breeds numbness. What we really have is a natural action of
$n\times n$ matrices on $\C^n$, and this induces their action on $\C^{n-1}$
(together with a codimension-2 set of `points at infinity') by interpreting
$\C^n$ as projective coordinates for $\C^{n-1}$. Specialising to $n=2$
gives us (2.3.1). In projective geometry, `parallel lines' intersect at
$\infty$. Projective coordinates allow one to treat `finite' and `infinite'
points on an equal footing.

Consider $\G:= {\rm SL}_2(\Z)$, the subgroup of SL$_2(\R)$ consisting of
the matrices with integer entries. It can be shown that it is generated by
$S$ and $T$ (in other words, every matrix $\alpha\in\G$ can be expressed as
a monomial in $S$ and $T$). For reasons that will be clear shortly, 
consider the {extended} upper-half plane $\overline{\H}
:=\H\cup\{\i\infty\}\cup\Q$ --- the extra points $\{\i\infty\}\cup\Q$ are
called {\it cusps}. $\Gamma$ acts on $\overline{\H}$ (e.g.\ $S$ interchanges 0
and $\i\infty$). By a {\it modular function} for $\G$,
we mean a meromorphic function $f:\overline{\H}\rightarrow\C$, symmetric with respect
to $\G$: i.e.\ $f(\alpha(\tau))=f(\tau)$ for all $\alpha\in\G$. Note that we
require $f$ to be meromorphic at the cusps (e.g.\ polynomials are meromorphic
at $\i \infty$, but $e^z$ is not).

It is not obvious why modular functions should be interesting, but in fact
they are one of the most fundamental notions in modern number theory
(see the last paragraph of \S1.6). For example,
consider the question of writing numbers as sums of squares. We can write
$5=1^2+(-2)^2=(-1)^2+1^2+0^2+1^2+(-1)^2$, to give a couple of trivial examples.
Let $N_n(k)$ be the number of ways we can write the integer $n$ as a sum
of $k$ squares, counting order and signs. For example $N_5(1)=0$ (since
5 is not a perfect square), $N_5(2)=8$ (since $5=(\pm 1)^2+(\pm 2)^2=(\pm 2)^2
+(\pm 1)^2$), $N_5(3)=24$, etc. Their
 generating functions are:\footnote{$^{27}$}{{\smal A fundamental principle
in math is: whenever you have a subscript with an infinite range, make a power
series (called a {\smit generating} {\smit function}) out of it.}}
$$\sum_{n=0}^\infty N_n(k)\, q^n=(\theta_3(q))^k\ ,$$
where
$$\theta_3(q)=1+2q+2q^4+\cdots=\sum_{n\in\Z}q^{n^2}$$
is called a {\it theta function}. It turns out that $\theta_3$ transforms nicely
with respect to $\G$, once we make the change-of-variables $q=\exp[\pi\i \tau]$.
This takes work to show. For example, $\theta_3$ is clearly invariant
under the action of $\left(\matrix{1&2\cr 0&1}\right)$, and a little work
(from e.g.\ Poisson summation)
shows that $\left(\matrix{0&-1\cr 1&0}\right)$ takes $\theta_3(\tau)$ to
$\sqrt{\tau\over \i}\,\theta_3(\tau)$. $\theta_3$ is not
precisely a modular function (it  is  a {`modular form of weight ${1\over 2}$'}
for $\Gamma_0(4)$), but this simple example illustrates the point that $\G$
 (and related groups) appear throughout number theory. More on this shortly.

That important change-of-variables $q=\exp[\pi\i\tau]$ was introduced by
Jacobi early last century, in his analysis of `elliptic integrals'.
The theory is beautiful and poorly remembered today, which is very
disappointing considering how much of modern math was touched by
 it. I strongly recommend the book [10], written over a century ago;
the style and motivation of math in our century is different from that in
Jacobi's, and we've lost a little in motivation what we've gained in
power. I'll briefly sketch Jacobi's theory.

Just as we could develop a theory of `circular functions' (i.e.\  sine etc.)
 starting from the integral $s(a)=\int_0^a{dx\over \sqrt{
1-x^2}}$, so can we develop a theory of `elliptic functions' starting
from the `elliptic integral' $F(k,a)=\int_0^a{dx\over
\sqrt{(1-x^2)(1-k^2x^2)}}$. Inverting $s(a)$ gives a function both more 
useful and with nicer properties than $s(a)$: we call it
$\sin(u)$. Similarly, for any  $k$ the elliptic function
${\rm sn}(k,u)$ is defined by $u=F(k,{\rm sn}(k,u))$. Just as we can define a numerical
constant $\pi$ by $\sin({1\over 2}\pi)=1$ (i.e.\ ${1\over 2}\pi=\int_0^1
{dx\over \sqrt{1-x^2}}$), we get a function $K(k)=\int_0^1{dx\over \sqrt{
(1-x^2)(1-k^2x^2)}}$. Just as $\sin(u)$ has period $4({1\over 2}\pi)$,
so has ${\rm sn}$ $u$-period $4K(k)$. ${\rm sn}$ also turns out
to have $u$-period $4\i\,K(k')$ where $k'=\sqrt{1-k^2}$ --- today we
take this as the
starting point and define an elliptic function to be
doubly periodic (see [42] or Cohen in [60]).

The theta functions aren't elliptic
functions, but they are closely related, and e.g.\ ${\rm sn}$ can be written
as a quotient of them. In Jacobi's language, we have
$$\theta_3({\i\,K(k')\over K(k)})=\sqrt{{2K(k)\over\pi}}\ .$$
The `modular transformation' $\tau\mapsto {-1\over\tau}$
corresponds to interchanging the `modulus' $k$ with the `complementary
modulus' $k'$, and thus is completely natural in Jacobi's theory. The
important formula $\theta_3({-1\over \tau})=\sqrt{{\tau\over \i}}\,
\theta_3(\tau)$ is trivial here.

A certain interpretation of modular functions also indicates their
usefulness,
and played an important role in Part 1. A {\it torus} is something
that looks like the surface of a bagel, at least as far as its topology is
concerned. For example, the Cartesian product $S^1\times S^1$ of circles
is a torus (think of one circle being the contact-circle of the bagel with
the table on which it rests, then from each point on that horizontal circle
imagine placing a vertical circle perpendicular to it,
like a rib; together all these ribs fill out the bagel's surface).
A more sophisticated example of a torus is an elliptic curve
 (a complex curve of the form $y^2=ax^3+bx^2+cx+d$ and a special point
on it playing the role of 0). A final
 example is the quotient $\C/\L$ of the complex plane $\C$ with a 2-dimensional
 { lattice} $\L$ (we saw lattices in \S 1.6; $\L$ here will be
a discrete doubly-periodic set of points in $\C$, containing 0). It turns out
 that certain equivalence classes of tori (e.g.\ with respect to conformal or
 complex-analytic equivalence) always contain a representative torus of the
 form $\C/\L$, where $\L$ consists of all points $\Z+\Z\, \tau$,
for some $\tau\in\H$. (Incidentally, the cusps correspond to degenerate tori.)
In other words, these equivalence
 classes are parametrized by complex numbers $\tau$ in $\H$. So if we have a
complex-valued  function $F$ on the set of all tori, which is e.g.\
conformally invariant (an example is the genus-one partition function ${\cal Z}$
in conformal field theories --- see \S1.1), then we can consider $F$ as a well-defined
function $F:\H\rightarrow\C$. However, it turns out that different points
$\tau$ in $\H$ correspond to the same equivalence class of tori: e.g.\ the lattice
for $\tau$ is the same as that for $\tau+1$, and these are a rescaling
of that for $-1/\tau$. Thus $F(\tau)=F(\tau+1)=F(-1/\tau)$, because
$\tau,\tau+1,-1/\tau$ all represent equivalent tori. Since $\tau\mapsto\tau+1$
and $\tau\mapsto -1/\tau$ generate PSL$_2(\Z)$, what in fact we
find is that $F$ has $\G$ as its group of symmetries. One often says that $\G$
is the `modular group of the torus', and that the orbit space
$\G\backslash \H$ is the
`moduli space' of (conformal equivalence classes of) tori. 
 $\H$ is called its `Teichm\"uller space' or `universal cover'. This is
 exactly analogous to $S^1=\R/\Z$: $\R$ is its universal cover and $\Z$
 is its `modular group' (or `mapping class group'). Another example:
 the Teichm\"uller space for (conformal equivalence classes of) `pair-of-pants', or equivalently a disc minus two
 open interior disks, is $\R^3$ (an ordered triple), while its modular group is
 the symmetric group ${\frak S}_3$ and its moduli space consists of unordered triples.
 Incidentally, we write $\G\backslash{\H}$
 instead of ${\H}/\G$ because the group $\G$ acts on $\H$ `on the left'.
A good introduction to the geometry here is [55].

In any case, a surprising number of innocent-looking questions in number
theory can be dragged (usually with effort) into the richly developed realm
of elliptic curves and modular functions, and it is there they are often solved.
For instance, we all know the ancient Greeks were interested in Pythagorean
triples: find all integer solutions $a,b,c$ to $a^2+b^2=c^2$, i.e.\ find all
integer (or if you prefer, rational) right-angle triangles. They solved this
by elementary means: choose any integers (or rationals) $x,y$ and put
$u={x^2-y^2\over x^2+y^2}$, $v={2xy\over x^2+y^2}$; then $u^2+v^2=1$ and
(multiplying by the denominator) this gives all Pythagorean triples.

There are two ways of extending this problem. One is to ask which $n\in\Z$
can arise as areas of these rational right-angle triangles. It turns out
$n=5$ is the smallest one: $a={3\over 2}$, $b={20\over 3}$, $c={41\over 6}$
works ($5={1\over 2}({3\over 2})({20\over 3})$ and $({3\over
2})^2+({20\over 3})^2=({41\over 6})^2$). This is a hard problem ---
just try to show $n=1$ cannot work. $n=157$ turns
out to work: the simplest triangle has $a$ and $b$ as quotients of integers
of size around $10^{25}$, and $c$ as the quotient of integers around $10^{47}$.
Although this problem was studied by the ancient Greeks and also by the
Arabs in the 10th century, it was finally cracked in the 1980s. It was solved
by first translating it into the question of whether the elliptic curve
$y^2=x^3-n^2x$ has infinitely-many rational points, and then applying all
the rich 20th century machinery to answering that question.

The other continuation of the Pythagorean triples question is more famous:
find all integer solutions to $a^n+b^n=c^n$ (or equivalently all rational
solutions to $a^n+b^n=1$). 350 years ago Fermat wrote in the margin of the
book he was reading (the book was describing the Greek solution to
Pythagorean triples) that he had found a ``truly marvelous'' proof that for
$n>2$ there are no nontrivial solutions, but that the margin was too narrow to
contain it. This result came to be known as `Fermat's
Last Theorem'\footnote{$^{28}$}{{\smal It was called his `Last Theorem' because
it was the last of his 48 margin notes to be proved by other mathematicians
--- another one is discussed in Section 1.8. The story of Fermat's Last
Theorem is a fascinating one, but alas this footnote is too small to do
it credit. See for instance the excellent book [53].}}
and despite considerable effort no one has succeeded in
rediscovering his proof. Most people today believe that Fermat soon
realised his `proof' wasn't valid, otherwise he would have alluded to it
in later letters. In any case, a very long and complicated proof was finally
achieved in the 1990s: the `Taniyama conjecture' says that a certain function
associated to any elliptic curve over $\Q$ will be modular; if $a^n+b^n=c^n$
for some $n>2$, then the elliptic curve $y^2=x^3+(a^n-b^n)x^2-a^nb^n$ will violate the Taniyama
conjecture; finally Wiles proved the Taniyama conjecture is true.

To most mathematicians, the `area-$n$ problem' and `Fermat's Last Theorem'
are interesting only because they can be related to elliptic curves and
modular forms  --- it's easy to ask hard questions in math, but most
questions tend to be stale. Number theory is infatuated with modular
stuff because (in increasing order of significance) (a) it's exceedingly rich, with lots of connections to other
areas of math and math phys; (b) it's a battleground on which many
innocent-looking but hard-to-crack problems can be slain; and (c) 
last generation's number theorists also worked on modular stuff.

In any case, modular functions turn out to be important for math (and mathematical
physics) even though they may at first glance look artificial. Poincar\'e
explained how to study them. He said to look at the orbits of $\overline{\H}$ with
respect to $\G$. For example, one orbit, hence
one point in $\G\backslash\overline{\H}$, contains all cusps. We write this
as $\G\backslash \overline{\H}$,
and give it the natural topological structure (i.e.\
2 points $[\tau],[\tau']\in\G\backslash \overline{\H}$ are considered `close'
if the 2 sets $\G\tau,\G\tau'$ nearly overlap). Note first that by applying
$T$ repeatedly,
every point in $\H$ corresponds to a point in the vertical strip $-{1\over 2}
\le {\rm Re}(\tau)\le {1\over 2}$ --- in fact to a unique point in that strip, if we
 avoid the two edges. $S$ is an inversion through the unit circle, so
 it permits us to restrict to those points
 in the vertical strip which are distance at least 1 from the origin. The
 resulting region $R$ is called a fundamental region for $\G$. Apart from the
 boundary of $R$, every $\G$-orbit will intersect $R$ in one and only one point.

What should we do about the boundary? Well, the edge Re$(\tau)=-{1\over 2}$
gets mapped by $T$ to the edge Re$(\tau)={1\over 2}$, so we should identify
(=glue together) these. The result is a cylinder running off to infinity,
with a strange lip at the bottom. $S$ tells us how we should close that lip:
identify $\i e^{\i \theta}$ and $\i e^{-\i\theta}$. This seals the bottom
of the cylinder, so we get an infinitely tall cup with a strangely
puckered base. In fact the top of this
cup is also capped off, by the cusp $\i\infty$. So what we have (topologically
speaking) is a {\it sphere}. It does not look like a smooth sphere, but
in fact it inherits the smoothness of $\H$.

Incidentally, topological manifolds of
dimension $\le 3$ always have a unique compatible smooth structure.
 `Topological structure' means you
can speak of continuity or closeness, `smooth structure' means you can
also do calculus. On the other hand, $\R^4$ has infinitely many smooth structures
compatible with its topological structure; mysteriously, all other
Euclidean spaces $\R^n$ have a unique smooth structure!
 Thus both mathematics and physics single out 4-space. Coincidence???

So anyways, what this construction of $\G\backslash\overline{\H}$
means is that a modular function can be reinterpreted
as a meromorphic complex-valued function on this sphere. This is very useful,
because our undergraduate complex variables class taught us all about meromorphic
complex-valued functions $f$ on the Riemann sphere $\C\cup\infty$. There
are many meromorphic functions on $\C$, but to also be meromorphic at $\infty$
forces $f$ to be {\it rational}, i.e.\ $f(w)={{\rm some\ polynomial}\ P(w)
\over {\rm some\ polynomial}\ Q(w)}$, where $w$ is the complex parameter
on the Riemann sphere. So our modular function $f(\tau)$ will simply be some
rational function $P/Q$ evaluated at the change-of-variables function
$w=c(\tau)$ which maps us from our sphere $\G\backslash \overline{\H}$ to the Riemann sphere.
There are many different choices for this function $c(\tau)$, but the standard
one is $c(\tau)=j(\tau)$, the $j$-function of (2.1.1)\footnote{$^{29}$}{{\smal
Historically, {\smit j} was the standard choice, but in Moonshine the preferred
choice would be the function {\smit J} = {\smit j} -- 744
with zero constant term.}}. Thus, any modular function
can be written as a rational function $f(\tau)=P(j(\tau))/Q(j(\tau))$ in the $j$-function.
Conversely, any such function will be modular.

This is  analogous
to saying that any function $g(x)$ periodic under $x\mapsto x+1$
can be thought of as a function on the unit circle $S^1\subset \C$ evaluated at
the change-of-variables function $x\mapsto e^{2\pi\i x}$, and hence has a
Fourier expansion $\sum_ng_n\,\exp[2\pi\i nx]$.

We can generalise this argument. Consider a subgroup $G$ of SL$_2(\R)$ which
is both not too big, and not too small. `Not too big' means it should be
{\it discrete}, i.e.\ the matrices in $G$ can only get so close to the
identity matrix $\left(\matrix{1&0\cr 0&1}\right)$. To make sure $G$ is
`not too small', it is enough to require that $G$ contains some subgroup
 of the form
$$\Gamma_0(N):=\{\left(\matrix{a&b\cr c&d}\right)\in {\rm SL}_2(\Z)\,|\,
c\equiv 0\ ({\rm mod}\ N)\,\}\ ,\eqno(2.3.2a)$$
i.e.\ $G$ must contain all matrices in $\G$ whose bottom-left entry is a
multiple of $N$. So $G$ must contain $T$, for example. We will also
be interested only in those $G$ which obey
$$\left(\matrix{1&t\cr 0&1}\right)\in G\ \Rightarrow\ t\in\Z\eqno(2.3.2b)$$
i.e.\ the only translations in $G$ are by integers. We will call a function
$f:\overline{\H}\rightarrow\C$ a {\it modular function for} $G$ if it is meromorphic
(including at the { cusps} $\Q\cup\{\i\infty\}$), and if
also $f$ is symmetric with respect to $G$: $f\circ\alpha=f$ for all $\alpha
\in G$. This implies we will be able to expand $f$ as a Laurent series in $q$.
We analyse this as before: look at the orbit space
$\Sigma=G\backslash \overline{\H}$;
because $G$ is not too big, $\Sigma$ will be a (Riemann) surface; because $G$ is
not too small, $\Sigma$ will be compact.

The compact Riemann surfaces have been classified (up to homeomorphism ---
i.e.\ considering only topology as relevant), and are characterised by
a number called the {\it genus}. Genus 0 is a sphere, genus 1 is a torus,
genus 2 is like two tori resting side-by-side, etc. For example, the surface
of a wine glass, or a fork, is topologically a sphere, while a coffee cup
and a key will (usually) be tori. Eye glasses with the lenses popped out
is a 2-torus, while a ladder with $n$ rungs on it has
genus $n-1$.

We will call $G$ `genus $g$' if its surface $\Sigma$ has genus $g$. For example,
 $G=\G_0(2)$ and $G=\G_0(25)$ are both genus 0, while $\G_0(50)$ is genus 2
and $\G_0(24)$ is genus 3. Once again, we are interested here in
 the genus 0 case. As before, this means that there is a change-of-variables
 function we'll denote $J_G$ which has the property that it's a modular
 function for $G$, and all other modular functions for $G$ can be written
 as a rational function in it. Because of (2.3.2), we can choose $J_G$
 to look like
$$J_G(\tau)=q^{-1}+a_1(G)\,q+a_2(G)\,q^2+\cdots$$
So $J_G$ plays exactly the same role for $G$ that $J:=j-744$  plays for
$\G$. $J_G$ is called the Hauptmodul for $G$. (Incidentally for genus $>0$,
two generators, not one, are needed.)

For example, $\G_0(2)$, $\G_0(13)$ and $\G_0(25)$ are all genus 0, with
Hauptmoduls 
$$\eqalignno{J_{2}(\tau)=&\,q^{-1}+276\,q-2048\,q^2+11202\,q^3-49152\,q^4+184024\,
q^5+\cdots&(2.3.3)\cr
J_{13}(\tau)=&\,q^{-1}-q+2\,q^2+q^3+2\,q^4-2\,q^5-2\,q^7-2\,q^8+q^9+\cdots&
(2.3.4)\cr
J_{25}(\tau)=&\,q^{-1}-q+q^4+q^6-q^{11}-q^{14}+q^{21}+\cdots&(2.3.5)\cr}$$
The smaller the modular group, the smaller the coefficients of
the Hauptmodul. In this sense, the $j$-function is optimally bad among the
Hauptmoduls: e.g.\ for it $a_{23}\approx 10^{25}$.

An obvious question is, how many genus 0 groups (equivalently, how many
Hauptmoduls) are there? It turns out that
$\G_0(p)$ is genus 0, for a prime $p$, iff $p-1$ divides 24. Thompson in
1980 proved that for any $g$, there are only finitely many genus $g$ groups obeying our two
conditions (2.3.2). In particular this means there are only finitely many Hauptmoduls.
Over 600 Hauptmoduls with integer coefficients $a_i(G)$ are presently known.

\bigskip\bigskip\noindent{\smcap 2.4. The Monstrous Moonshine Conjectures}
\bigskip

We are now ready to make precise the main conjecture of Conway and Norton [12].
(We should emphasise though that there have been several other conjectures,
some of which turned out to be partially wrong.)

They conjectured that for each element $g$ of the Monster $\M$, there is
a Hauptmodul
$$J_g(\tau)=q^{-1}+\sum_{n=1}^\infty a_n(g)\,q^n\eqno(2.4.1)$$
for a genus 0 group $G_g$ such that each coefficient $a_n(g)$ is an integer, and
for each $n$ the map $g\mapsto a_n(g)$ is a character of $\M$. They also
conjectured that $G_g$ contains $\G_0(N)$ as a normal subgroup, for
some  $N$ depending on the order of $g$.

Another way of saying this is that there exists an infinite-dimensional
graded representation       $V=V_{-1}\oplus\bigoplus_{n=1}^\infty V_n$
of $\M$ such that the McKay-Thompson series $T_g(\tau)$ in (2.1.4) 
is a Hauptmodul.

There are around $8\times 10^{53}$ elements to the Monster, so naively we
may expect around  $8\times 10^{53}$ different Hauptmoduls $J_g=T_g$. However
the character of  a representation evaluated at $g$ and at $hgh^{-1}$ will
always be the same, so $J_g=J_{hgh^{-1}}$. Hence the relevant quantity is the number of conjugacy
classes, which for $\M$ is only 194. Moreover, a character evaluated at $g^{-1}$
will always be the complex conjugate of its value at $g$, but here all
character values $\chi_{V_n}(g)$ are integers (according to the conjecture).
Thus  $J_{g}=J_{g^{-1}}$. 
The total number of distinct Hauptmoduls $J_g$ arising in Monstrous Moonshine
turns out to be only 171.

For example, if we choose $g$ to be the identity, we recover $T_{id.}=J$.
It turns out that there are precisely 2 different conjugacy classes of
order 2 elements,  one of them giving the Hauptmodul $J_2$ in (2.3.3).
Similarly for 13, but $J_{25}$ doesn't correspond to any conjugacy class
of $\M$.

Moonshine provides an explanation for a forgotten mystery of classical
mathematics: why are the
coefficients of the $j$-function {\it positive} integers? On the other hand,
that they are {\it integers} has long been important to number theory
(complex multiplication, class field theory --- see e.g.\ [16]).

There are lots of other less important conjectures. One which played  a
role in ultimately proving the main conjecture involves the {\it replication
formulae}. Conway-Norton want to think of the Hauptmoduls $J_g$ as
being intimately connected with $\M$; if so, then the group structure
of $\M$
should somehow directly relate different $J_g$. In particular, consider the
power map $g\mapsto g^p$. Now, it was well-known that $j(\tau)$ has the property that
$j(p\tau)+j({\tau\over p})+j({\tau+1\over p})+\cdots+j({\tau+p-1\over p})$
equals a polynomial in $j$, for any prime $p$ ({\it sketch of proof:}
 it's a modular
function for $\G$, and hence equals a rational function of $j$; since its
only poles will be at the cusps, the denominator polynomial must be trivial).
Hence the same will hold for $J$. Explicitly we get
$$\eqalignno{J(2\tau)+J({\tau\over 2})+J({\tau+1\over 2})=&\,J^2(\tau)-2a_1
&(2.4.2a)\cr
J(3\tau)+J({\tau\over 3})+J({\tau+1\over 3})+J({\tau+2\over 3})=&\,J^3(\tau)-
3a_1J(\tau)-3a_2&(2.4.2b)\cr}$$
where $J(\tau)=\sum_ka_kq^k$. Slightly more complicated formulas hold in fact
for any composite $n$. Conway and Norton conjectured that these formulas
have an analogue for the Moonshine functions $J_g$ in (2.4.1). In particular,
(2.4.2) become for any $g\in\M$
$$\eqalignno{J_{g^2}(2\tau)+J_g({\tau\over 2})+J_g({\tau+1\over 2})=&\,J^2_g
(\tau)-2a_1(g)&(2.4.3a)\cr
J_{g^3}(3\tau)+J_g({\tau\over 3})+J_g({\tau+1\over 3})+J_g({\tau+2\over 3})=
&\,J_g^3(\tau)-3a_1(g)\,J_g(\tau)-3a_2(g).&(2.4.3b)\cr}$$
These are examples of  the replication formulae.

`Replication' concerns the power map $g\mapsto g^n$ in $\M$. Can
Moonshine see more of the group structure of $\M$? A step in this
direction was made by Norton [47], who associated a Hauptmodul to
commuting elements $g,h$ in $\M$. Physically [19], this corresponds to
orbifold traces, i.e.\ the $V^\natural$ RCFT with boundary conditions
twisted by $g$ and $h$ in the `time' and `space' directions. Still, we
would like to see more of $\M$ in Moonshine. 

An important part of the Monstrous Moonshine conjectures came a few years
after [12]. Frenkel-Lepowsky-Meurman [23] constructed a graded infinite-dimensional
representation $V^{\natural}$ of $\M$ and conjectured (correctly) that it
is the representation in (2.1.4). $V^\natural$ has a very rich algebraic
structure, which will be discussed in \S 2.6.

A major claim of [23]
was that $V^\natural$ is a `natural' structure (hence their notation).
To see what they mean by that, it's best to view another simpler
 example of a natural construction: that of the {\it Leech lattice} $\L_{24}$.
Recall the discussion of (root) lattices in \S1.6.

$\L_{24}$ is one  of the most interesting lattices, and
is related to Moonshine.
It can be defined using `laminated lattices'. Start with the
0-dimensional lattice $\L_0=\{0\}$, which consists of just a single point.
Use it to construct a 1-dimensional lattice, with minimal (nonzero) norm
4, built out of infinitely many copies of $\L_0$ laid side by side.
The result of course is simply the even integers $2\Z$, which we will call
here $\L_1$. Now construct a 2-dimensional lattice,
of minimum norm 4, built out of infinitely many copies of $\L_1$ laid next
to each other. There are lots of ways to do this, but choose the densest
lattice possible. The result is unique: it is the hexagonal lattice $A_2$
scaled
by a factor of $\sqrt{2}$: call it $\L_2$. Continue in this way: $\L_3$, $\L_4$,
$\L_5$, $\L_6$, $\L_7$, and $\L_8$ will be the root lattices $A_3$, $D_4$,
$D_5$, $E_6$, $E_7$ and $E_8$, respectively, all scaled by $\sqrt{2}$.
See [13] chapter 6 for a more complete treatment of laminated lattices.

The 24th repetition of this construction yields the 
Leech lattice. It is the unique 24-dimensional self-dual lattice with no
norm-2 vectors, and provides among other things the densest known packing of
23-dimensional spheres in $\R^{24}$. Many of its
properties are discussed throughout [13].
So lamination provides us with a sort of
no-input construction of the Leech lattice, and a good example of the
mathematical meaning of `natural'.
After dimension 24, it seems chaos results from the lamination procedure (there
are 23 different 25-dimensional lattices that have an equal
right to be called $\L_{25}$, and over $75\,000$ are expected for $\L_{26}$).

It is natural to ask about Moonshine for other
groups. There is a partial Moonshine for the Mathieu groups 
$M_{24}$ and $M_{12}$ (which have about $2\times 10^8$ and $10^5$
elements resp.), the automorphism group .0 of 
$\L_{24}$ (which has about $8\times 10^{18}$ elements), and a few others
--- see e.g.\ [49]. These groups are either simple or almost simple
(e.g.\ .0 is the direct product of $\Z_2$ with the simple group .1). More 
generally, 
there will be some sort of Moonshine for any group which is the
automorphism group of a vertex operator algebra; the finite simple groups of
Lie type should be automorphism groups of VOAs closely related to the
affine algebras except defined over fields like $\Z_p$.

There is a geometric side to Moonshine, associated to names like
Lian-Yau and Hirzebruch. In particular, Hirzebruch's `prize question'
asks for the construction of a 24-dimensional manifold on which $\M$
acts, whose twisted elliptic genus are the McKay-Thompson series. This
is still open.

It should be emphasised that Monstrous Moonshine is a completely unexpected
connection between finite groups and modular functions. Although there
has been enormous progress in our understanding of this connection (so much
so, that Richard Borcherds won the 1998 Fields medal for his work on this),
there still is mystery at its heart. In particular, 
that $\M$ is associated with {\it modular functions} can be explained
mathematically by it being the automorphism group of the Moonshine VOA $V^{\natural}$,
and physically by the associated RCFT, but what is so special about $\M$ that
these modular functions should be genus 0? We will come back to this
in \S2.9.

\bigskip\bigskip\noindent{{\smcap 2.5. Formal Power Series}}
\bigskip

{\it Vertex (operator) algebras} (VOAs) are a mathematically precise formulation of
the notion of {\it W-algebra} or {\it chiral algebra}\footnote{$^{30}$}{{\smal An alternate
(and much more complicated) mathematical formulation of chiral algebra
is due to Beilinson and Drinfeld, and belongs to algebraic geometry.
See [28] for a good --- but still difficult --- review.}} which is so
central  to conformal field theory (see \S1.1). VOAs were first
defined by Borcherds, and their theory has since been developed by a
number of people (Frenkel, Lepowsky, Meurman, Zhu, Dong, Li, Mason,
Huang, ...). Because our primary
motivation here is Moonshine, I will only focus on one aspect of their
theory (the connection with Lie algebras). Useful
to consult while reading this review are the notes [27]
 --- they take a more analytic approach to
many of the things we discuss, and their approach (namely that of
CFT) motivates beautifully much of VOA theory.

In quantum field theory the basic object is the quantum field, which roughly
speaking is a choice of operator $\hat{A}(x)$ at each space-time point $x$.
`Operator' means something that `operates on' functions or vectors. E.g.\ an
indefinite integral is an operator, as is a derivative. 
The operators in the QFT act on the space spanned
by the states $|\star\rangle$, and together form an infinite-dimensional vector
space (e.g.\ a  $C^*$ algebra) --- this infinite-dimensionality
of QFT is a major source of its mathematical difficulties, and QFT still
has not been put on completely satisfactory mathematical grounds.

But another difficulty is that the quantum field $\hat{A}$ really isn't
an operator-valued {\it function} of space-time. `Function' is too narrow a
concept. For example,
one of the most familiar `functions' in quantum mechanics is the Dirac
delta $\delta(x)$. You see it for example in the canonical commutation
relations: e.g.\ for a scalar field $\hat{\varphi}$, we have $[\hat{\varphi}
(\vec{x},t),{\partial\over\partial t}\hat{\varphi}(\vec{y},t)]=
\i\hbar\delta^3 (\vec{x} -\vec{y})$. $\delta(x)$ has the property
that for any other smooth function $f$,
$$\int_{-1}^1 f(y)\,\delta(y)\, {\rm d}y=f(0)\ , \qquad \int_{-1}^1 f(y)\,\delta'(y)\,{\rm
d}y=f'(0)\ ,$$
etc. The problem is that $\delta(x)$ isn't a function --- no function could
possibly have those properties.

One way to make sense of `functions' like the Dirac delta and its derivatives
is distribution theory. Although it was first informally used in
physics, it was rigourously developed around 1950
by Laurent Schwartz, and uses the idea of test functions. See e.g.\ [15].

What I will describe now is an alternate approach, {\it algebraic} as
opposed to {\it analytic}. These two approaches are not equivalent: you can do
some things in one approach which you can't do in  the other. But the
algebraic approach is considerably simpler technically
--- no calculus or convergence to worry about ---
and it is remarkable how much can still be captured. This approach is the
starting point for the VOA story described next section,
and was first created around 1980 by Garland and Date-Kashiwara-Miwa.
Keep in mind that what we are trying to capture is an operator-valued
`function' on space-time. Space-time in CFT is 2-dimensional, and so
we can think of it (at least locally) as being on the complex plane $\C$
(more precisely, 
we will usually associate the space-time point $(x,t)$ with the complex number
$z=e^{t+\i x}$). Good introductions to the material in this section are
[23,39,31].

Let $W$ be any vector space. We are most interested in it being an
infinite-dimensional space of matrices (i.e.\ operators on an infinite-dimensional
 space), but forget that for now.
Define $W[[z,z^{-1}]]$ to be the set of all formal series $\sum_{n=-\infty}^\infty
w_n z^n$, where the coefficients $w_n$ lie in our space $W$.
We don't ask here whether a given series converges or diverges --- 
$z$ is merely a formal variable.
We will also be interested in  $W[z,z^{-1}]$ (Laurent polynomials).
We can add these formal series in the usual way, and multiply them by
numbers (scalars) in the usual way. 

Remember our ultimate aim here: we want to capture quantum fields. So we
want our formal series to be operator-valued. The way to accomplish this
is to choose $W$ to be a vector space of operators, or matrices if you
prefer. A fancy way to say this is `$W={\rm End}(V)$', which means the
things in $W$ operate on vectors in $V$. If we take $V=
\C^m$, then we can think of $W$ as being the space of all $m\times m$ complex
matrices.
We are ultimately interested in the case $m=\infty$, but we won't lose much
now by taking $m=1$, which would mean formal power series
with numerical coefficients.

Because our coefficients $w_n$ are operators,
we can multiply our formal series. We define multiplication in the
usual way. For example, consider
$W=V=\C$, and take $c(z)=z^{21}-5z^{100}$ and $d(z)=\sum_{n=-\infty}^\infty z^n$.
Then
$$c(z)\,d(z)= \sum_{n=-\infty}^\infty z^{n+21}-5\sum_{n=-\infty}^\infty
z^{n+100}=\sum_{n=-\infty}^\infty z^n-5\sum_{n=-\infty}^\infty z^n
=-4d(z)\ .$$
So far so good. Now try to compute the square $d(z)^2$. You get infinity.
So the lesson is: you can't always
multiply in $W[[z,z^{-1}]]$. We'll come back to this later.

But first, look again at that first product: $c(z)\,d(z)=-4d(z)$. One thing
it tells us is that {\it we can't always divide} (certainly $c(z)$ and $-4$ are
two very different power series!). But there's another lesson here: if you work
out a few more multiplications of this kind, what you'll find is that
$f(z)\,d(z)=f(1)\,d(z)$ for any $f$, at least for those $f$ for which
$f(1)$ exists (e.g.\ any $f\in W[z,z^{-1}]$). Thus $d(z)$ is what we would
call the Dirac delta $\delta(z-1)$! (You can think of it as the Fourier
expansion of the Dirac delta, followed by a change of variables).
Unfortunately, the standard notation here is to write it without the `$-1$':
$$\delta(z):=\sum_{n=-\infty}^\infty z^n$$
and that is the notation we will also adopt. Similarly, $\delta(az)$ and
$\delta'(z)$ etc (which are the formal series defined in the obvious way) 
act on $W[z,z^{-1}]$ in the way one would expect: $f(z)\,\delta(az)= f({1\over a})
\,\delta(az)$
and $f(z)\,\delta'(z)=f'(1)\,\delta'(z)$. So of course it makes perfect sense that
we couldn't work out $d(z)^2$: we were trying to square the Dirac delta,
which we know is impossible!

A similar theory can be developed for several variables $z_i$,
with identities such as
$f(z_1,z_2)\,\delta(z_1/z_2)=f(z_2,z_2)\,\delta(z_1/z_2)=f(z_1,z_1)\,
\delta(z_1/z_2)$.

But we must not get too overconfident:\smallskip

\noindent{{\bf Paradox 1.}} Consider the following product:
$$\delta(z)=\bigl[(\sum_{n=0}^\infty z^n)\,(1-z)\bigr]\,\delta(z)=
(\sum_{n=0}^\infty z^n)\,\bigl[(1-z)\,\delta(z)\bigr]=(\sum_{n=0}^\infty 
z^n)\,\bigl[0\,\delta(z)\bigr]=0\ .$$

When physicists are confronted with `paradoxes' such as this, they tend
to respond by keeping them in the back of their mind, by treading with care
when they are involved in a calculation which reminds them of one of the
paradoxes, and otherwise trusting their instincts. Mathematicians typically
over-react: they kick themselves for getting overconfident and walking
head-first into a `paradox', and then they devise some rule which will
absolutely guarantee that that paradox will always be safely avoided in the
future. We will follow the
mathematicians' approach, and in the next few paragraphs will describe
their rule for avoiding Paradox 1: to forbid certain innocent-looking
products. 

Remember that we are actually interested in the vector space $W={\rm End}(V)$.
Suppose we have infinitely many matrices $w_i\in {\rm End}(V)$. We will 
call them {\it summable}
if for every column vector $v\in V$, only finitely many products $w_i(v)\in V$
are different from 0. In other words, only finitely many of the matrices
$w_i$ have a nonzero first column, only finitely many have a nonzero
second column, $\ldots$.

We will certainly have a well-defined sum $\sum_i w_i(z)$ if for each fixed
$n$, the set $\{w_i(n)\}$ (as $i$ varies) of matrices is summable. All
other sums are forbidden. We will certainly have a
well-defined product\footnote{$^{31}$}{{\smal {\smit m} here will be finite:
we permit infinite sums but only finite products.}}
$\prod_{i=1}^m w_i(z)$ if for each $n$, the set $\{w_1(n_1)\,
w_2(n_2)\,\cdots\,w_m(n_m)\}_{\sum n_i=n}$ (vary the $n_i$ 
subject to $\sum_i n_i=n$) is
summable. All other products are forbidden. This is
reasonable because the sum of those matrix products $w_1(n_1)\,
\cdots\,w_m(n_m)$ will precisely equal the $n$th coefficient of the
product $\prod_{i=1}^mw_i(z)$.

Note that there are certainly more general ways to have a well-defined
product (or sum). For example, according to our rule,
we cannot even add $\sum_n 2^{-n}$!
This way has the advantage of not touching the more
complicated realm of convergence issues. We are doing {\it algebra}
here, not {\it analysis}. The way out of Paradox 1 is that
$(\sum z^n)(1-z)$ doesn't equal 1 --- rather, it's a forbidden product.

An interesting  consequence of the fact that we are doing
algebra instead of analysis is that the product $z^{{1\over 2}}\delta(z)$
here does not and cannot equal $1^{{1\over 2}}\delta(z)=\delta(z)$ --- their
formal power series are very different. In hindsight this `failing' is
understandable: algebraically, it seems artificial to prefer the positive root
 of 1 over the negative root.

\smallskip\noindent{{\bf Paradox 2.}} Expand ${1\over 1-z}$ in a formal power
series in $z$ to get $\sum_{n\ge 0}z^n$. Next, expand ${1\over 1-z}=
{-z^{-1}\over 1-z^{-1}}$ in a formal power series in $z^{-1}$ to get
$-\sum_{n<0}z^n$. Subtract these; we presumably should get 0, but we actually
get $\delta(z)\,$!\smallskip

The analytic explanation is that the first expression converges only for $|z|<1$,
and the second for $|z|>1$, so it would be naive to expect their difference
to be 0. We see from this `paradox' that {\it it really
matters in which variable we expand rational functions}.
For instance, at first glance the identity
$$z_0^{-1}\delta\left({z_1-z_2\over z_0}\right)-z_0^{-1}\delta\left({z_2-
z_1\over -z_0}\right)=z_2^{-1}\delta\left({z_1-z_0\over z_2}\right)$$
is nonsense; it only holds if you expand the terms in
positive powers of $z_2$, $z_1$, and $z_0$ respectively. The procedure
of expanding a function in positive and negative powers of a variable
and then subtracting the results, yields what are called {\it expansions of zero}; it
is possible to show that expansions of zero will always be linear combinations
of Dirac deltas $\delta(az)$ and their various derivatives $\delta^{(k)}(az)$,
as we saw in Paradox 2.

\bigskip\bigskip\noindent{{\smcap 2.6. Ingredient \#3: Vertex Operator
Algebras}}\bigskip

We are now prepared to introduce the important new structure called
vertex operator algebras (VOAs).  They are essentially the chiral algebras
of RCFTs --- see [26,27] for excellent motivation of the 7 axioms below.
A more detailed treatment of the basic theory of VOAs is provided by e.g.\
[23,39,31]. Although VOAs are natural from the CFT perspective and
appear to be an important and rapidly developing area in math, their
definition is not easy: Borcherds is known to have said that you
either know what they are, or you don't want to know.

A VOA is a (infinite-dimensional) graded vector space $V=\oplus_{n\in\Z}V_n$
with infinitely
many bilinear products $u*_nv$ respecting the grading (in particular
$V_k*_n V_\ell\subseteq V_{k+\ell-n-1}$), which obey
infinitely many constraints. `Bilinear' means that for any $a,a',b,b'\in
\C$ and $u,u',v,v'\in V$, $(au+a'u')*_n(bv+b'v')
=ab\,u*_nv+ab'\,u*_nv'+a'b\,u'*_nv+a'b'\,u'*_nv'$ --- i.e.\ that the products
are compatible with the vector space structure of $V$. The subspaces
$V_n$ must all be finite-dimensional, and they must be trivial (i.e.\
$V_n=\{0\}$) for all sufficiently small $n$ (i.e.\ for $n\approx -\infty$).
Note that we can collect all these products into one generating
function: a linear map $Y:V\rightarrow({\rm End} V)[[z,z^{-1}]]$. That is,
to each vector $u\in V$ we associate the formal power series (called 
a {\it vertex operator}) $Y(u,z)=\sum_{n\in\Z}u_nz^{-n-1}$. For each $u$, the
coefficients $u_n$ will be functions from $V$ to $V$. The idea is
that the product $u*_nv$ will now be written $u_nv:=u_n(v)$. The bilinearity
of $*_n$ translates into two things in this new language: that $Y(\star,z)$ is
linear, and that each function $u_n$ is itself linear (i.e.\ they are 
endomorphisms).

The constraints are:\smallskip

\item{{\bf VOA 1.}} ({\it regularity}) $u_nv=0$ for all $n>N(u,v)$;

\item{{\bf VOA 2.}} ({\it vacuum}) there is a vector $1\in V$ such that $Y(1,z)$
is the identity (i.e.\ $1_nv=\delta_{n,-1}v$);

\item{{\bf VOA 3.}} ({\it state-field correspondence}) $Y(u,0)1=u$;

\item{{\bf VOA 4.}} ({\it conformal}) there is a vector $\omega\in V$, called
the {\it conformal vector}, such that $L_{n}:=\omega_{n+1}$ gives us a
representation of the {Virasoro algebra} ${\cal V}$, with central term
$C\mapsto cI$ for some $c\in\C$;

\item{{\bf VOA 5.}} ({\it translation generator}) $Y(L_{-1}u,z)={d\over dz}Y(u,z)$;

\item{{\bf VOA 6.}} ({\it conformal weight}) $L_0u=nu$ whenever $u\in V_n$;

\item{{\bf VOA 7.}} ({\it locality}) $(z-w)^M[Y(u,z),Y(v,w)]=0$ for some integer
$M=M(u,v)$.\smallskip

We saw the Virasoro algebra in Part 1 (see (1.2.7)). The
number $c$ in {\bf VOA 4} is called the {\it central
charge}, and is an important invariant of ${V}$. The peculiar-looking
{\bf VOA 7} simply says that the commutator
$[Y(u,z),Y(v,w)]$ of two vertex operators will be a finite linear
combination of derivatives
of various orders of the Dirac delta centred at $z-w$. A recommended
exercise for the reader is to show that $M=4$ works in {\bf VOA 7} for $u=v=
\omega$. Note that in a VOA,
 any $Y(u,z)v$ will be a finite sum --- i.e.\ the series $Y(u,z)$ is
summable (defined last section). It is a consequence of the axioms that
 $1\in V_0$ and $\omega\in V_2$: for instance, {\bf VOA 7} says all $u_n1=0$
 for any  $n\ge 0$, so $L_01=\omega_11=0$ and hence $1\in V_0$.

In RCFT, $V$ would be the `Hilbert space of states' (more carefully,
$V$ will be a dense subspace of it), and $z=e^{t+\i x}$ would be a local
complex coordinate on a Riemann surface. $L_0$ generates
time translations, and so its eigenvalues (the {\it conformal weights})
can be identified with energy.
Physically, the requirement that $V_n\rightarrow 0$ for $n\rightarrow
-\infty$ corresponds to
 the energy of the RCFT being bounded from below. Also, $z=0$ in {\bf VOA 3}
 corresponds to the time limit $t\rightarrow-\infty$. For each state $u$,
the vertex operator $Y(u,z)$ is a holomorphic (chiral) quantum field.
 The vector 1 is the
 vacuum $|0\rangle$, and $Y(\omega,z)$ is the stress-energy tensor $T$.
The most important axiom, {\bf VOA 7}, says that vertex operators commute up
to a possible pole at
 $z=w$, and so are local quantum fields. It is equivalent to the duality
 axiom of many treatments of CFT.
In the physics literature, there is a minor notational difference: for
$u\in V_k$, $Y(u,z)=\sum u_nz^{-n-1}$ is written $\sum
u_{(n)}z^{-n-k}$. (Physicists prefer this because it cleans up some formulas
a little; mathematicians abhor it because it artificially prefers the
`homogeneous' vectors $u\in V_k$.)

In Segal's language (see \S1.1), $Y(u,z)$ appears quite naturally. Consider
the physical event of two strings combining to form a third. To first order
(i.e.\ the tree-level Feynman diagram), this would correspond in Segal's
language to a `pair-of-pants', or a sphere with three punctures, two of
which are positively oriented (corresponding to the incoming strings) and
the other being negatively oriented. We can think of the sphere as the
Riemann sphere $\C\cup\{\infty\}$; put the punctures at $\infty$ (outgoing)
and $z$ and 0 (the incomings). Segal's functor ${\cal T}$ will associate to this a
$z$-dependent homomorphism $\varphi_z:V\times V\rightarrow V$. We write $
\varphi_z(u,v)\in V$ as $Y(u,z)v$. Incidentally, the symbol `$Y$' was
chosen because of this `pair-of-pants' picture (time flows from the
top of the `$Y$' to the bottom), as was the name `vertex
operator'.

 The original axioms by Borcherds were a little more complicated and
 general:  he didn't require dim$(V_n)<\infty$ nor the
 $V_n\rightarrow 0$ condition, and he only considered $L_0$ and $L_{-1}$ rather
than the full Virasoro algebra. The resulting generalisation is called
 a {\it vertex algebra}.

{\bf VOA 7} can be rewritten in the form (usually called the {\it Jacobi
identity} for the VOA)
$$\eqalignno{z_0^{-1}\delta({z_1-z_2\over z_0})\,Y(u,z_1)\,Y(v,z_2)-z_0^{-1}\delta({
z_2-z_1\over -z_0})&\,Y(v,z_2)\,Y(u,z_1)&(2.6.1)\cr =\,&z_2^{-1}
\delta({z_1-z_0\over z_2})\,Y(Y(u,z_0)v,z_2)\ ,&\cr}$$
where the formal series are expanded in the appropriate way. This is the
embodiment of commutativity and associativity in the VOA, as we will see.
To bring it into a more useful form, hit it with $t\in V$ and expand out
into  $z_0^\ell z_1^mz_2^n$: we obtain
$$\sum_{i\ge 0}(-1)^i\left({\ell\atop i}\right)\,(u_{\ell+m-i}\circ v_{n+i}
-(-1)^\ell v_{\ell+n-i}\circ u_{m+i})=\sum_{i\ge 0}\left({m\atop i}\right)\,
(u_{\ell+i}v)_{m+n-i}\ ,\eqno(2.6.2)$$
where for any $k\in\Z$, $j\in\Z_\ge$, $\left({k\atop j}\right):=
{k\,(k-1)\cdots(k-j+1)\over j!}$.
For instance, specialising (2.6.2) to $\ell=0$ and $m=0$, resp., gives us
$$\eqalignno{[u_m,v_n]=&\,\sum_{i\ge 0}\left({m\atop i}\right)\,(u_iv)_{m+n-i}
&(2.6.3)\cr (u_\ell v)_n=&\,\sum_{i\ge 0}(-1)^i\left({\ell\atop i}\right)\,
(u_{\ell-i}\circ v_{n+i}-(-1)^\ell v_{\ell+n-i}\circ u_i)\ .&(2.6.4)\cr}$$

Why is (2.6.1) called the Jacobi identity? Put $\ell=m=n=0$ in (2.6.2):
we get
$u_0(v_0t)-v_0(u_0t)=(u_0v)_0t$. If we now formally write $[xy]:=x_0y$,
then this becomes $[u[vt]]-[v[ut]]=[[uv]t]$, which is one of the forms
of the Lie algebra Jacobi identity (1.2.1b).
Even though $[xy]\ne -[yx]$ here, this
formal little trick will turn out to be quite important next section.
            
The simplest examples of  VOAs correspond to any even positive-definite
lattice $\L$; for their construction see e.g.\ [27,39]. Physically, they correspond to
a bosonic string compactified on the torus $\R^n/\L\cong S^1\times\cdots
\times S^1$ (where $n$ is the dimension of $\L$); the central charge
$c=n$. Other important examples, first
constructed by Frenkel-Zhu (again see e.g.\ [27,39]), correspond to affine
Kac-Moody algebras $X_\ell^{(1)}$ at level $k\in\Z_>$, and physically to WZW theories
on simply-connected compact group manifolds. (We discussed affine
algebras in \S 1.4.) These have central charge $c={k\,{\rm dim}(X_\ell)
\over k+h^\vee}$.

In 1984 Frenkel-Lepowsky-Meurman [23] constructed the {\it Moonshine module}
$V^{\natural}$. It is a VOA with $c=24$, with
$V^{\natural}=V^\natural_0\oplus V^\natural_1\oplus V^\natural_2\oplus\cdots$,
where $V^\natural_0=\C 1$ is 1-dimensional, $V^\natural_1=\{0\}$ is trivial, and
$V^\natural_2=(\C\omega)\oplus({\rm Griess\ algebra})$ is $(1+196883)$-dimensional.
 Its automorphism
group (=symmetry group) is precisely the Monster  $\M$.
Each graded piece $V^\natural_n$ is a finite-dimensional representation of $\M$;
Borcherds proved that in fact $V^{\natural}$ is the McKay-Thompson
infinite-dimensional representation
of $\M$. It can  be regarded as the most natural representation of
$\M$ --- it is rather surprising that important aspects of a finite group need
 to be studied via an infinite-dimensional representation.

$V^{\natural}$ has an elegant physical interpretation. First
construct the bosonic string on $\R^{24}/\L_{24}$ (recall that $\L_{24}$
is the Leech lattice). The resulting $c=24$ VOA has partition function
(=graded dimension) $J(\tau)+24$, but although its graded pieces (at
least for $n>0$) have the right dimensions, they don't carry a
natural representation of $\M$ and so can't qualify for the
McKay-Thompson representation. To get $V^{\natural}$, orbifold this
$\L_{24}$ VOA by the order-2 automorphism of $\L_{24}$ sending 
$\vec{x}\mapsto-\vec{x}$.
$V^{\natural}$ thus corresponds to a holomorphic $c=24$ RCFT, and Moonshine
is related to physics. Most of Moonshine can be
interpreted physically, except perhaps  the 
 genus 0 property of the McKay-Thompson series $T_g$.

There is a formal parallel between e.g.\ lattices and VOAs.
For example,
the Leech lattice $\L_{24}$ and the Moonshine module $V^{\natural}$ play
analogous roles: $\L_{24}$ is the unique even lattice which (i) is self-dual,
(ii) contains
no norm 2 vectors, and (iii) has dimension 24; $V^{\natural}$ is
believed to be the unique VOA
which (i) possesses only one irreducible representation (namely itself), (ii)
contains no conformal weight 1 elements, and (iii) has central charge $c=24$.
Analogies of these kinds
are always useful as they suggest new directions to explore, and the
history of math blooms with them. The battlecry `Why invent what
can be profitably copied' is not only heard in Hollywood.

We will end this section on a more speculative note. Witten (1986) said that to
understand string theory conceptually, we need a new analogue of Riemannian
geometry. Huang (1997) has pushed this thought a little further, saying that
there is a more classical `particle-math' and a more modern `string-math'.
According to Huang we  have the real numbers (particle physics) vrs
the complex numbers (string theory); Lie algebras vrs VOAs; and the representation
theory of Lie algebras vrs RCFT, etc. What are the stringy analogues of
calculus, ordinary differential equations, Riemannian manifolds, the
Atiyah-Singer Index theorem,...? At present these are all unknown. However,
Huang suggests that just as we could imagine Moonshine
as a mystery which is explained in some way by RCFT, perhaps
 the stringy version of calculus would similarly explain the mystery of 2-dimensional
gravity, stringy ODEs would explain the mystery of infinite-dimensional
integrable systems, stringy Riemannian manifolds would help explain the
mystery of mirror symmetry, and the stringy index theorem would help
explain the elliptic genus.

\bigskip\bigskip\noindent{{\smcap 2.7. Ingredient \#4: Generalised
Kac-Moody algebras}}\bigskip

In this section we investigate Lie algebras arising from VOAs. These
Lie algebras are an interesting generalisation of Kac-Moody algebras.
See e.g.\ [5,6,38 Chapter 11.13, 29].

Much of Lie theory (indeed much of algebra) is developed by analogy with
simple properties of integers. In \S 2.2 I invited you to think of
a finite group as a massive generalisation
of the concept of whole number. Specifically, the number $n$ can be
identified with the cyclic group $\Z_n$ with $n$ elements. A divisor
$d$ of $n$ generalises to a normal subgroup of a group.
A prime number then corresponds to a simple group. Multiplying numbers
corresponds to taking the semidirect product of groups (more
generally, taking  extensions of groups). Then we find that
every group has a unique set of simple building blocks (although unlike
numbers, different groups can have the same list of building blocks).

For a finite-dimensional Lie algebra, a divisor is called an ideal; a prime is called simple;
and multiplying corresponds to semidirect sum. Lie algebras behave simpler
than groups but not as simple as numbers, and the analogy sketched above
is a reasonably satisfactory one. In particular, simple Lie algebras are
important for similar reasons that simple groups are, and as mentioned
in \S 1.2 can also be classified (with {\it much} less effort). 
A good treatment of this important classification (over $\C$) is provided by
[36].
The proof is now reaching the state of perfection of the formulation
of classical mechanics. One unobvious discovery is that the best way to
capture the structure of a simple Lie algebra is by an integer matrix, called
the {\it Cartan matrix}, or equivalently but more effectively (since
most entries in the Cartan matrix are 0's) by using a graph called the
{\it (Coxeter-)Dynkin diagram}. For instance the Dynkin diagram for $A_\ell$
consists of  $\ell$ nodes connected sequentially in a line. See Figure
6 in [59].

More precisely, define a {\it symmetrised Cartan matrix} to be a symmetric
real matrix $A=(a_{ij})_{i,j\le\ell}$ such that $a_{ij}\le 0$ if $i\ne j$,
$a_{ii}>0$, each $2{a_{ij}\over a_{ii}}\in\Z$, and $A$ is positive-definite.
Examples of $2\times 2$ symmetrised Cartan matrices are\footnote{$^{32}$}{{\smal Note
that our Cartan matrices differ from the usual definition, in which every
diagonal entry equals 2.}}
$$\left(\matrix{2&-1\cr -1&2}\right)\ ,\quad
\left(\matrix{2&0\cr 0&2}\right)\ ,\quad
\left(\matrix{1&-1\cr -1&2}\right)\ ,\quad
\left(\matrix{2&-3\cr -3&6}\right)$$
The Dynkin diagram corresponding to $A$ consists of $\ell$ nodes; the $i$th
and $j$th nodes are connected with $4a_{ij}^2/a_{ii}a_{jj}$ lines, and
if $a_{ii}\ne a_{jj}$, then we put an arrow over those lines pointing
to $i$ if $a_{ii}<a_{jj}$. The Dynkin diagrams corresponding
to those four Cartan matrices are respectively
$$\circ\!\!-\!\!\circ\qquad ,\qquad \circ\ \circ\qquad ,\qquad \circ\!\!=
\!\!\!\!<\!\!\!\circ\qquad ,\qquad \circ\!\!\equiv\!\!\!\!<\!\!\!\circ$$
We may
without loss of generality require $A$ to be {\it indecomposable},
or equivalently that the Dynkin diagram be connected. Of the 4 given above,
only the second is decomposable.

To any $\ell\times \ell$ symmetrisable Cartan matrix, we can construct the
corresponding Lie
algebra ${\frak g}$ in the following way. For each $i$, create 3 generators $e_i,f_i,
h_i$ (so there are a total of $3\ell$ generators). The relations these
generators obey are given by the following brackets: $[e_if_j]=\delta_{ij}
h_i$, $[h_ie_j]=a_{ij}e_j$, $[h_if_j]=-a_{ij}f_j$, 
and for $i\ne j$
${\rm ad}(e_i)^ne_j={\rm ad}(f_i)^nf_j=0$ where $n=1-2{a_{ij}\over a_{ii}}$.
By `${\rm ad}(e)$' here I mean the function ${\frak g}\rightarrow {\frak g}$ defined by
${\rm ad}(e)f=[ef]$. So ${\rm ad}(e)^2f=[e[ef]]$, ${\rm ad}(e)^3f=[e[e[ef]]]$, etc.

To get a better feeling for these relations, consider a fixed $i$.
The generators $e=\sqrt{{2\over a_{ii}}}\,e_i,f=\sqrt{{2\over a_{ii}}}\,f_i,
h={2\over a_{ii}}\,h_i$ obey the relations (1.2.2b).
In other words, every node in the Dynkin diagram
corresponds to a copy of the $A_1$ Lie algebra. The lines connecting
these nodes tells how these $\ell$ copies of $A_1$ intertwine.

For instance consider the first Cartan matrix given above. It corresponds
to the Lie algebra $A_2$, or sl$_3(\C)$. The two $A_1$ subalgebras which
generate it (corresponding to the 2 nodes of the Dynkin diagram)
can be chosen to be the trace-zero matrices of the form
$$\left(\matrix{\star&\star&0\cr \star&\star&0\cr 0&0&0\cr}\right)\ ,\qquad
\left(\matrix{0&0&0\cr 0&\star&\star\cr 0&\star&\star\cr}\right)\ .$$

It can be shown that the Lie algebra corresponding to an
indecomposable symmetrised Cartan matrix will be finite-dimensional
and simple, and conversely that any finite-dimensional simple Lie algebra
corresponds to an indecomposable symmetrisable Cartan matrix in this way.

A confusion sometimes arises between the terms `generators' and `basis'.
Both generators and basis vectors build up the whole
algebra; the difference lies in which operations you are  permitted to
use. For a basis, you are only allowed to use linear combinations (i.e.\
addition of vectors and multiplication by numbers), while for generators
you are also permitted multiplication of vectors (or the bracket, in the
Lie case). `Dimension' refers to basis, while `rank' usually refers in
some way to generators. For instance the (commutative associative)
algebra  of polynomials
in one variable $x$ is infinite-dimensional --- any basis needs infinitely
many vectors. However, the single polynomial $x$ is enough to generate it
(so we could say that its rank is 1).
Although those Lie algebras have $3\ell$ generators, their dimensions
in general will be greater.

 From the point of view of generators and relations,
the step from `finite-dimensional simple' to `symmetrisable Kac-Moody' is rather
easy: the only difference is that we drop the `positive-definite' condition 
(which was responsible for finite-dimensionality).
Kac-Moody (KM) algebras are also generated by (finitely many) $A_1$ subalgebras,
and their theory is quite parallel to that of the simple algebras.
Compare Figures 6 and 9 in [59].

Now, it is easy
to generalise something; the challenge is to generalise it in a rich and
interesting direction. One natural and appealing strategy for generalisation
was followed instinctively by a grad student named Robert Moody.
 Moody's original motivation for developing the theory of Kac-Moody
algebras was the Weyl group. If there were Lie
algebras for finite Coxeter groups, he asked, why not also for the Euclidean
(=affine) ones?  For another example of this style of generalisation, consider
the question: What is the analogue of calculus (or manifolds) over
weird fields --- fields (like $\Z_p$) for which the usual limit
definitions make no sense? This question leads to the riches of
algebraic geometry. Nevertheless this generalisation strategy, even in the
hands of a master, will not always be successful. For instance, consider all
the trouble the following metaphor has caused: my watch has a maker, so so
should the Universe.

Recently Borcherds produced a further generalisation of finite-dimensional
simple Lie algebras, which is rather less obvious than that of Kac-Moody
algebras. It is easy to associate a Lie algebra to a
matrix $A$, but which class of matrices will yield a
deep theory? Borcherds found such a class by holding in his hand
a single algebra (the fake Monster Lie algebra, see \S 2.9) which acted a lot
 like a KM algebra, even though it had `imaginary simple roots'.

By a {\it generalised symmetrised Cartan matrix} $A=(a_{ij})$ we will mean a symmetric
real matrix (possibly infinite), such that $a_{ij}\le 0$ if $i\ne j$,
and if $a_{ii}>0$ then $2{a_{ij}\over a_{ii}}\in\Z$ for all $j$. By a
{\it universal generalised Kac-Moody algebra} (universal GKM) or {\it
universal
Borcherds-Kac-Moody algebra} ${\frak g}$ we mean the algebra\footnote{$^{33}$}{{
\smal As with KM algebras, usually we want to extend it by some derivations;
enough derivations are added so that the simple roots are linearly independent.}}
 with generators $e_i,f_i,
h_{ij}$, and with relations: $[e_if_j]=h_{ij}$; $[h_{ij}e_k]=\delta_{ij}
a_{ik}e_k$; $[h_{ij}f_k]=-\delta_{ij}a_{ik}f_k$; if $a_{ii}>0$ and $i\ne j$
then ${\rm ad}(e_i)^ne_j={\rm ad}(f_i)^nf_j=0$ for $n=1-2{a_{ij}\over a_{ii}}$;
and if $a_{ij}=0$ then $[e_ie_j]=[f_if_j]=0$.

For example the Heisenberg algebra (1.2.3) corresponds to the choice $A=(0)$,
while any other $1\times 1$ $A=(a)$ corresponds to $A_1$.
A universal GKM algebra differs from a KM algebra
in that it is built up from Heisenberg algebras as well as $A_1$, and
these subalgebras intertwine in more complicated ways. Nevertheless
{\it much of
the theory for finite-dimensional simple Lie algebras continues to find
an analogue in this much more general setting} (e.g.\ root-space decomposition,
Weyl group, character formula,...). This unexpected fact is the point
of GKM algebras.

To get a feel for these algebras, let us prove a few simple results
concerning the $h_{ij}$. Note first that, using the above relations together
with the Jacobi identity, we obtain
$[h_{ij}h_{k\ell}]=\delta_{ij}(a_{jk}-a_{j\ell})h_{k\ell}$. Comparing this
with $[h_{k\ell}h_{ij}]=-[h_{ij}h_{k\ell}]$, we
see that bracket must always equal 0. Hence all $h$'s pairwise commute,
and $h_{ij}=0$ unless the $i$th and $j$th columns of $A$ are identical.
An easy exercise now is to show that when $i\ne j$,
$h_{ij}$
will lie in the centre of the algebra (i.e.\ $h_{ij}$ will commute with
all other generators).

Although the definition of universal GKM algebra is more
natural, it turns out that an equivalent form can be more useful in practice.
(It's simpler to describe over $\R$, so in most
expositions the reals are used, but alas it's far too late for us to
switch loyalties now.)
By a {\it generalised Kac-Moody algebra} (or {\it Borcherds-Kac-Moody
 algebra}) ${\frak g}$, we mean a (complex!) Lie algebra which is: \smallskip

\item{$-$} $\Z$-graded, i.e.\ ${\frak g}=\oplus_i {\frak g}_i$, where $[{\frak g}_i
{\frak g}_j]\subseteq
{\frak g}_{i+j}$;

\item{$-$} ${\frak g}_i$ is finite-dimensional for $i\ne 0$;

\item{$-$} ${\frak g}$ has an antilinear involution $\omega$ (i.e.\ $\omega(
kx+y)=k^*\omega(x)+\omega(y)$,
$[\omega(x)\omega(y)]=\omega([xy])$, and $\omega\circ\omega=id.)$ which maps
${\frak g}_i$
to ${\frak g}_{-i}$ and acts as multiplication by $-1$ on some basis of
${\frak g}_0$;

\item{$-$} ${\frak g}$ has an invariant symmetric bilinear form $(\star,\star)$
(i.e.\
$([xy],z)=(x,[yz])$ and $(y,z)=(z,y)\in\C$), obeying
$(\omega(x),\omega(y))=(x,y)^*$, such that $({\frak g}_i,
{\frak g}_j)=0$ if
$i\ne-j$;

\item{$-$} the Hermitian form defined by $(x|y):=-(\omega(x),y)$ is
positive-definite on ${\frak g}_{i\ne 0}$.\smallskip

Note that  for some basis $x_i$ of ${\frak g}_0$, the third condition tells us
$-[x_ix_j]=[(-x_i)(-x_j)]=[x_ix_j]$, i.e.\ ${\frak g}_0$ has a trivial bracket.
It plays the role of the Cartan subalgebra ${\frak h}$ in the theory.

For example, let ${\frak g}={\rm sl}_2(\C)$ and recall (1.2.5). Then
${\frak g}_1=\C e$, ${\frak g}_0=\C h$, ${\frak g}_{-1}=\C f$ is the root-space
decomposition. $\omega(x)=-x^{\dag}$ is the Cartan involution. $(x,y)=
{\rm tr}(xy)$ is the Killing form.

It turns out [5] that any universal GKM algebra is a GKM algebra, and any
GKM algebra can be constructed from a unique universal GKM algebra (by
quotienting out part of the centre and adding derivations), so
in that sense the two structures are equivalent. This theorem is
important, because it tells us that {\it GKM algebras are the ultimate
generalisation
of simple Lie algebras,} in the sense that any further generalisation will lose
some basic structural ingredient.

We know simple Lie algebras (and groups) arise in both classical and
quantum physics, and the affine KM algebras are important in CFT, as
we saw in Part 1. GKM algebras have recently appeared in the physics literature
(see Harvey-Moore) in the context of BPS states in string theory.

How do GKM algebras arise in VOAs? If we define $[xy]:=
x_0y$, then as mentioned in \S 2.6 we get from the VOA Jacobi identity
the equation $[x[yz]]-[y[xz]]=[[xy]z]$. Thus our bracket will be anti-associative
if it is anti-commutative. But is it anti-commutative? It can be shown
$$u_nv=\sum_{i=0}^\infty {1\over i!}(-1)^{i+n+1}(L_{-1})^i
(v_{n+i}u)\eqno(2.7.1)$$
so $u_0v\equiv -v_0u$ if we look at things mod $L_{-1}V$.

Since our bracket is clearly bilinear, we thus get a Lie algebra structure
on $V/L_{-1}V$.
Similarly, we get a symmetric bilinear product on $V/L_{-1}V$, given by
$\langle u,v\rangle:=u_1v$.

We would like $\langle\star,\star\rangle$ to respect the Lie algebra structure,
i.e.\ be $[\star\star]$-invariant. We compute from (2.6.4)
$$\langle[uv],t\rangle=-[v\langle u,t\rangle]+\langle u,[vt]\rangle\ .
\eqno(2.7.2)$$
Since we would like  to identify $\langle\star,\star\rangle$
with the bilinear form in the GKM algebra definition, we also would like
it to be number-valued (i.e.\ have 1-dimensional range).

There is a simple way to satisfy both of these. First, restrict attention
to $V_1$, i.e.\ the conformal weight 1 vectors: $V_1\cap(V/L_{-1}V)=
V_1/L_{-1}V_0$. Then $\langle u,v\rangle\in V_0$. Assume $V_0$ is
1-dimensional: i.e.\ $V_0=\C 1$. Then $\langle u,v\rangle$ will equal
a number times $1$, so call $(u,v)$ that number. Also, $[\langle u,t\rangle
v]=(u,t)1_0v=0$, so $\langle\star,\star\rangle$ and hence $(\star,\star)$, 
will be invariant. Of course, when $V_0=\C 1$, $L_{-1}V_0=\{0\}$.

$V_1/L_{-1}V_0$ is generally too large in practice to be useful; a subalgebra can be
defined as follows. Let $P_n$ be the `primary states with conformal weight
$n$', i.e.\ the $u\in V_n$ killed by $L_m$ for all $m>0$. Then ${\frak g}(V)
:=P_1/L_{-1}P_0$ will be a subalgebra of $V_1/L_{-1}V_0$. Through the
assignment $u\mapsto u_0$, ${\frak g}(V)$ acts on $V$ and this action
commutes with that of $L_i$. This association
of a Lie algebra to a VOA is due to Borcherds (1986).

Similar arguments show that when $V_0$ is 1-dimensional and $V_1$ is
 0-dimensional, then $V_2$ will necessarily be a commutative nonassociative
  algebra
with product $u\times v:=u_1v\in V_2$ and identity element ${1\over 2}\omega$
({\it proof}: $\omega\times u=L_0u=2u$).
Now, those conditions on $V_0,V_1$ are satisfied by the Moonshine module
$V^{\natural}$. We find that  $V^{\natural}_2$ is none other than the
Griess algebra extended by an identity element.

\vfill\eject\noindent{{\smcap 2.8. Ingredient \#5: Denominator identities}}
\bigskip

In \S 1.3, we discussed the representation theory of Lie
algebras. An important invariant of a representation is its {\it character}.
Simple Lie algebras possess a very useful formula for their characters,
due to Weyl:
$${\rm ch}_\la(z):=\sum_{\mu}{\rm dim}(L_\la(\mu))\,e^{\mu\cdot
z}=e^{-\rho\cdot z}\,{\sum_{w\in W}\pm e^{w(\la+\rho)\cdot z}\over 
\prod_{\alpha\in\Delta_+} (1-e^{-\alpha\cdot z})}\ ,\eqno(2.8.1)$$
where $W$ is the Weyl group, $\Delta_+$ the positive roots, and where
$\oplus_\mu L_\la(\mu)$ is the weight-space decomposition of $L_\la$ ---
i.e.\ the simultaneous eigenspaces of the $h_i$.
Here $z$ belongs to the Cartan subalgebra ${\frak h}$; the character is complex-valued.
Analogous statements hold for all GKM algebras.

It is rare indeed when a trivial special case of a theorem or formula
produces something interesting. But that is what happens here. Consider
the trivial representation: i.e.\ $x\mapsto 0$ for all $x\in X_\ell$. Then the
character is identically 1, by definition: $ch_0\equiv 1$. Thus the
character formula tells us that a certain alternating  sum over a Weyl group, equals
a certain product over positive roots. These formulas, called {\it
denominator formulas}, are nontrivial even in the finite-dimensional cases.

Consider for instance the smallest simple algebra, $A_1$. Here the
identity indeed is too trivial: it reads $e^{z/2}-e^{-z/2}=e^{z/2}
(1-e^{-z})$. For $A_2$ we get a sum of 6 terms equalling a product of
3 terms, and the complexity continues to rise from there.

Around 1970 Macdonald tried to generalise these finite denominator
identities to infinite identities, corresponding to the extended Dynkin diagrams.
These were later reinterpreted by Kac, Moody and others as denominator
identities for affine nontwisted KM algebras. The simplest one was known
classically as the Jacobi triple product identity:
$$\sum_{n=-\infty}^\infty(-1)^nx^{n^2}y^n=
\prod_{m=1}^\infty(1-x^{2m})(1-x^{2m-1}y)(1-x^{2m-1}y^{-1})$$
We now know it to be the denominator identity for the simplest infinite-dimensional
KM algebra, $A_1^{(1)}$.

Freeman Dyson is famous for his work in quantum field theory, but he started
as an undergraduate in  number theory and
still enjoys it as a hobby. Dyson [20] found a curious formula for
the Ramanujan $\tau$-function, which can be defined by the generating
function $\sum_{n=1}^\infty
\tau(n)x^n=\eta^{24}(x):=x\prod_{m=1}^\infty(1-x^m)^{24}$.  Dyson found the
remarkable formula
$$\tau(n)=\sum{(a-b)(a-c)(a-d)(a-e)(b-c)(b-d)(b-e)(c-d)(c-e)(d-e)\over
1!\ 2!\ 3!\ 4!}$$
where the sum is over all 5-tuples $(a,b,c,d,e)\equiv(1,2,3,4,5)$ (mod 5)
obeying $a+b+c+d+e=0$ and $a^2+b^2+c^2+d^2+e^2=10n$. Using this, an
analogous formula can be found for $\eta^{24}$. Dyson knew that
similar-looking formulas were also known for $\eta^d$ for the values
$d=3,8,10,14,15,21,24,26,28,35,36,\ldots$.

What was ironic was that Dyson found that formula at the same time that
Macdonald was finding the Macdonald identities. Both were at Princeton
then, and would often chat a little when they bumped into each other
after dropping off their daughters at school. But they never discussed
work. Dyson didn't realise that his strange list of numbers has a
simple interpretation: they are precisely the dimensions of the 
simple Lie algebras! $3={\rm dim(A_1)}$, $8={\rm dim}(A_2)$, $10={\rm dim}(
C_2)$, $14={\rm dim}(G_2)$, etc. In fact these formulas for $\eta^d$
are none other
than (specialisations of) the Macdonald identities. For example, Dyson's formula is the
denominator formula for $A_4^{(1)}$ ($24={\rm dim}(A_4)$). If they had
spoken, they would probably have anticipated the affine denominator identity
interpretation.

One curiousity apparently still has no algebraic 
interpretation: No
simple Lie algebra has dimension 26, so the formula for $\eta^{26}$ 
can't correspond to any Macdonald identity.

Macdonald didn't close the book on denominator identities. More recently
Kac and Wakimoto [40] have used denominator identities for Lie superalgebras
 to obtain nice formulas for various generator functions
involving sums of squares, sums of triangular numbers (triangular
numbers are numbers of the form ${1\over 2}k(k+1)$), etc. For instance,
the number of ways $n$ can be written as a sum of 16 triangular numbers
is
$${1\over 3\cdot 4^3}\sum ab\,(a^2-b^2)^2$$
where the sum is over all odd positive integers $a,b,r,s$ obeying
$ar+bs=2n+4$ and $a>b$. 

Another example of denominator identities is Borcherds' use of them in
proving the Moonshine conjectures. In particular this motivated his
introduction of the GKM algebras.
The denominator identities for other GKM algebras were used by Borcherds
to obtain results on the moduli spaces of e.g.\ families of K3 surfaces.
They are also often turned-around now and used for learning about the
positive roots in a given GKM.

\bigskip\bigskip\noindent{{2.9. \smcap Proof of the Moonshine conjectures}}
\bigskip

The main Conway-Norton conjecture was proved almost immediately. Thompson
showed that if $g\mapsto a_n(g)$ is a character for all $n\le 1200$, then
it will be for all $n$. He also showed that if certain congruence conditions
hold for a certain number of $a_n(g)$ (all with $n\le 100$), then all
$g\mapsto a_n(g)$ will be {\it virtual} characters (i.e.\ a linear combination
over $\Z$ of irreducible characters of $\M$; only if all coefficients
are nonnegative will it be a true character).
Atkin-Fong-Smith [54] used that to prove on a computer
that indeed all were virtual characters. But their work didn't say
anything about the underlying (possibly virtual) representation $V$. The real
challenge was to construct (preferably in a natural manner) the representation
which works. Frenkel-Lepowsky-Meurman [23] constructed a candidate for
it (the Moonshine module $V^\natural$); it was Borcherds who finally
proved $V^\natural$ obeyed the Conway-Norton conjecture. A good overview
of Borcherds' work on Moonshine is provided in [32].

We want to show that the McKay-Thompson series $T_g(\tau):=q^{-1}{\rm tr}_{V^{\natural}}
(gq^{L_0})$ of (2.1.4) equals the Hauptmodul $J_g(\tau)$ in (2.4.1) 
(the `fudge factor' $q^{-1}
=q^{-c/24}$ is familiar to e.g.\ KM algebras and CFT and was discussed
at the end of \S 1.2). Borcherds' strategy was to
bring in Lie theory and to use the corresponding denominator identity
to provide useful combinatorial data. The first guess for this `Monster
Lie algebra' was the Kac-Moody algebra whose Dynkin diagram is
essentially the Leech lattice (i.e.\ a node for each vector in
$\L_{24}$, and 2 nodes are connected by a number of edges depending 
on the value of a certain dot product). It was eventually discarded because
 some of the critical data (namely positive root multiplicities)
needed in order to write down its denominator identity were too complicated.
But looking at that failed attempt led Borcherds to a second candidate, 
now called the {\it fake Monster Lie algebra} ${\frak g}_{\M}'$. In order to
{\it construct} it, he developed the theory of VOAs; and in order to {\it understand}
it, he developed the theory of GKM algebras. We will define it shortly.
 ${\frak g}_{\M}'$ also turned out to be inadequate for proving the Moonshine
conjectures; however it directly led him to the GKM algebra now called
the true {\it Monster Lie algebra} ${\frak g}_{\M}$. And that
directly led to the proof of Moonshine.

\medskip\noindent{{\it Step 1}}: Construct ${\frak g}_{\M}$ from $V^{\natural}
=V^\natural_0\oplus V^\natural_1\oplus\cdots$. For later convenience, 
reparametrise these
subspaces $V^i:=V^\natural_{i+1}$. Recall the even indefinite lattice $II_{1,1}$
defined in (1.6.1). Of course the direct choice ${\frak g}(V^{\natural})$
is 0-dimensional because $V_1^{\natural}$ is trivial, so we must
modify $V^\natural$ first.

The Monster Lie algebra is (essentially) defined to be the Lie algebra 
${\frak g}(V^\natural\otimes V_{II_{1,1}})$ associated to
the vertex algebra $V^{\natural}\otimes V_{II_{1,1}}$ (strictly
speaking, more of ${\frak g}$ is quotiented away). By contrast, the
fake Monster is the Lie algebra associated to the vertex algebra $V_{\L_{24}}
\otimes V_{II_{1,1}}\cong V_{II_{25,1}}$. Both of these are vertex algebras
as opposed to VOAs, because of the presence of the indefinite lattices,
but this isn't important here. ${\frak g}_{\M}$ inherits a $II_{1,1}$-grading
from $V_{II_{1,1}}$: the piece of grading $(m,n)$
is isomorphic (as a vector space) to $V^{mn}$, if $(m,n)\ne (0,0)$;
the (0,0) piece is isomorphic to $\R^2$. Borcherds uses the No-Ghost
Theorem of string theory to show that  the homogeneous pieces of ${\frak g}_{
\M}$ are those of $V^{\natural}$.

Both ${\frak g}_{\M}$ and ${\frak g}'_{\M}$  are GKM algebras; for instance
the $\Z$-grading of ${\frak g}_{\M}$ is given by $({\frak g}_{\M})_k=
\oplus_{m+n=k} V^{mn}$ for $k\ne 0$, while the 0-part is $V^{-1}\oplus
V^{-1}$. 
Although ${\frak g}_{\M}'$ is not used in the proof of the Monstrous Moonshine
conjectures, it is related to some kind of Moonshine for the finite
simple group $.1$, which is `half' of the automorphism group $.0$ of the
Leech lattice $\L_{24}$.
       ${\frak g}_{\M}$ corresponds to the Cartan matrix
$$\left(\matrix{B_{-1,-1}&B_{-1,1}&B_{-1,2}&\cdots\cr B_{1,-1}&B_{1,1}&
B_{1,2}&\cdots\cr B_{2,-1}&B_{2,1}&B_{2,2}&\cdots\cr \vdots&\vdots&\vdots
&\ddots\cr}\right)\ ,$$
where $B_{i,j}$ for $i,j\in\{-1,1,2,3,\ldots\}$ is the $a_i\times a_j$
block with fixed entry $-i-j$ (the $a_i$ as usual are the coefficients
$\sum_n a_nq^n$ of $j-744$).

\medskip
\noindent{{\it 2nd step}}: Compute the denominator identity of ${\frak g}_{\M}$:
we get
$$p^{-1}\prod_{{m>0\atop n\in\Z}}(1-p^mq^n)^{a_{mn}}=j(z)-j(\tau)\eqno(2.9.1)$$
where $p=e^{2\pi \i z}$.
The result are various formulas involving the coefficients $a_i$,
 for instance $a_4=a_3+(a_1^2-a_1)/2$. It turns out to be possible to `twist'
(2.9.1) by each $g\in\M$, obtaining
$$p^{-1}\exp[-\sum_{k>0}\sum_{{m>0\atop
n\in\Z}}a_{mn}(g^k){p^{mk}q^{nk}\over k}]=T_g(z)-T_g(\tau)\ .\eqno(2.9.2)$$
This looks a lot more complicated, but you can glimpse the Taylor
expansion of ln$(1-p^mq^n)$ there and in fact for $g=id$ (2.9.2) reduces
to (2.9.1). This formula gives more generally identities like
$a_4(g)=a_2(g)+(a_1(g)^2-a_1(g^2))/2$, where $T_g(\tau)=\sum_ia_i(g)q^i$.
These formulas involving the McKay-Thompson coefficients are equivalent to
the replication
formulae conjectured in \S 2.4.

\medskip\noindent{{\it 3rd step}}: It was known earlier that all of the
Hauptmoduls also obey those replication formula, and that anything obeying
them will lie in a finite-dimensional manifold which we'll call $R$.
In particular, if $B(q)=q^{-1}+\sum_{n>0}b_nq^n$ and $C(q)=q^{-1}+\sum_{n>0}
c_nq^n$ both lie in $R$, and $b_n=c_n$ for $n\le 23$, then $B(q)=C(q)$. In fact,
it turns out that if we verify for each conjugacy class $[g]$ of $\M$
that the first, second, third, fourth and sixth coefficients of the
McKay-Thompson series $T_g$ and the corresponding Hauptmodul $J_g$ agree,
then $T_g=J_g$, and we are done.

That is precisely what Borcherds then did: he compared finitely many
coefficients, and as they all equalled what they should, this concluded
the proof of Monstrous Moonshine!\medskip

However there was a disappointing side to his proof. While no one disputed
its {\it logical} validity, it did seem to possess a disappointing
{\it conceptual gap}. In particular, the Moonshine conjectures were made
in the hope that proving them would help explain what the Monster had
to do with the $j$-function and the other Hauptmoduls. A good proof
says much more than `True' or `False'. The case-by-case verification
 occurred at the critical point where the McKay-Thompson series
were being compared directly to the Hauptmoduls. The proof showed that
indeed the Moonshine module establishes some sort of relation between
$T_g$ and $J_g$ (namely, they must lie in the
same finite-dimensional space), but  why couldn't it be just a happy
meaningless 
accident that they be equal? Of course we believe it's more than merely an
accident, so our proof should reflect this: we want a more conceptual
explanation.

This conceptual gap has since been filled [17] --- i.e.\ the
case-by-case verification has been replaced with a general theorem.
It turns out that
something obeying the replicable formulas will also obey something called
{\it modular equations}. A modular equation for a function $f$ is a polynomial
identity obeyed by $f(x)$ and $f(nx)$. The simplest examples come from
the exponential and cosine functions: note that for any $n>0$,
$\exp(nx)=(\exp(x))^n$ and $\cos(nx)=T_n(\cos(x))$ where $T_n$ is a
Tchebychev polynomial. A more interesting example of a modular
equation is obeyed by $J(\tau)=j(\tau)-744$: put $X=J(\tau)$ and 
$Y=J(2\tau)$, then
$$\eqalignno{(X^2-Y)(Y^2-X)=&\,393768\,(X^2+Y^2)+42987520\,XY+40491318744\,(X+Y)
&\cr&\,-120981708338256\ .&\cr}$$

Finding modular equations (for various elliptic functions) was a passion
of the great mathematician Ramanujan. His notebooks are filled with them.
See e.g.\ [7] for an application of Ramanujan's modular equations to
computing the first billion or so digits of $\pi$. Many modular equations
are also studied in [10]. For more of their applications, see e.g.\ [16].
It can be shown that the only functions $f(\tau)=q^{-1}+a_1q+\cdots$
which obey modular equations for all $n$, are $J(\tau)$ and the `modular
fictions' $q^{-1}$ and $q^{-1}\pm q$ (which are essentially exp, cos, and
sin).

It was proved in [17] that, roughly speaking, a function $B(\tau)
=q^{-1}+\sum_{n>0}b_nq^n$ which obeys enough modular equations, will
either be of the form $B(\tau)=q^{-1}+b_1q$, or will
necessarily be a Hauptmodul for a modular group containing some
$\Gamma_0(N)$. The converse is also true: for instance, a modular equation
for the Hauptmodul $J_{25}$ of $\G_0(25)$ given in (2.3.5) is
$$(X^2-Y)(Y^2-X)=-2\,(X^2+Y^2)+4\,(X+Y)-4\ ,$$
where $X=J_{25}(\tau)$ and $Y=J_{25}(2\tau)$. To eliminate the conceptual
gap, this result should then replace step 3. Steps 1 and 2 are still
required, however.

This conceptual gap should not take away from what was a remarkable
accomplishment by Borcherds: not only the proof of the Monstrous Moonshine
conjectures, but also the definition of two new and important algebraic
structures. I hope the preceding sections give the reader some indication
of why Borcherds was awarded one of the 1998 Fields Medals.

Another approach to the Hauptmodul property is by Tuite [57], who
related it to the (conjectured) uniqueness of $V^{\natural}$. Norton
 has suggested that the reason $\M$ is associated to {\it
genus-zero} modular functions could be what he calls its
`6-transposition'  property [47].

So has Moonshine been explained? According to Conway, McKay, and many others,
it hasn't. They consider VOAs in general, and $V^\natural$ in
particular, to be too complicated to be God-given. The progress,
though impressive, has broadened not lessened the fundamental mystery,
they would argue.

For what it's worth, I don't completely agree. Explaining away a
mystery is a little like grasping a bar of soap in a bathtub, or
quenching a child's curiousity. Only extreme measures like pulling the
plug, or enrollment in school, ever really work. True progress means
displacing the mystery, usually from the particular to the
general. Why is the sky blue? Because of how light scatters in
gases. Why are Hauptmoduls attached to each $g\in\M$? Because of
$V^\natural$. Mystery exists wherever we can ask `why' --- like beauty it's
in the beholder's eye.

Moonshine is now `leaving the nest'. We are entering a consolation
phase, tidying up, generalising, simplifying, clarifying, working out
more examples. Important and interesting discoveries will be made in the
next few years, and yes there still is mystery, but no longer does
a Moonshiner feel like an illicit distiller: Moonshine is now a day-job!

\bigskip\noindent{{\it Acknowledgments.}}\quad I warmly thank the Feza Gursey
Institute in Istanbul, and in particular Teoman Turgut, for their
invitation  to the workshop and hospitality
during my month-long stay. These notes are based on 16 lectures I gave
there in Summer 1998.
I've also benefitted from numerous conversations with Y.\ Billig,
A.\ Coste, C.\ Cummins,
M.\ Gaberdiel, J.\ McKay, M.\ Tuite, and M.\ Walton --- Mark Walton in
particular made a very careful
reading of the manuscript (and hence must share partial blame for any errors
still remaining). My appreciation as well goes to J.-B.\ Zuber and P.\ Ruelle
for sharing with me their personal stories behind the discoveries of, respectively,
 A-D-E in $A_1^{(1)}$ and Fermat in $A_2^{(1)}$.
The research was supported in part by NSERC.


\vfill\eject\noindent{{\smcap References}} \bigskip

\item{1.} V.\ I.\ Arnold, {\it Catastrophe Theory}, 2nd edn.\
(Springer, Berlin, 1997);

\item{} M.\ Hazewinkel, W.\ Hesselink, D.\ Siersma, and F.D.\ Veldkamp,
{\it Nieuw Arch.\ Wisk.}\ {\bf 25} (1977) 257;

\item{} P.\ Slodowy, in:  Lecture Notes in Math
 1008, J.\ Dolgachev (ed.) (Springer, Berlin, 1983).

\item{2.} {M.\ Bauer, A.\ Coste, C.\ Itzykson, and P.\ Ruelle},
{\it J.\ Geom.\ Phys.}\ {\bf 22} (1997) 134.

\item{3.} {D.\ Bernard}, 
{\it Nucl.\ Phys.} {\bf B288} (1987) 628.

\item{4.} J.\ B\"ockenhauer and D.\ E.\ Evans, {\it Commun.\ Math.\ Phys.}
{\bf 200} (1999) 57; ``Modular invariants, graphs
and $\alpha$-induction for nets of subfactors III'', 
hep-th/9812110.

\item{5.} {R.\ E.\ Borcherds},
{\it Invent.\ math.}\ {\bf 109} (1992) 405.

\item{6.} R.\ E.\ Borcherds, ``What is moonshine?'', math.QA/9809110.

\item{7.} J.\ M.\ Borwein, P.\ B.\ Borwein and D.\ H.\ Bailey,
{\it Amer.\ Math\ Monthly} {\bf 96} (1989) 201.

\item{8.} {A.\ Cappelli, C.\ Itzykson, and J.-B.\ Zuber,}
{\it Commun.\ Math.\ Phys.} {\bf 113} (1987) 1.

\item{9.} R.\ Carter, G.\ Segal and I.\ M.\ Macdonald, {\it Lectures
on Lie Groups and Lie algebras} (Cambridge University Press, Cambridge, 1995).

\item{10.} A.\ Cayley, {\it An Elementary Treatise on Elliptic Functions},
2nd edn (Dover, New York, 1961).
             
\item{11.} {J.\ H.\ Conway}, 
{\it Math.\ Intelligencer} {\bf 2}
(1980) 165.

\item{12.} {J.\ H.\ Conway and S.\ P.\ Norton},
{\it Bull.\ London Math.\ Soc.}\ {\bf 11} (1979) 308.

\item{13.} {J.\ H.\ Conway and N.\ J.\ A.\ Sloane}, {\it Sphere
packings, lattices and groups}, 3rd edn (Springer, Berlin, 1999).

\item{14.} A.\ Coste and T.\ Gannon, {\it Phys.\ Lett.}\ {\bf B323}
(1994) 316.

\item{15.} R.\ Courant and D.\ Hilbert, {\it Methods of Mathematical
Physics II} (Wiley, New York, 1989).

\item{16.} D.\ Cox, {\it Primes of the form} $x^2+ny^2$, (Wiley, New York,
1989).

\item{17.} {C.\ J.\ Cummins and T.\ Gannon}, 
{\it Invent.\ math.} {\bf 129} (1997) 413.

\item{18.} P.\ Di Francesco, P.\ Mathieu and D.\ S\'en\'echal, {\it
Conformal Field Theory} (Springer, New York, 1996).

\item{19.} {L.\ Dixon, P.\ Ginsparg, and J.\ Harvey},
{\it Commun.\ Math.\ Phys.}\ {\bf 119} (1988) 221.

\item{20.} F.\ Dyson, {\it Bull.\ Amer.\ Math.\ Soc.}\ {\bf 78} (1972) 635.

\item{21.} F.\ J.\ Dyson, {\it Math.\ Intelligencer} {\bf 5} (1983) 47.

\item{22.} {D.E.\ Evans and Y.\ Kawahigashi}, {\it Quantum symmetries
on operator algebras} (Oxford University Press, Oxford, 1998).

\item{23.} {I.\ Frenkel, J.\ Lepowsky, and A.\ Meurman}, {\it
Vertex operator algebras and the Monster} (Academic Press, San Diego,
1988).

\item{24.} J.\ Fuchs and C.\ Schweigert, {\it Symmetries, Lie
algebras, and representations} (Cambridge University Press, Cambridge,
1997).

\item{25.} {W.\ Fulton and J.\ Harris}, {\it Representation Theory:
A first course} (Springer, New York, 1996).

\item{26.} M.\ R.\ Gaberdiel and P.\ Goddard, ``Axiomatic conformal field
theory'', hep-th/9810019.

\item{27.} M.\ R.\ Gaberdiel and P.\ Goddard, ``An introduction to
meromorphic conformal field theory and its representations'', lecture
notes in this volume.

\item{28.} D.\ Gaitsgory, ``Notes on two dimensional conformal field
theory and string theory'', math.AG/9811061.

\item{29.} T.\ Gannon, ``The Cappelli-Itzykson-Zuber A-D-E classification'',
math.QA/9902064.

\item{30.} K.\ Gawedzki, ``Conformal field theory: a case study'',
lecture notes in this volume.

\item{31.} {R.\ W.\ Gebert},
{\it Internat.\ J.\ Mod.\ Phys.}\ {\bf A8} (1993) 5441.

\item{32.} P.\ Goddard, ``The work of R.E.\ Borcherds'', math.QA/9808136.

\item{33.} D.\ Gorenstein, {\it Finite Simple Groups} (Plenum, New York,
1982).

\item{34.} D.\ Gorenstein, R.\ Lyons, and R.\ Solomon, {\it The Classification
of the Finite Simple Groups} (AMS, Providence, 1994).

\item{35.} {A.\ Hanany and Y.-H.\ He}, ``Non-abelian finite gauge
theories'', hep-th/9811183.

\item{36.} J.\ E.\ Humphreys, {\it Introduction to Lie algebras and
representation theory} (Springer, New York, 1994).

\item{37.} V.\ G.\ Kac, In: Lecture Notes in Math 848 (Springer, New York,
1981).

\item{38.} {V.\ G.\ Kac}, {\it Infinite Dimensional Lie algebras}, 
3rd edn (Cambridge University Press, Cambridge, 1990).

\item{39.} V.\ G.\ Kac, {\it Vertex Algebras for Beginners}, 2nd edn (AMS, Providence,
1998).

\item{40.} {V.\ G.\ Kac and M.\ Wakimoto}, 
In: {\it Lie Theory
and Geometry in Honor of Bertram Kostant}, Progress in Math.\ {\bf 123}
(Birkh\"auser, Boston, 1994).

\item{41.} S.\ Kass, R.\ V.\ Moody, J.\ Patera and R.\ Slansky, {\it
Affine Lie algebras, weight multiplicities, and branching rules}, Vol.\ 1
(University of California Press, Berkeley, 1990).

\item{42.} S.\ Lang, {\it Elliptic Functions}, 2nd edn (New York,
Springer, 1997).

\item{43.} F.\ W.\ Lawvere and S.\ H.\ Schanuel, {\it Conceptual Mathematics:
A first introduction to Categories} (Cambridge University Press, Cambridge,
1997).

\item{44.} H.\ Minc, {\it Nonnegative matrices} (Wiley, New York, 1988).

\item{45.} {E.J.\ Mlawer, S.G.\ Naculich, H.A.\ Riggs,  and 
H.\ J.\ Schnitzer}, 
{\it Nucl.\ Phys.} {\bf B352} (1991) 863.

\item{46.} W.\ Nahm, {\it Commun.\ Math.\ Phys.}\ {\bf 118} (1988) 171.

\item{47.} S.\ P.\ Norton, In: {\it Proc.\ Symp.\ Pure Math.}\ {\bf 47}
(1987) 208.

\item{48.} A.\ Ocneanu, ``Paths on Coxeter diagrams: From Platonic 
solids and singularities to minimal models and subfactors'' 
(Lectures given at Fields Institute (1995), notes recorded by S.\ Goto).

\item{49.} L.\ Queen, {\it Math.\ of Comput.}\ {\bf 37} (1981) 547.

\item{50.} {A.N.\ Schellekens and S.\  Yankielowicz}, 
{\it Nucl.\ Phys.}\ {\bf B327} (1989) 673.

\item{51.} M.\ Schottenloher, {\it A Mathematical Introduction to Conformal
Field Theory} (Springer, Berlin, 1997).

\item{52.} G.\ Segal, In: {\it IXth Proc.\ Int.\ Congress Math.\ Phys.}\
(Hilger, 1989).

\item{53.} S.\ Singh, {\it Fermat's Enigma} (Penguin Books, London, 1997).

\item{54.} S.\ D.\ Smith, In: {\it Finite Groups -- Coming of Age},
Contemp.\ Math.\ {\bf 193} (AMS, Providence, 1996).


\item{55.} W.\ Thurston, {\it Three-dimensional geometry and topology}, vol.\ 1
(Princeton, 1997).

\item{56.} L.\ Toti Rigatelli, {\it Evariste Galois} (Birkh\"auser, Basel, 1996).

\item{57.} M.\ Tuite, {\it Commun.\ Math.\ Phys.}\  {\bf 166} (1995) 495.

\item{58.} V.\ G.\ Turaev, {\it Quantum invariants of knots and 3-manifolds}
(de Gruyter, Berlin, 1994).

\item{59.} M.\ A.\ Walton, ``Affine Kac-Moody algebras and the
Wess-Zumino-Witten model'', lecture notes in this volume. 

\item{60.} M.\ Waldschmidt, P.\ Moussa, J.-M.\ Luck, and C.\
Itzykson (ed.), {\it From Number Theory to Physics} (Berlin, Springer, 1992).

\item{61.} D.\ Zagier, {\it Amer.\ Math.\ Monthly} {\bf 97} (1990) 144.

\item{62.} J.-B.\ Zuber, {\it Commun.\ Math.\ Phys.}\ {\bf 179} (1996)
265.

\item{63.} J.-B.\ Zuber, private communication, February 1999.

\end